\documentclass[a4paper]{amsart}
\usepackage{enumerate}
\usepackage{amssymb,bm}
\usepackage{graphicx, subfig}
\usepackage{tikz}
\usetikzlibrary{decorations.pathreplacing}
\definecolor{Bittersweet}{HTML}{FE6F5E}
\definecolor{garnet}{rgb}{0.545098,0.133333,0.321569}
\newcommand{\leftlimitsRH}[3]{\kern-0.45em\bgroup\array[t]{l}
      \displaystyle #1 \\[#2]
      \scriptstyle\subarray{l}#3\endsubarray
   \endarray%
\egroup\kern-0.45em}
\newcommand{\LSleftlimits}[2]{\leftlimitsRH{#1}{0pt}{#2}}
\newcommand{\leftlimits}[2]{\leftlimitsRH{#1}{-1ex}{#2}}
\newcommand{\inidemopartfree}[1]{\par\par\medskip\noindent\textbf{#1}\ }
\newcommand{\inidemopart}[1]{\inidemopartfree{(#1)}}
\newtheorem{MainTheorem}{Theorem}

\newtheorem{theorem}{Theorem}[section]
\newtheorem{proposition}[theorem]{Proposition}
\newtheorem{lemma}[theorem]{Lemma}

\newtheorem*{ST}{Sharkovski\u{\i} Theorem for maps from $\cSO$}

\theoremstyle{definition}
\newtheorem{definition}[theorem]{Definition}
\newtheorem{remark}[theorem]{Remark}
\newtheorem{example}[theorem]{Example}
\newcommand{\openthickbox}{\leavevmode\hbox to.77778em{%
  \hfil\vrule width.175em
  \vbox to.675em{\hrule width.26em height.175em\vfil\hrule height.175em}%
  \vrule width.175em\hfil}}
\newcommand{\qedple}{\leavevmode\unskip\penalty9999\hbox{}\nobreak\hfill\quad\hbox{\openthickbox}}
\makeatletter
\newcommand{\TeoremaAmbFinalMarcat}[1]{%
  \expandafter\gdef\csname end#1\endcsname{\qedple\@endtheorem}}
\TeoremaAmbFinalMarcat{definition}
\TeoremaAmbFinalMarcat{remark}
\TeoremaAmbFinalMarcat{example}
\newenvironment{autocase}[2][Case]{\case{#1 #2.}}{\endcase}
\newenvironment{case}[1]{%
  \trivlist \item[\hskip\labelsep{\bfseries #1}]\begin{em}}{
  \end{em}\endtrivlist}
\mathchardef\textupuparrows="\hexnumber@\symAMSa14
\def\@map#1#2[#3]{\mbox{$#1 \colon #2 \longrightarrow #3$}}
\def\map#1#2{\@ifnextchar [{\@map{#1}{#2}}{\@map{#1}{#2}[#2]}}
\def\@ls#1#2{\@mathmeasure\z@\displaystyle{#2}
    \@mathmeasure\@ne\displaystyle{#1}\box\@ne\box\z@}
\def\gtso#1{\mathrel{\@ls{_{#1}}{>}}}
\def\geso#1{\mathrel{\@ls{_{#1}}{\ge}}}
\def\ltso#1{\mathrel{<_{#1}}}
\def\leso#1{\mathrel{\le_{#1}}}
\newenvironment{labeledlist}[1]{\begin{list}{}{\def\makelabel##1{##1:\hfill}%
         \setlength\labelsep{0.5em}\rightmargin\z@\itemindent\z@\leftmargin\labelsep%
         \setbox\z@\hbox{\makelabel{#1}}\labelwidth\wd\z@\advance\leftmargin by \labelwidth%
         \itemsep=2pt\parsep=0pt\topsep=3pt plus 1pt minus 1 pt}}{\end{list}}
\makeatother
\def\Sho{\mbox{\tiny\textup{Sh}}}
\newcommand{\cball}[2]{\ensuremath{B_{#2}\mathchoice{}{\kern-0.199996em}{}{}\left[#1\right]}}
\newcommand{\ball}[2]{\ensuremath{B_{#2}\mathchoice{}{\kern-0.199996em}{}{}\left(#1\right)}}
\def\wings{\raise 0.249995em\hbox{\usebox{\wingspic}}}
\newsavebox{\wingspic}
\sbox{\wingspic}{\begin{tikzpicture}[x=0.1em,y=0.1em]
       \filldraw (0,0.7) .. controls (1, 1) and (1.97, 0.65) .. (2.4,0)
                         .. controls (2.97, 0.65) and (3.9, 1) .. (4.9, 0.7)
                         .. controls (3.9, 1.5) and (2.97, 1) .. (2.4, 0.79)
                         .. controls (1.97, 1) and (1, 1.5) .. (0,0.7);
\end{tikzpicture}}
\newcommand{\wingcball}[2]{\ensuremath{B^{\wings}_{#2}\mathchoice{}{\kern-0.199996em}{}{}\left[#1\right]}}
\newcommand{\wingball}[2]{\ensuremath{B^{\wings}_{#2}\mathchoice{}{\kern-0.199996em}{}{}\left(#1\right)}}
\usepackage{ifthen}
\makeatletter
\newcommand{\AlternativeIfOneCharacter}[3]{\ifthenelse{\equal{\@secondoftwo#1XYZ}{XYZ}}{#2}{#3}}
\makeatother
\newcommand{\sstarplain}[1]{\ensuremath{#1^{*}}}
\newcommand{\sstar}[1]{\AlternativeIfOneCharacter{#1}{\sstarplain{#1}}{\sstarplain{(#1)}}}
\newcommand{\BSG}[3][\alpha]{\cball{\sstar{#2}}{#1_{#3}}}
\newcommand{\OBG}[3][\alpha]{\ball{\sstar{#2}}{#1_{#3}}}
\newcommand{\basint}[2][\alpha]{\BSG[#1]{#2}{#2}}
\newcommand{\basintabs}[2][\alpha]{\BSG[#1]{#2}{\abs{#2}}}
\newcommand{\basintneg}[2][\alpha]{\AlternativeIfOneCharacter{#2}{\cball{\sstarplain{(-#2)}}{#1_{#2}}}{\cball{\sstarplain{(-(#2))}}{#1_{#2}}}}
\newcommand{\obasint}[2][\alpha]{\OBG[#1]{#2}{#2}}
\newcommand{\obasintabs}[2][\alpha]{\OBG[#1]{#2}{\abs{#2}}}
\newcommand{\wbasint}[1]{\wingcball{\sstar{#1}}{#1}}
\newcommand{\wobasint}[1]{\wingball{\sstar{#1}}{#1}}
\newcommand{\wbasband}[1]{\ensuremath{\mathsf{V}^{\kern1pt\wings}_{\AlternativeIfOneCharacter{#1}{\sstar{#1}}{\sstar{(#1)}}}}}
\newcommand{\sstarset}[1]{\ensuremath{\{\sstar{#1}\}}}
\newcommand{\Zstar}{\sstarplain{Z}}
\newcommand{\jstar}{\sstarplain{j}}
\newcommand{\lstar}{\sstarplain{\ell}}
\newcommand{\lstarset}{\sstarset{\ell}}
\newcommand{\all}{\abs{\ell}}
\newcommand{\istar}{\sstarplain{i}}
\newcommand{\istarset}{\sstarset{i}}
\newcommand{\iistarset}{\sstarset{i+1}}
\newcommand{\ai}{\abs{i}}
\newcommand{\aii}{{\abs{i+1}}}
\newcommand{\kstar}{\sstarplain{k}}
\newcommand{\kstarset}{\sstarset{k}}
\newcommand{\mstar}{\sstarplain{m}}
\newcommand{\mstarset}{\sstarset{m}}
\newcommand{\qstarset}{\sstarset{q}}
\newcommand{\ak}{\abs{k}}
\newcommand{\aq}{\abs{q}}
\newcommand{\akk}{\abs{k+1}}
\newcommand{\am}{\abs{m}}

\newcommand{\gams}[1]{\ensuremath{\gamma_{_{#1}}}}
\newcommand{\cR}{\ensuremath{\mathcal{R}}}
\newcommand{\win}[1]{\ensuremath{\cR(#1)}}
\newcommand{\basicbox}[1]{\win{\sstar{#1}}}
\newcommand{\wwin}[1]{\ensuremath{\cR^{\wings}(#1)}}
\newcommand{\wbasicbox}[1]{\wwin{\sstar{#1}}}
\DeclareMathOperator{\dep}{\mathsf{depth}}
\newcommand{\bt}[2][m]{\ensuremath{\mathsf{b}^{\mathchoice{\wings}{\kern-2pt\wings}{\kern-0.1em\wings}{\kern-0.15em\raise 0.23em\hbox{\scalebox{0.7}{\usebox{\wingspic}}}}}\mathchoice{\negmedspace}{\!}{\kern-1pt}{\kern-1.5pt}\left(#2, #1\right)}}
\newcommand{\led}[2][m]{\ensuremath{\mathsf{led}\mathchoice{\!}{}{}{}\left(#2, #1\right)}}
\newcommand{\DS}[1][m]{\ensuremath{\mathfrak{D}_{_{#1}}}}
\newcommand{\IndSetWWings}[2]{\ensuremath{\mathbb{#1}^{\kern1pt\wings}_{_{#2}}}}
\newcommand{\IBD}[1][m]{\ensuremath{\mathbb{B}_{_{#1}}}}

\newcommand{\WDB}[1][m]{\ensuremath{\mathbb{WDB}_{_{#1}}}}
\newcommand{\IW}[2][m]{\ensuremath{\mathbb{IW}_{_{#1, #2}}}}
\newcommand{\wIBD}[1][m]{\IndSetWWings{B}{#1}}
\newcommand{\wIVD}[1][m]{\IndSetWWings{V}{#1}}
\newcommand{\wEIBD}[1][m]{\IndSetWWings{EB}{#1}}
\newcommand{\WDS}[1][m]{\ensuremath{\mathfrak{WFD}_{_{#1}}}}
\newcommand{\WIB}[1][m]{\ensuremath{\mathbb{WIB}_{_{#1}}}}
\newcommand{\WB}[1][m]{\IndSetWWings{WB}{#1}}
\def\is{\mathsf{i}}

\newcommand{\N}{\ensuremath{\mathbb{N}}}
\newcommand{\Z}{\ensuremath{\mathbb{Z}}}
\newcommand{\Q}{\ensuremath{\mathbb{Q}}}
\newcommand{\R}{\ensuremath{\mathbb{R}}}
\newcommand{\I}{\ensuremath{\mathbb{I}}}
\newcommand{\SI}{\ensuremath{\mathbb{S}^1}}
\newcommand{\A}{\mathsf{A}}
\newcommand{\pcs}{\mbox{\large $\bm{\mathcal{A}}$}}
\newcommand{\dist}[1]{\ensuremath{\mathsf{d}_{_{#1}}}}
\newcommand{\dom}{\dist{\Omega}}
\newcommand{\dinf}{\dist{\infty}}
\newcommand{\dSI}{\dist{\SI}}
\DeclareMathOperator{\diam}{\mathsf{diam}}
\DeclareMathOperator{\Int}{\mathsf{Int}}
\DeclareMathOperator{\Bd}{\mathsf{Bd}}
\DeclareMathOperator{\Graph}{\mathsf{Graph}}
\newcommand{\cSO}{\ensuremath{\mathcal{S}(\Omega)}}
\newcommand{\C}{\ensuremath{\mathcal{PCG}}}
\newcommand{\SC}{\ensuremath{\mathcal{PC}}}
\providecommand{\abs}[1]{\ensuremath{\left\lvert#1\right\rvert}}
\providecommand{\norm}[1]{\ensuremath{\left\lVert#1\right\rVert}}
\newcommand{\evalat}[1]{\bigr\rvert_{{#1}}}
\let\setminus\backslash
\newcommand{\set}[2]{\ensuremath{\{#1 \,\colon #2\}}}
\newcommand{\tsfR}{\textsf{R}}
\newcommand{\Orbom}{\ensuremath{\sstarplain{O}(\omega)}}
\newcommand{\pc}[1][\varphi,G]{\ensuremath{\A_{_{(#1)}}}}
\newcommand{\andq}[1][and]{\quad\text{#1}\quad}
\def\liftedsetsymbol{%
    \mathchoice{\raise0.45ex\hbox{$\displaystyle\boldsymbol{\textupuparrows}$}\kern-0.15em}
               {\raise0.45ex\hbox{$\boldsymbol{\textupuparrows}$}\kern-0.15em}
               {\raise0.315ex\hbox{$\scriptstyle\boldsymbol{\textupuparrows}$}\kern-0.15em}
               {\raise0.225ex\hbox{$\scriptscriptstyle\boldsymbol{\textupuparrows}$}\kern-0.15em}%
}
\DeclareMathOperator{\setsilift}{\liftedsetsymbol}
\newcommand{\setfib}[2]{\ensuremath{#1^{\setsilift#2}}}
\newcommand{\setfibpt}[2]{\ensuremath{#1^{#2}}}
\newcommand{\setfibth}[1]{\setfibpt{#1}{\theta}}
\newcommand{\setfibbb}[1]{\setfibpt{\basicbox{#1}}{\sstar{#1}}}
\newcommand{\setfibball}[3][\alpha]{\setfib{\basicbox{#2}}{\BSG[#1]{#2}{#3}}}
\title[A skew-product without fixed-curves]{%
     A quasiperiodically forced skew-product on the cylinder without fixed-curves}
\author[Ll. Alsed\`a]{Llu\'{\i}s Alsed\`a}
\address{Departament de Matem\`atiques and Centre de Recerca Matem\`atica
Edifici Cc,
Universitat Aut\`onoma de Barcelona,
08913 Cerdanyola del Vall\`es,
Barcelona,
Spain}
\email{alseda@mat.uab.cat}

\author[F. Ma\~{n}osas]{Francesc Ma\~{n}osas}
\address{Departament de Matem\`atiques,
Edifici Cc,
Universitat Aut\`onoma de Barcelona,
08913 Cerdanyola del Vall\`es,
Barcelona,
Spain}
\email{manyosas@mat.uab.cat}

\author[L. Morales]{Leopoldo Morales}
\address{Departament de Matem\`atiques,
Edifici Cc,
Universitat Aut\`onoma de Barcelona,
08913 Cerdanyola del Vall\`es,
Barcelona,
Spain}
\email{mleo@mat.uab.cat}
\thanks{The authors have been partially supported by MINECO grant
numbers MTM2008-01486, MTM2011-26995-C02-01 and MTM2014-52209-C2-1-P}
\subjclass{Primary: 37C55, 37C70}
\keywords{Quasiperiodically forced systems on the cylinder, invariant strips}
\date{March 3, 2016}

\begin{document}
\begin{abstract}
In \cite{FJJK} the Sharkovski\u{\i} Theorem was extended to periodic
orbits of strips of quasiperiodic skew products in the cylinder.

In this paper we deal with the following natural question that arises
in this setting: \emph{Does Sharkovski\u{\i} Theorem holds when
restricted to curves instead of general strips?}

We answer this question in the negative by constructing a
counterexample: We construct a map having a periodic orbit of period 2
of curves (which is, in fact, the upper and lower circles of the
cylinder) and without any invariant curve.

In particular this shows that there exist quasiperiodic skew products
in the cylinder without invariant curves.
\end{abstract}
\maketitle
\section{Introduction}
We consider the coexistence and implications
between periodic objects of maps on the cylinder
$\Omega = \SI\times \I,$ of the form:
\[
\map{F}{\begin{pmatrix} \theta \\ x\end{pmatrix}}[
       {\begin{pmatrix} R_\omega(\theta)\\ \zeta(\theta,x)\end{pmatrix}}
    ],
\]
where $\SI = \R / \Z$, $\I$ is an interval of the real line,
$R_\omega(\theta) = \theta + \omega \pmod{1}$ with
$\omega \in \R \setminus \Q$
and $\zeta(\theta,x) = \zeta_{\theta}(x)$ is continuous on both variables.
The class of all maps of the above type will be denoted by \cSO.

In this setting a very basic and natural question is the following:
\emph{is it true that any map in the class $\cSO$ has an invariant
curve?}

In \cite{FJJK}, the authors created an appropriate topological
framework that allowed them to obtain the following extension of
the Sharkovski\u{\i} Theorem to the class {\cSO}\footnote{%
As already remarked in \cite{FJJK}, instead of $\SI$ we could take any
compact metric space~$\Theta$ that admits a minimal homeomorphism
$\map{R}{\Theta}$ such that $R^{\ell}$ is minimal
for every $\ell > 1.$
However, for simplicity and clarity we will remain in the class $\cSO.$}.

Let $X$ be a compact metric space.
We recall that a subset $G \subset X$ is \emph{residual} if it
contains the intersection of a countable family of open dense subsets
in $X.$

In what follows, $\map{\pi}{\Omega}[\SI]$ will denote the
standard projection from $\Omega$ to the circle.
Given a set $B \subset \SI,$ for convenience we will use the following
notation:
\[
 \setsilift{B} := \pi^{-1}(B) = B \times \I \subset \Omega
\]
In the particular case when $B = \{\theta\},$ instead of
$\setsilift{\{\theta\}}$ we will simply write $\setsilift{\theta}.$
Also, given $A \subset \Omega,$ we will denote by $\setfib{A}{B}$ the set
\[
 A \cap \setsilift{B} = \set{(\theta,x) \in \Omega}{\theta \in B \text{ and } (\theta,x) \in A}.
\]
In the particular case when $B = \{\theta\},$ instead of
$\setfib{A}{\theta}$ we will simply write $\setfibth{A}.$

Instead of periodic points we use objects that project over the
whole~$\SI,$ called \emph{strips} in \cite[Definition~3.9]{FJJK}.
A set $B \subset \Omega$ such that $\pi(B) = \SI$
(i.e., $B$ projects on the whole $\SI$) will be called a \emph{circular set}.
\begin{definition}
A \emph{strip in $\Omega$} is a compact circular set $B \subset \Omega$
such that $\setfibth{B}$ is a closed interval (perhaps degenerate to
a point) for every $\theta$ in a residual set of $\SI.$
\end{definition}

Given two strips $A$ and $B,$ we will write $A < B$ and $A \le B$
(\cite[Definition~3.13]{FJJK}) if there exists a residual set
$G \subset \SI,$ such that
for every $(\theta,x) \in \setfib{A}{G}$ and
$(\theta,y) \in \setfib{B}{G}$
it follows that $x < y$ and,
respectively, $x \le y$.
We say that the strips $A$ and $B$ are \emph{ordered}
(respectively \emph{weakly ordered})
if either $A < B$ or $A > B$
(respectively $A \le B$ or $A \ge B$).

\begin{figure}[ht]
\begin{center}
\hfill \subfloat[$3.28 x(1-x) + \tfrac{4}{100}\cos(2\pi\theta)$]{\includegraphics[width=0.46\textwidth]{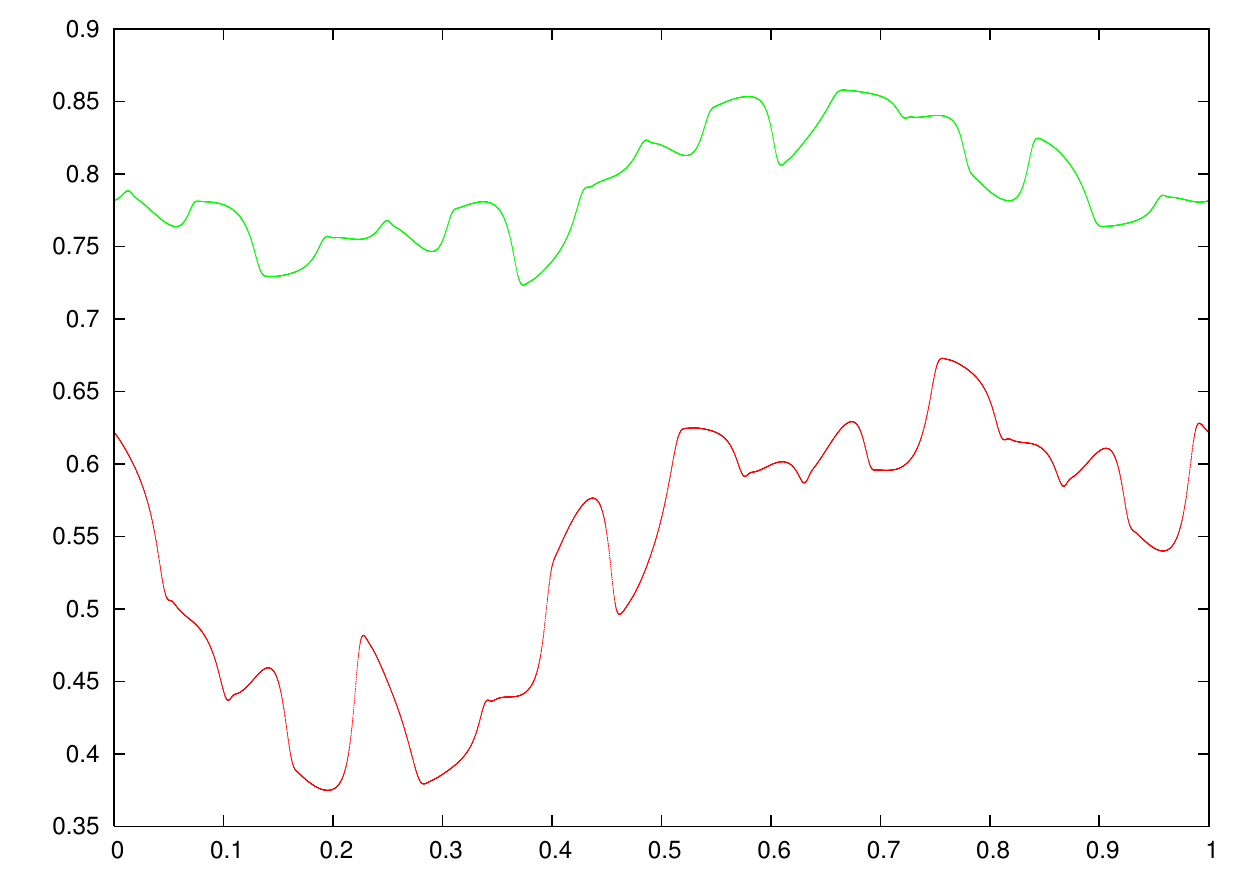}}
\hfill \subfloat[$3.85 x(1-x)(1 + \frac{111}{10^5}\cos(2\pi\theta))$]{\includegraphics[width=0.46\textwidth]{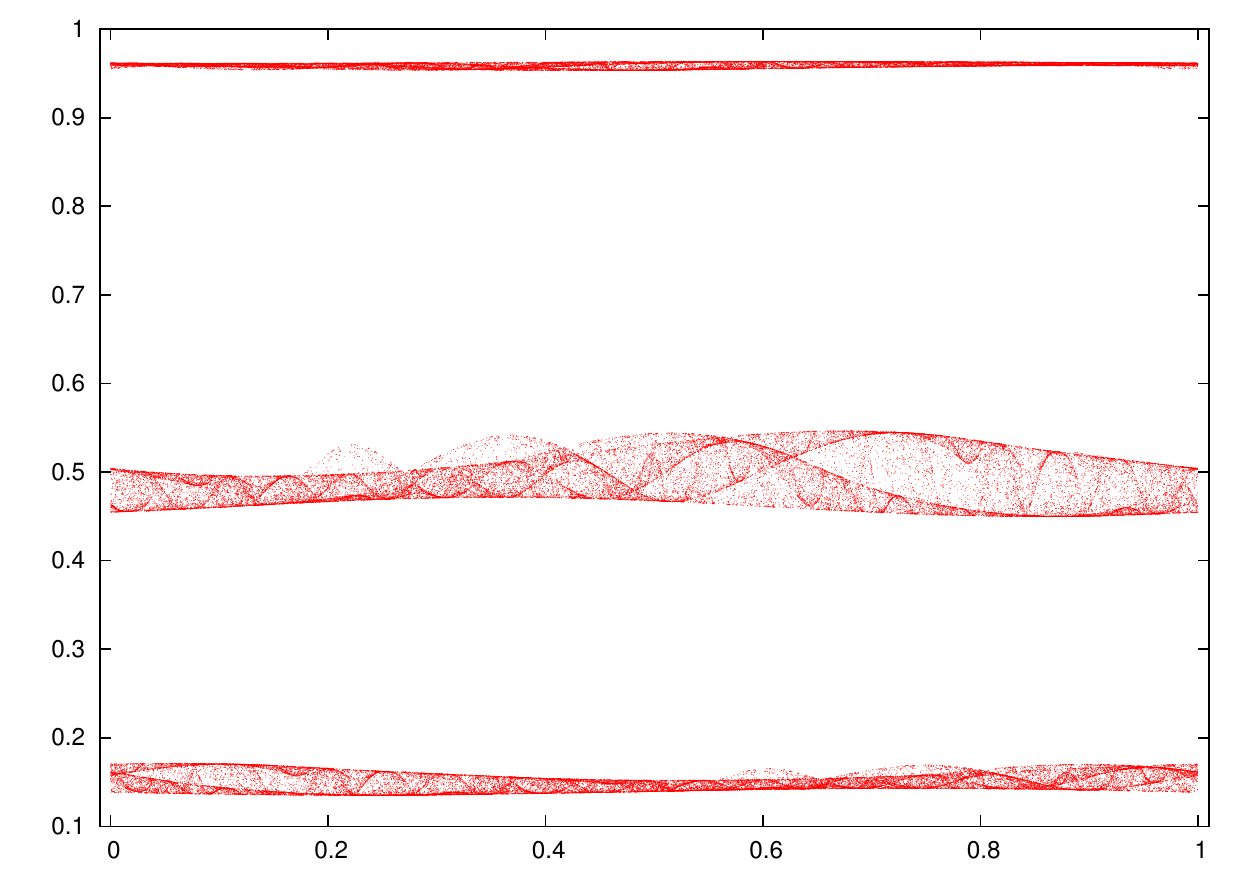}}
\hfill \strut
\end{center}
\caption{In the left picture we show an example two periodic orbit of curves,
and in the second we show a possible example of a three periodic orbit solid strips.}\label{fig-example-periodicorbits}
\end{figure}

\begin{definition}[\protect{\cite[Definition~3.15]{FJJK}}]
A strip $B \subset \Omega$ is called \emph{$n$-periodic} for $F \in \cSO$
if $F^{n}(B) = B$ and the image sets
$B,\ F(B),\ F^{2}(B),\dots, F^{n-1}(B)$
are pairwise disjoint and pairwise ordered
(see Figure~\ref{fig-example-periodicorbits} for examples).
\end{definition}

To state the main theorem of \cite{FJJK} we need to recall the
\emph{Sharkovski\u{\i} Ordering} (\cite{Shar, Shartrans}).
The \emph{Sharkovski\u{\i} Ordering} is a linear ordering of $\N$
defined as follows:
\begin{align*}
& 3 \gtso{\Sho} 5 \gtso{\Sho} 7 \gtso{\Sho} 9 \gtso{\Sho} \dots \gtso{\Sho} \\
& 2 \cdot 3 \gtso{\Sho} 2 \cdot 5 \gtso{\Sho} 2 \cdot 7 \gtso{\Sho} 2 \cdot 9 \gtso{\Sho} \dots \gtso{\Sho} \\
& 4 \cdot 3 \gtso{\Sho} 4 \cdot 5 \gtso{\Sho} 4 \cdot 7 \gtso{\Sho} 4 \cdot 9 \gtso{\Sho} \dots \gtso{\Sho}\\
& \hspace*{7em} \vdots \\
& 2^n \cdot 3  \gtso{\Sho}  2^n \cdot 5 \gtso{\Sho} 2^n \cdot 7 \gtso{\Sho} 2^n \cdot 9 \gtso{\Sho} \dots \gtso{\Sho} \\
& \hspace*{7em} \vdots \\
& \cdots \gtso{\Sho} 2^n \gtso{\Sho} \dots \gtso{\Sho}
  16 \gtso{\Sho} 8 \gtso{\Sho} 4 \gtso{\Sho} 2 \gtso{\Sho} 1.
\end{align*}
In the ordering $\geso{\Sho}$ the least element is 1 and the largest
one is 3. The supremum of the set $\{1,2,4,\dots,2^n,\dots\}$ does not
exist.

\begin{ST}[\cite{FJJK}]\label{SharTheo}
Assume that the map $F \in \cSO$ has a $p$-periodic strip.
Then $F$ has a $q$-periodic strip for every $q \ltso{\Sho} p.$
\end{ST}

In view of this result, the new following natural question
(that is stronger that the previous one) arises:
\emph{Does Theorem~\ref{SharTheo} holds when restricted to curves?}
where a curve is defined as the graph of a continuous map from $\SI$ to $\I$.
More precisely, \emph{is it true that if $F$ has a $q$-periodic
curve and $p\leso{\Sho} q$ then does there exists a $p$-periodic curve of $F$?}

The aim of this paper is to answer both of the above questions in the negative
by constructing a counterexample.
This is done by the following result which is the main result of the paper.

\begin{MainTheorem}\label{MainTh}
There exists a map $T \in \cSO$ with
$f(\theta,\cdot)$ non-increasing for every $\theta \in \SI,$
such that $T$ permutes the upper and lower circles of $\Omega$
(thus having a periodic orbit of period two of curves),
and $T$ does not have any invariant curve.
\end{MainTheorem}

The construction will be done in two steps.
First, in Section~\ref{pseudo-curve}, we construct a strip $A$
which is a pseudo-curve which is not a curve.
This strip is obtained as a \emph{limit} of sets defined inductively
by using of a collection of \emph{winged boxes} $\wbasicbox{i} \subset \Omega.$
Second, we construct a Cauchy sequence $\{T_m\}_{m=0}^{\infty}$ that gives
as a limit the function $T$ from Theorem~\ref{MainTh} having $A$ as invariant set.
To this end, in Section~\ref{FunctionsGi} we define a collection of auxiliary
functions $G_i$ defined on the winged boxes $\wbasicbox{i}$.
Next, in Section~\ref{stratification} we introduce a notion of \emph{depth}
in the set of winged boxes $\wbasicbox{i}$ which defines a convenient
stratification in the set of winged boxes $\wbasicbox{i}.$
In Section~\ref{BoxesintheWings} we study the wings of box
and its interaction with boxes of higher depth.
In Section~\ref{skew-product},
by using the auxiliary functions from Section~\ref{FunctionsGi},
the stratification from Section~\ref{stratification} and the technical results
from Section~\ref{BoxesintheWings}
we construct the Cauchy sequence $\{T_m\}_{m=0}^{\infty} \subset \cSO,$
we define the map $T = \lim_{m\to\infty} T_m,$
and we prove Theorem~\ref{MainTh}.

For clarity, we omit the proofs of all results from Section~\ref{skew-product}.
These proofs will be provided in Sections~\ref{proofofT0mapProperties},
\ref{proofofseqTmProperties} and \ref{proofofdistTmTm-1}.
Section~\ref{definitions} is devoted to introduce the necessary
definitions and, in particular, to introduce the notion of
pseudo-curve and some necessary results on the space of pseudo-curves.

\section{Definitions and preliminary results}\label{definitions}

The main aim of this section is to introduce the definition and
basic results about pseudo-curves.

Given $G \subset \SI$ and a map $\map{\varphi}{G}[\I]$,
$\Graph(\varphi)$ denotes the \emph{graph of $\varphi$}.
Also, given a set $A$ we will denote the closure of $A$
by $\overline{A}$.

\begin{definition}[Pseudo-curve]\label{PCDefinition}
Let $G$ be a residual set of $\SI$ and let $\map{\varphi}{G}[\I]$
be a continuous map from $G$ to $\I.$
The set $\overline{\Graph(\varphi)},$ denoted by $\pc,$
will be called a \emph{pseudo-curve}.
Notice that every pseudo-curve is a compact circular set.

Also, $\pcs$ will denote the class of all pseudo-curves.
\end{definition}

A set $A \subset \Omega$ is \emph{$F$-invariant}
(respectively \emph{strongly $F$-invariant})
if $F(A) \subset A$ (respectively $F(A) = A$).
Observe that if $F \in \cSO,$ every compact
$F$-invariant set is circular.
A closed invariant set is called \emph{minimal}
if it does not contain any proper closed invariant set.

An \emph{arc of a curve} is the graph of a continuous
function from an arc of $\SI$ to $\I$.

The pseudo-curves have the following properties which are easy to prove:
\begin{lemma}\label{PC-properties}
Given a pseudo-curve $\pc \in \pcs$ the following statements hold.
\begin{enumerate}[(a)]
\item $\setfibth{\pc}$ consists of a single point for  every $\theta \in G.$
      Consequently, \[ \setfib{\pc}{G} = \Graph(\varphi). \]
\item Every circular compact set contained in a pseudo-curve
      coincides with the pseudo-curve.
\item $\pc = \overline{\Graph(\varphi\evalat{\widetilde{G}})}$
      for every $\widetilde{G} \subset G$ dense in $\SI.$
\item If $\pc$ contains a curve then it is a curve.
\end{enumerate}
\end{lemma}

\begin{proof}
We start by proving (a).
By the definition of a pseudo-curve we have
$\Graph(\varphi) \subset \setfib{\pc}{G}.$
To prove the other inclusion fix $\theta \in G$ and $x \in \I$
such that $(\theta,x) \in \pc.$
Then, there exists a sequence
$
\{(\theta_n,\varphi(\theta_n))\}_{n=1}^{\infty} \subset \Graph(\varphi)
$
such that
$\lim_{n\to\infty}(\theta_n,\varphi(\theta_n)) = (\theta,x).$
The continuity of $\varphi$ in $G$ (and hence in $\theta$)
implies $x = \varphi(\theta)$ and, therefore,
$(\theta,x)\in \Graph(\varphi).$

Now we prove (b).
Assume that $B \subset \pc$ is a circular compact set.
From the assumptions and statement (a) we get
$\setfib{\pc}{G} = \setfib{B}{G}.$
Hence,
\[
 \pc = \overline{\Graph(\varphi)} = \overline{\setfib{\pc}{G}} =
       \overline{\setfib{B}{G}} \subset B.
\]

Now (d) follows directly from (b) and the fact that a curve
is compact since it is the graph of a continuous function.
Statement (c) also follows from (b) because
$
\overline{\Graph\left(\varphi\evalat{\widetilde{G}}\right)} \subset \pc
$
and
$
\overline{\Graph\left(\varphi\evalat{\widetilde{G}}\right)}
$
is a circular set (since $\widetilde{G}$ is dense in $\SI$).
\end{proof}

We also will be interested in the pseudo-curves as a possible invariant
objects of maps from $\cSO.$ The next lemma studies their properties
in this case.

\begin{lemma}\label{PC-properties-invariant}
Let $F \in \cSO$ and assume that $\pc \in \pcs$
is an $F$-invariant pseudo-curve. Then,
\begin{enumerate}[(a)]
\item $\pc$ is strongly $F$-invariant and minimal.
\item If $\pc$ contains an arc of a curve then it is a curve.
\end{enumerate}
\end{lemma}

\begin{proof}
We start by proving (a).
Let $B \subset \pc$ be a closed invariant set.
We have that $B$ is circular and, by Lemma~\ref{PC-properties}(b),
$B = \pc.$ Hence, $\pc$ is minimal.

On the other hand, $F(\pc) \subset \pc$ implies
$F^{2}(\pc) \subset F(\pc)$ and, hence,
$F(\pc)$ is a compact $F$-invariant set.
Therefore, by the part already proven,
$F(\pc) = \pc.$

Now we prove (b).
Let $S$ be an (open) arc of $\SI$ and let
{\map{\xi}{S}[\I]} be a continuous map such that
$\Graph(\xi) \subset \pc.$
Clearly, there exists $m\in \N$ such that
$\bigcup_{i=0}^m R^i_\omega(S) = \SI.$
Now we set $\xi_0 := \xi$ and, for $i=1,2,\dots,m,$
we define {\map{\xi_i}{R^i_\omega(S)}[\I]} by
\[
\xi_i(\theta) :=
  f\left(R^{-1}_\omega(\theta),\xi_{i-1}\left(R^{-1}_\omega(\theta)\right)\right).
\]
The continuity of $f$ implies that every $\xi_i$
is an arc of a curve and
$\Graph(\xi_i) = F(\Graph(\xi_{i-1})).$
Hence,
\[
\bigcup_{i=0}^m \Graph(\xi_i) = \bigcup_{i=0}^m F^i(\Graph(\xi)) \subset \pc
\]
because $\pc$ is $F$-invariant.

In view of Lemma~\ref{PC-properties}(d) we only have to show that
$\bigcup_{i=0}^m \Graph(\xi_i)$
is a curve. We will prove prove this by induction.

Assume that $\emptyset \ne M \varsubsetneq \{0,1,2,\dots,m\}$
verifies that
$S_M:= \bigcup_{i \in M} R^i_\omega(S)$
is an (open) arc of $\SI$ and
$\bigcup_{i \in M} \Graph(\xi_i)$ is an arc of a curve
(initially we can take $M$ to be any unitary subset of $\{0,1,2,\dots,m\}$).
Then, there exists a continuous map {\map{\xi_{_M}}{S_M}[\I]}
such that $\Graph(\xi_{_M}) = \bigcup_{i \in M} \Graph(\xi_i).$

Clearly, there exists $j \in \{0,1,2,\dots,m\}\setminus M$
such that
$
S_{M,j} := S_M \cap R^j_\omega(S) \ne \emptyset.
$
The set $S_{M,j}$ is an open arc of $\SI$ and, by Lemma~\ref{PC-properties}(a),
$\xi_{_M}\evalat{S_{M,j}\cap G} = \xi_j\evalat{S_{M,j}\cap G}$
because
$\Graph(\xi_{_M}), \Graph(\xi_j) \subset \pc.$
Since $S_{M,j}\cap G$ is dense in $S_{M,j},$
given $\theta \in S_{M,j}\setminus G,$
there exists  a sequence
$\{\theta_n\}_{n=0}^\infty \subset S_{M,j}\cap G$
converging to $\theta.$
The continuity of $\xi_{_M}$ and $\xi_j$ on $S_{M,j}$ implies that
$
 \xi_{_M}(\theta) = \lim_{n\to\infty} \xi_{_M}(\theta_n) =
 \lim_{n\to\infty} \xi_j(\theta_n) = \xi_j(\theta).
$
Consequently, $\xi_{_M}\evalat{S_{M,j}} = \xi_j\evalat{S_{M,j}}$
and $\Graph(\xi_{_M}) \cup \Graph(\xi_j)$ is an arc of a curve
(defined on the open arc $S_M \cup R^j_\omega(S)$).
By redefining $M$ as $M \cup \{j\}$ and iterating this procedure until
$M \cup \{j\} = \{0,1,2,\dots,m\}$ we see that the whole
$\bigcup_{i=0}^m \Graph(\xi_i)$ is a curve.
\end{proof}

Next we will introduce and study the space of pseudo-curves.

\begin{definition}\label{PCG-dinfinito}
We define the \emph{space of pseudo-curve generators} as
\[
\C := \set{(\varphi,G)}{\text{$G$ is a residual set in $\SI$ and {\map{\varphi}{G}[\I]} is a continuous map}}.
\]
On $\C$ we also define the \emph{supremum pseudo-metric}
$\map{\dinf}{\C\times\C}[\R^+]$
by:
\[
  \dinf\bigl((\varphi,G),(\varphi',G')\bigr) :=
     \sup_{\theta\in G\cap G'} \abs{\varphi(\theta)-\varphi'(\theta)}.
\]
Clearly,
$\dinf((\varphi,G),(\varphi',G')) = 0$
if and only if $\varphi\evalat{G\cap G'} = \varphi'\evalat{G\cap G'}$
and, hence, $\dinf$ is a pseudo-metric.
\end{definition}

The next lemma will be useful in using the metric $\dinf.$

\begin{lemma}\label{dinfinito-equiv}
Let $(\varphi,G),(\varphi',G') \in \C.$
Then,
\[
 \dinf\bigl((\varphi,G),(\varphi',G')\bigr) =
   \sup_{\theta\in \widetilde{G}} \abs{\varphi(\theta)-\varphi'(\theta)}
\]
for every $\widetilde{G} \subset G \cap G'$ dense in $\SI.$
\end{lemma}

\begin{proof}
Set
$
\dist{\infty, \widetilde{G}}\bigl((\varphi,G),(\varphi',G')\bigr) :=
   \sup_{\theta\in \widetilde{G}} \abs{\varphi(\theta)-\varphi'(\theta)}.
$
With this notation, we clearly have
$
\dist{\infty, \widetilde{G}}\bigl((\varphi,G),(\varphi',G')\bigr) \le
\dinf\bigl((\varphi,G),(\varphi',G')\bigr).
$

To prove the reverse inequality take
$\theta \in (G \cap G')\setminus \widetilde{G}.$
Since $\widetilde{G}$ is dense in $\SI,$
there exists  a sequence
$\{\theta_n\}_{n=0}^\infty \subset \widetilde{G}$
converging to $\theta.$
On the other hand, by definition, the maps $\varphi$ and $\varphi',$
are continuous in $G \cap G'$ (and, hence, in $\theta$).
Consequently,
$
\abs{\varphi(\theta),\varphi'(\theta)} =
  \lim_{n\to\infty} \abs{\varphi(\theta_n)-\varphi'(\theta_n)} \le
  \dist{\infty, \widetilde{G}}\bigl((\varphi,G),(\varphi',G')\bigr).
$
This ends the proof of the lemma.
\end{proof}

As it is customary we will introduce an equivalent relation
in the space of pseudo-curve generators so that the quotient space
will be a metric space.

\begin{definition}\label{relation}
Two pseudo-curve generators
$(\varphi,G),(\varphi',G')\in \C$ are said to be equivalent,
denoted by $(\varphi,G) \sim (\varphi',G')$
if and only if $\pc = \pc[\varphi',G'].$
Clearly $\sim$ is an equivalence relation in $\C$.
The $\sim$-equivalence class of $(\varphi,G)\in \C$
will be denoted by $[\varphi,G].$
\end{definition}

\begin{remark}\label{equiv-relation}
From Lemma~\ref{PC-properties}(a,c) it follows that
$(\varphi,G) \sim (\varphi',G')$ if and only if
$\varphi\evalat{\widetilde{G}} = \varphi'\evalat{\widetilde{G}}$
for every $\widetilde{G} \subset G \cap G'$ dense in $\SI.$
In particular, by taking $\widetilde{G} = G \cap G',$
we get that $\dinf((\varphi,G),(\varphi',G')) = 0$
if and only if $(\varphi,G) \sim (\varphi',G').$
\end{remark}

\begin{definition}
The space $\C/\kern-3pt\sim$ will be called the
\emph{space of pseudo-curves generator classes}
and denoted by $\SC.$
Also, on $\SC$ we define the \emph{supremum metric}, also denoted
$\map{\dinf}{\SC\times\SC}[\R^+]$
by abuse of notation,
in the following way.
Given $A = [\varphi_A,G_A], B=[\varphi_B,G_B] \in \SC$ we set
\[
 \dinf(A,B) :=  \dinf\bigl((\varphi_A,G_A),(\varphi_B,G_B)\bigr).
\]
Note that $\dinf$  is well defined. To see this take
$[\varphi_A,G_A]=[\varphi_{'A},G_{A'}], [\varphi_B,G_B] \in \C.$
Then, by Lemma~\ref{dinfinito-equiv} and Remark~\ref{equiv-relation}
applied to $\widetilde{G} = G_A \cap G_{A'} \cap G_B$
we get
$
\dinf\bigl((\varphi_A,G_A),(\varphi_B,G_B)\bigr) =
\dinf\bigl((\varphi_{A'},G_{A'}),(\varphi_B,G_B)\bigr).
$
\end{definition}

The next result establishes the basic properties of the
space of pseudo-curves generator classes $(\SC, \dinf).$

\begin{proposition}\label{Ccompleto}
The space of pseudo-curves generator classes $\SC$ is a complete metric space.
\end{proposition}

\begin{proof}
The fact that $\dinf$ is a metric in $\SC$
follows from Remark~\ref{equiv-relation}.

Now we prove that $\SC$ is complete.
Assume that $\{[\varphi_n,G_n]\}_{n=1}^{\infty}$
is a Cauchy sequence in $\SC.$
We have to see that $\lim_{n\to\infty} [\varphi_n,G_n] \in \SC.$

Set, $G := \cap_{i=1}^{\infty} G_n.$
Since this intersection is countable, $G$ is still a residual set.
The definition of $\dinf$ implies that the sequence
$\{\varphi_n(\theta)\}_{n=1}^{\infty} \subset \I$
is a Cauchy sequence in $\I$ for every $\theta\in G.$
So, it is convergent and we can define a map $\map{\varphi}{G}[\I]$
by $\varphi(\theta) := \lim_{n\to\infty} \varphi_n(\theta).$

If $(\varphi,G)\in \C$ we have $[\varphi,G] \in \SC$ and,
from the definition of $\varphi$ it follows that
\[
\lim_{n\to\infty} \dinf([\varphi,G], [\varphi_n,G_n]) =
\leftlimits{\sup}{\theta\in G\cap G_n} \lim_{n\to\infty}  \abs{\varphi(\theta)-\varphi_n(\theta)} = 0.
\]
Consequently, $[\varphi,G] = \lim_{n\to\infty} [\varphi_n,G_n].$
Since $\varphi$ is the uniform limit of a sequence of continuous functions on $G,$
it is continuous on $G.$ That is, $(\varphi,G)\in \C.$
\end{proof}

In what follows we want to look at the space $\pcs$ as a metric space and
relate this metric space with $(\SC, \dinf).$

Let $\rho$ denote the euclidean metric in $\Omega.$
Then, the space $(\Omega, \rho)$ is a compact metric space.
We recall that the \emph{Hausdorff metric} is defined
in the space of compact subsets of $(\Omega, \rho),$ by
\[
H_{\rho}(\A,\mathsf{B}) = \max\left\{%
         \leftlimits{\max}{(\theta,x) \in \A} \rho((\theta,x),\mathsf{B}),
         \leftlimits{\max}{(\theta,x) \in \mathsf{B}} \rho((\theta,x),\A)
\right\}.
\]
Then, $(\pcs, H_{\rho})$ is a metric space. To study the relation between
$(\SC, \dinf)$ and $(\pcs, H_{\rho})$ we need a couple of simple technical
results.

\begin{lemma}\label{HdimCCS}
Let $\A,\mathsf{B} \subset \Omega$ be compact circular sets.
Then,
\[
H_{\rho}(\A,\mathsf{B}) \le
  \max_{\theta \in \SI} H_{\rho}\bigr(\setfibth{\A},\setfibth{\mathsf{B}}\bigl).
\]
\end{lemma}

\begin{proof}
It follows directly from the definitions:
\begin{align*}
H_{\rho}\left(\A, \mathsf{B}\right)
   &\le \max\left\{
           \leftlimits{\sup}{(\theta,x) \in \A} \rho\Bigl((\theta,x), \setfibth{\mathsf{B}}\Bigr),
           \leftlimits{\sup}{(\theta,x) \in \mathsf{B}} \rho\Bigl((\theta,x), \setfibth{\A}\Bigr)
        \right\}\\
  &= \max\left\{
           \sup_{\theta \in \SI}\leftlimits{\max}{\set{x\in \I}{(\theta,x) \in \A}} \rho\Bigl((\theta,x), \setfibth{\mathsf{B}}\Bigr),\right.\\
  &\hspace*{15.5em}         \left.\sup_{\theta \in \SI}\leftlimits{\max}{\set{x\in \I}{(\theta,x) \in \mathsf{B}}} \rho\Bigl((\theta,x), \setfibth{\A}\Bigr)
        \right\}\\
  &= \sup_{\theta \in \SI} \max\left\{
           \leftlimits{\max}{\set{x\in \I}{(\theta,x) \in \A}} \rho\Bigl((\theta,x), \setfibth{\mathsf{B}}\Bigr),
           \leftlimits{\max}{\set{x\in \I}{(\theta,x) \in \mathsf{B}}} \rho\Bigl((\theta,x), \setfibth{\A}\Bigr)
        \right\}\\
  &= \sup_{\theta \in \SI} H_{\rho}\left(\setfibth{\A}, \setfibth{\mathsf{B}}\right).
\end{align*}
\end{proof}

\begin{proposition}\label{Hdimdinfty}
Let $(\varphi,G), (\widetilde{\varphi},\widetilde{G}) \in \C.$
Then,
\[
 H_{\rho}\left(\pc, \pc[\widetilde{\varphi},\widetilde{G}]\right) \le
 \sup_{\theta \in \SI} H_{\rho}\left(
           \setfibth{\pc},
           \setfibth{\pc[\widetilde{\varphi},\widetilde{G}]}
 \right) =
 \dinf\bigl((\varphi,G),(\widetilde{\varphi},\widetilde{G})\bigr).
\]
\end{proposition}

\begin{proof}
The first inequality follows from Lemma~\ref{HdimCCS}.

Now we prove the second equality.
By Lemma~\ref{PC-properties}(a),
\[
\dinf\bigl((\varphi,G),(\widetilde{\varphi},\widetilde{G})\bigr)
  = \leftlimits{\sup}{\theta \in G \cap \widetilde{G}} \abs{\varphi(\theta) - \widetilde{\varphi}(\theta)}
  = \leftlimits{\sup}{\theta \in G \cap \widetilde{G}} H_{\rho}\left(
             \setfibth{\pc},
             \setfibth{\pc[\widetilde{\varphi},\widetilde{G}]}
    \right).
\]
So, to end the proof of the lemma, we have to see that
\[
    H_{\rho}\left(
             \setfibth{\pc},
             \setfibth{\pc[\widetilde{\varphi},\widetilde{G}]}
    \right) \le
 \dinf\bigl((\varphi,G),(\widetilde{\varphi},\widetilde{G})\bigr)
 \andq[for every]
 \theta \in \SI \setminus (G \cap \widetilde{G}).
\]
Fix $\theta \in \SI \setminus (G \cap \widetilde{G}).$
From the definition of the Hausdorff metric it follows that there exist
$x,y \in \I$ such that
$
H_{\rho}\left(
    \setfibth{\pc},
    \setfibth{\pc[\widetilde{\varphi},\widetilde{G}]}
\right) = \abs{x-y},
$
$(\theta,x) \in \setfibth{\pc},$ and
$(\theta,y) \in \setfibth{\pc[\widetilde{\varphi},\widetilde{G}]}.$

Since $G \cap \widetilde{G}$ is residual (and thus dense) in $\SI,$
from Lemma~\ref{PC-properties}(a,c) it follows that there exists sequences
$
\{(\theta_n,\varphi(\theta_n))\}_{n=0}^\infty,\
\{(\theta_n,\widetilde{\varphi}(\theta_n))\}_{n=0}^\infty \subset \setsilift{(G \cap \widetilde{G})}
$
such that
$\lim_{n\to\infty} (\theta_n,\varphi(\theta_n)) = (\theta,x)$ and
$\lim_{n\to\infty} (\theta_n,\widetilde{\varphi}(\theta_n)) = (\theta,y).$
Hence,
\[
H_{\rho}\left(
    \setfibth{\pc},
    \setfibth{\pc[\widetilde{\varphi},\widetilde{G}]}
\right)
  = \abs{x-y}
  = \lim_{n\to\infty} \abs{\varphi(\theta_n) - \widetilde{\varphi}(\theta_n)}
  \le \dinf\bigl((\varphi,G),(\widetilde{\varphi},\widetilde{G})\bigr).
\]
\end{proof}

Proposition~\ref{Hdimdinfty} tells us that that if
$\{[\varphi_n,G_n]\}_{n=1}^{\infty}$ is a Cauchy sequence in $\SC$
then $\pc[\varphi_n,G_n]$ is a Cauchy sequence in $(\pcs, H_{\rho}),$
and if $[\varphi,G] = \lim_{n\to\infty} [\varphi_n,G_n]$ then
$\pc = \lim_{n\to\infty} \pc[\varphi_n,G_n].$
Unfortunately the space $(\pcs, H_{\rho})$ is not complete as the
following simple example shows.

\begin{example}[The space $(\pcs, H_{\rho})$ is not complete]
Consider continuous maps {\map{\xi_n}{\SI}[\I]} with $n\in \N,\ n \ge 2,$
defined by
\[
\xi_n(\theta) = \begin{cases}
   2n\theta     & \text{if $\theta \in [0,\tfrac{1}{2n}]$,}\\
   2(1-n\theta) & \text{if $\theta \in [\tfrac{1}{2n}, \tfrac{1}{n}]$,}\\
   0            & \text{if $\theta \ge \tfrac{1}{n}$.}
\end{cases}
\]
Clearly, $(\xi_n,\SI) \in \C$ and
$
H_{\rho}(\pc[\xi_n,\SI], \pc[\xi_m,\SI]) \le \tfrac{1}{\min\{n,m\}}.
$
Hence, $\{\pc[\xi_n,\SI]\}$ is a Cauchy sequence in $\pcs.$
However, the sequence $\{\pc[\xi_n,\SI]\}$ has no limit in $\pcs.$
Indeed,
$
  \lim_{n\to\infty} \pc[\xi_n,\SI] = \mathsf{L} =
  (\SI \times \{0\}) \cup (\{0\} \times [0,1]),
$
which is not the closure of the graph of a continuous map on a
residual set of $\SI$ (in other words, $\mathsf{L} \notin \pcs$).
This is consistent with the fact that, clearly, $\{[\xi_n,\SI]\}$
is not a Cauchy sequence in $(\SC, \dinf).$
\end{example}

\section{Construction of a connected pseudo-curve}\label{pseudo-curve}
The aim of this subsection is to construct a strip $\A = \pc[\gamma,G]$
as a connected pseudo-curve with certain topological
properties that will allow us to define the map $T \in \cSO$ having this
pseudo-curve as the only proper invariant object.
The pseudo-curve $\pc[\gamma,G]$ will be obtained as a limit in $\SC$
of a sequence of pseudo-curves that will be constructed recursively.

We will start by introducing the necessary notation.

In what follows, for simplicity, we will take the interval $\I$ as the
interval $[-2,2].$ Also, fix $\omega \in [0,1]\setminus\Q.$
For any $\ell \in \Z$ set
$\lstar = \ell\omega \pmod{1}$ and
$\Orbom = \set{\lstar}{\ell\in\Z}.$
That is, $\Orbom$ is the orbit of $0$ by the rotation of angle $\omega.$

We will denote by $\dSI$ the arc distance on $\SI = \R/\Z$.
That is, for $\theta_1, \theta_2 \in \SI,$ we set
\[
\dSI(\theta_1, \theta_2) := \begin{cases}
    \theta_2 - \theta_1 & \text{when $\theta_1 \le \theta_2$, and}\\
    (\theta_2 + 1) - \theta_1 & \text{when $\theta_1 > \theta_2$.}
\end{cases}
\]
The closed arc of $\SI$ joining $\theta_1$ and $\theta_2$ in the natural direction
will be denoted by $[\theta_1, \theta_2].$
That is,
\[
[\theta_1, \theta_2] = \begin{cases}
   \set{t \pmod{1}}{\theta_1 \le t \le \theta_2}
      & \text{when $\theta_1 \le \theta_2$, and}\\
   \set{t \pmod{1}}{\theta_1 \le t \le \theta_2+1}
      & \text{when $\theta_1 > \theta_2$.}
\end{cases}
\]
The open arc of $\SI$ joining $\theta_1$ and $\theta_2$
will be denoted by
$(\theta_1, \theta_2) = [\theta_1, \theta_2] \setminus \{\theta_1, \theta_2\},$
and is defined analogously with strict inequalities
Given an arc $B \subset \SI$, $\Bd(B)$ will denote the set of endpoints of $B.$

We will denote the open (respectively closed) ball (in $\SI$)
of radius $\delta$ centred at $\theta \in \SI$ by
$\ball{\theta}{\delta}$ (respectively $\cball{\theta}{\delta}$):
\begin{align*}
\ball{\theta}{\delta} &=
    \set{\widetilde{\theta} \in \SI}{\dSI(\theta, \widetilde{\theta}) < \delta}
    = (\theta - \delta \pmod{1}, \theta + \delta \pmod{1}),\text{ and}\\
\cball{\theta}{\delta} &=
     \overline{\ball{\theta}{\delta}} =
     \set{\widetilde{\theta} \in \SI}{\dSI(\theta, \widetilde{\theta}) \le \delta}
     = [\theta - \delta \pmod{1}, \theta + \delta \pmod{1}].
\end{align*}

We consider the space $\Omega$ endowed the metric induced by the maximum of $\dSI$
and the absolute value on $\I$. That is, given $(\theta,x), (\nu, y) \in \Omega$
we set
\[
  \dom((\theta,x), (\nu, y)) := \max\left\{\dSI(\theta,\nu), \abs{x-y}\right\}.
\]
Then, given $A \subset \Omega$ we will denote the
\emph{interior of $A$} by $\Int(A)$ and
$\diam(A)$ will denote the \emph{diameter of $A$} whenever $A$ is compact.

To define the sequence of pseudo-curves that will converge to
$\pc[\gamma,G]$ we first need to construct an auxiliary family
$\{\basicbox{\ell}\}_{\ell\in\Z}$ of compact regions in $\Omega$
and a family of compact sets
$\{\Gamma\varphi_{_{\lstar}}\}_{\ell\in\Z}$
such that, for every $\ell \in \Z,$
$
\Gamma\varphi_{_{\lstar}} \subset \basicbox{\ell}
$
and it is the restriction of a pseudo-curve generator to $\pi(\basicbox{\ell}).$
To do this we define the auxiliary functions
$\map{\beta}{[-1,1]}$
and $\map{\phi}{[-1,1]\setminus \{0\}}[{[-1,1]}]$
by (see Figure~\ref{fig1}):
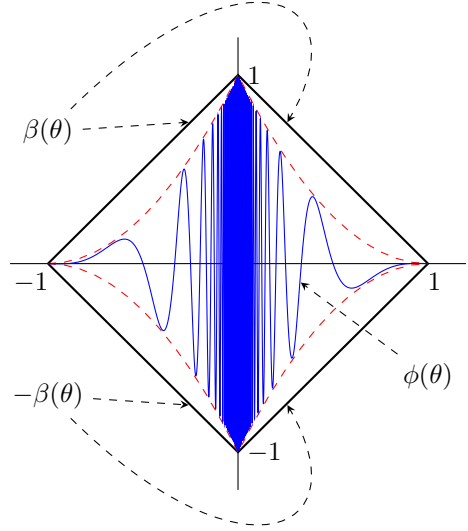
\begin{figure}[tb]
\begin{tikzpicture}[scale=2.5, domain=-1.2:1.2]
\draw[thick] (-1,0) -- (0,1) -- (1,0) -- (0,-1) -- (-1,0); 
\draw (-1.2,0) -- (1.2,0); 
\draw (0,-1.2) -- (0,1.2); 
\node at (0,1) [right] {$1$}; \node at (0,-1) [right] {$-1$}; 
\node at (0.95,0) [anchor=north west] {$1$}; \node at (-0.95,0) [anchor=north east] {$-1$}; 
\node (phi) at (1,-0.6) {$\phi(\theta)$}; \draw[->, >=stealth, dashed] (phi) -- (0.33,-0.1); 
\draw[color=red, dashed, domain=-1:1] plot (\x,{(1-abs(\x))^2}); 
\node (beta) at (-1,0.7) {$\beta(\theta)$}; 
\draw[->, >=stealth, dashed] (beta) -- (-0.26,0.75);
\draw[->, >=stealth, dashed] (beta) .. controls +(1,1) and (0.66,1.45) .. (0.26,0.75);
\draw[color=red, dashed, domain=-1:1] plot (\x,{-(1-abs(\x))^2}); 
\node (mbeta) at (-1,-0.7) {$-\beta(\theta)$}; 
\draw[->, >=stealth, dashed] (mbeta) -- (-0.26,-0.75);
\draw[->, >=stealth, dashed] (mbeta) .. controls +(1,-1) and (0.66,-1.45) .. (0.26,-0.75);
\color{blue}
  \draw plot file {PseudoCurve_AlsManMor-base.table};
  \draw[xscale=-1, yscale=-1] plot file {PseudoCurve_AlsManMor-base.table};
\end{tikzpicture}
\caption{The graphs of the functions $\phi$ (in \textcolor{blue}{blue}) and $\pm\beta$ in thick black.
The \textcolor{red}{red} dashed curve is $(1-\abs{x})^2.$}\label{fig1}
\end{figure}
\[
\beta(x) := 1 - \abs{x}
\qquad\text{and}\qquad
\phi(x) := (1-\abs{x})^2 \sin\left(\frac{\pi}{x}\right).
\]
Note that $-\beta(x) < \phi(x) < \beta(x),$
for all $x\in [-1,1]\setminus\{0\}$ and the  graphs of $-\beta$ and $\beta$
intersect the closure of the graph of $\phi$ only at the points
$(0,-1),\ (0,1), (-1,0)$  and $(1,0).$

To define the families $\{\basicbox{\ell}\}_{\ell\in\Z}$
and $\{\Gamma\varphi_{_{\lstar}}\}_{\ell\in\Z}$
we use the following \emph{generic boxes}.

For every $\theta\in\SI$ and $\delta < \tfrac{1}{2},$
$\map{\vartheta_{_\theta}}{[-\delta,\delta]}[\SI]$
denotes the map defined by $\vartheta_{_\theta}(x) = x + \theta \pmod{1}.$
Clearly $\vartheta_{\theta}$ is a homeomorphism between
$[-\delta,\delta]$ and $\cball{\theta}{\delta}.$
Finally
$
\map{\vartheta^{-1}_{\theta}}{\cball{\theta}{\delta}}[{[-\delta,\delta]}]
$
denotes the inverse homeomorphism of $\vartheta_{\theta}.$

\begin{definition}[Generic boxes]\label{GenericBoxes}
Fix $\ell,n\in\Z,\ n \ge \all,\ \alpha \in (0,2^{-n}),\ \delta \in (0, \alpha),$
$a\in[-1,1]$ and $a^+,a^- \in \ball{2^{-n}\beta(\delta)}{a}$
(see Figure~\ref{fig-boxes}).
Now we consider the Jordan closed curve in $\Omega,$
formed by the graphs of the functions
\[
a + 2^{-n}(\beta\circ\vartheta^{-1}_{_{\lstar}}) \evalat{\cball{\lstar}{\delta}}
\andq
a- 2^{-n}(\beta\circ\vartheta^{-1}_{_{\lstar}}) \evalat{\cball{\lstar}{\delta}},
\]
together with the four segments that join the points:
\begin{align*}
& \text{$(\lstar -\alpha,a^-)$ with $\left(\lstar-\delta,a-2^{-n}\beta(-\delta)\right)$},\\
& \text{$(\lstar -\alpha,a^-)$ with $\left(\lstar-\delta,a + 2^{-n}\beta(-\delta)\right)$},\\
& \text{$(\lstar + \alpha,a^+)$ with $\left(\lstar + \delta,a-2^{-n}\beta(\delta)\right)$, and}\\
& \text{$(\lstar + \alpha,a^+)$ with $\left(\lstar + \delta,a + 2^{-n}\beta(\delta)\right)$}.
\end{align*}
We denote the closure of the connected component of the complement of the above
Jordan curve in  $\Omega$ that contains the point $(\lstar,a)$
by $\win{\lstar,n,\alpha,\delta,a,a^+,a^-}$
(the coloured region in Figure~\ref{fig-boxes}).
Observe that
$\pi\left(\win{\lstar,n,\alpha,\delta,a,a^+,a^-}\right),$
the projection of $\win{\lstar,n,\alpha,\delta,a,a^+,a^-}$ to $\SI,$
is $\cball{\lstar}{\alpha} = [\lstar-\alpha, \lstar+\alpha].$

\begin{figure}
\begin{tikzpicture}[scale=2.5]
\def\CentRegBound{0.4}
\filldraw[draw=Bittersweet, fill=Bittersweet!50] (-1,-0.25) -- (-\CentRegBound, 1-\CentRegBound) -- (0,1) --
         (\CentRegBound, 1-\CentRegBound) -- (1, 0.15) -- (\CentRegBound, \CentRegBound-1) -- (0,-1) --
         (-\CentRegBound, \CentRegBound-1)  -- (-1,-0.25); 
\draw (-1,-1) -- (-1,1) -- (1,1) -- (1,-1) -- (-1,-1); 

\foreach \x in {-\CentRegBound,0,\CentRegBound} { \draw[dashed, thin] (\x,-1) -- (\x,1); }
\node at (-1,-1) [anchor=north] {$\lstar-\alpha$};
\node at (-\CentRegBound,-1) [anchor=north] {$\lstar-\delta$};
\node at (0,-1) [anchor=north] {$\lstar$};
\node at (\CentRegBound,-1) [anchor=north] {$\lstar+\delta$};
\node at (1,-1) [anchor=north] {$\lstar+\alpha$};

\draw[dashed, thin] (-1,0) -- (1,0); \node at (-1,0) [left] {$a$};

\node (beta) at (-1,1.2) {\footnotesize $a+\tfrac{1}{2^n}(\beta\circ\vartheta^{-1}_{_{\lstar}})(\theta)$} ;
\draw[->, dashed, >=stealth] (-1,1.15) -- (-0.26,0.75);
\draw[->, dashed, >=stealth] (beta) .. controls +(0.7,0.5) and (0.56,1.25) .. (0.26,0.75);
\node (mbeta) at (-0.87,-1.4) {\footnotesize $a-\tfrac{1}{2^n}(\beta\circ \vartheta^{-1}_{_{\lstar}})(\theta)$};
\draw[->, dashed, >=stealth] (-0.92,-1.32) -- (-0.26,-0.75);
\draw[->, dashed, >=stealth] (mbeta) .. controls +(1.6,-0.7) and (0.9,-0.8) .. (0.35,-0.65);

\filldraw[blue] (-1,-0.25) circle(0.6pt) (1,0.15) circle(0.6pt);
\node at (-1,-0.25) [left] {$a^-$};
\node at (1,0.15) [right] {$a^+$};

\node (graph) at (0.6,-0.65) [anchor=west] {$\Gamma\varphi_{_{(\lstar,n,\alpha,\delta,a,a^+,a^-)}}$};
\draw[->, >=stealth, dashed] (graph.north west)+(0.03,-0.03) -- (0.33,-0.17383);

\color{blue}
\begin{scope} \clip (-\CentRegBound,-1) rectangle (\CentRegBound,1);
      \draw plot file {PseudoCurve_AlsManMor-base.table};
      \draw[xscale=-1, yscale=-1] plot file {PseudoCurve_AlsManMor-base.table};
\end{scope}
\draw (-1,-0.25) -- (-\CentRegBound, -0.35006); 
\draw (1,0.15) -- (\CentRegBound, 0.35006); 
\end{tikzpicture}\vspace*{-8ex}
\caption{The region $\win{\lstar,n,\alpha,\delta,a,a^+,a^-}$
is the \textcolor{Bittersweet}{colored filled area}, delimited in the rectangle
$\setsilift{\cball{\lstar}{\delta}}$
by the graphs of the functions $a\pm\tfrac{1}{2^n}(\beta\circ\vartheta^{-1}_{_{\lstar}})(\theta).$
In \textcolor{blue}{blue} the set $\Gamma\varphi_{_{(\lstar,n,\alpha,\delta,a,a^+,a^-)}}$
inductively defining the pseudo-curve.}\label{fig-boxes}
\end{figure}
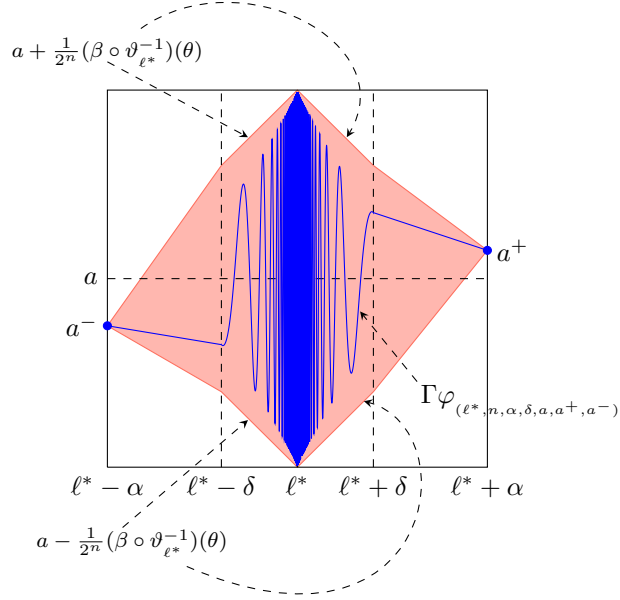

We denote by
\[
\map{ \varphi_{_{\lstar}} = \varphi_{_{(\lstar,n,\alpha,\delta,a,a^+,a^-)}}}{
     \cball{\lstar}{\alpha}\setminus\lstarset}[\I]
\]
the continuous map defined as follows:
\begin{enumerate}[(i)]
\item
$
 \varphi_{_{\lstar}}\evalat{\cball{\lstar}{\delta}\setminus\lstarset} =
     a + (-1)^{\ell} 2^{-n}(\phi \circ \vartheta^{-1}_{_{\lstar}}).
$

\item
$
  \varphi_{_{\lstar}}(\lstar-\alpha) =  a^-
$
and
$
  \varphi_{_{\lstar}}(\lstar + \alpha) = a^+.
$

\item
$
 \varphi_{_{\lstar}}\evalat{[\lstar-\alpha,\lstar-\delta]}
$
and
$
\varphi_{_{\lstar}}\evalat{[\lstar+\delta, \lstar+\alpha]}
$
are affine.
\end{enumerate}
We also denote by
$
  \Gamma\varphi_{_{(\lstar,n,\alpha,\delta,a,a^+,a^-)}}\subset
  \win{\lstar,n,\alpha,\delta,a,a^+,a^-}
$
the closure in $\Omega$ of the graph of
$\varphi_{_{(\lstar,n,\alpha,\delta,a,a^+,a^-)}}.$
\end{definition}

\begin{remark}\label{propiedadesR}
The region
$
\win{\lstar,n,\alpha,\delta,a,a^+,a^-}
$
and the set
$
\Gamma\varphi_{_{(\lstar,n,\alpha,\delta,a,a^+,a^-)}}
$
satisfy the following properties:
\begin{enumerate}[(1)]
\item
$
\win{\lstar,n,\alpha,\delta,a,a^+,a^-} \subset
  \cball{\lstar}{\alpha}\times[a-2^{-n}, a+2^{-n}].
$

\item
$
 \diam(\win{\lstar,n,\alpha,\delta,a,a^+,a^-}) =
 \diam(\setfibpt{\win{\lstar,n,\alpha,\delta,a,a^+,a^-}}{\lstar}) =
 2 \cdot 2^{-n}.
$

\item
The sets $\Gamma\varphi_{_{(\lstar,n,\alpha,\delta,a,a^+,a^-)}}$ and
$\partial\win{\lstar,n,\alpha,\delta,a,a^+,a^-}$
only intersect at the points $(\lstar, a-2^{-n}),(\lstar, a+2^{-n}),$
$(\lstar-\alpha,a^-)$ and $(\lstar + \alpha, a^+).$

\item $
      \setfibpt{\left(\Gamma\varphi_{_{(\lstar,n,\alpha,\delta,a,a^+,a^-)}}\right)}{\lstar}
      =  \win{\lstar,n,\alpha,\delta,a,a^+,a^-}^{\lstar}
      $
      is an interval.

\item Let $\win{\lstar,n,\alpha,\delta,a,a^+,a^-}$ and
          $\win{\kstar,\widetilde{n},\widetilde{\alpha},\widetilde{\delta},\widetilde{a},\widetilde{a}^+,\widetilde{a}^-}$
      be two regions, then
      $\cball{\lstar}{\alpha} \cap \cball{\kstar}{\widetilde{\alpha}} = \emptyset$
      implies
      \[\win{\lstar,n,\alpha,\delta,a,a^+,a^-} \cap \win{\kstar,\widetilde{n},\widetilde{\alpha},\widetilde{\delta},\widetilde{a},\widetilde{a}^+,\widetilde{a}^-} = \emptyset.\]
\end{enumerate}
\end{remark}

For every $j\in \Z^+,$ we set
\begin{align*}
Z_{j}  &:= \set{i\in\Z}{\ai \le j\} = \{-j,-j + 1,\ldots,-1,0,1,\ldots,j-1,j}\text{ and}\\
\Zstar_{j} & := \set{\istar}{i\in Z_{j}}.
\end{align*}

With the help of the sets
$\win{\lstar,n,\alpha,\delta,a,a^+,a^-}$ and
$\Gamma\varphi_{_{(\lstar,n,\alpha,\delta,a,a^+,a^-)}},$
which are the ``bricks'' of our construction we are ready to define
the sequence of pseudo-curve generators
$\{(\gams{j},\SI\setminus \Zstar_{j})\}_{j=0}^{\infty}$
that we are looking for.

To do this, for every $j \ge 0$ we define
\begin{itemize}
\item a strictly increasing sequence $\{n_j\}_{j=0}^\infty \subset \N,$
\item a strictly decreasing sequence $\{\alpha_j\}_{j=0}^\infty$ such that $2^{-n_{j+1}} < \alpha_j < 2^{-n_{j}}$
\item and a sequence $\{\delta_j\}_{j=0}^\infty$ with $2^{-n_{j+1}} < \delta_{j} < \alpha_j$
\end{itemize}
verifying some technical properties that we will make explicit below,
and we define a sequence of boxes
$\basicbox{j}  := \win{\jstar,n_{j},\alpha_{j},\delta_{j},a_j,a_j^+,a_j^-}$ and
$\basicbox{-j} := \win{\sstar{-j},n_{j},\alpha_{j},\delta_{j},a_{-j},a_{-j}^+,a_{-j}^-}$
(for $j=0$ both sets coincide) with projections
\[
\pi\left(\basicbox{j}\right) = \basint{j}
\andq
\pi\left(\basicbox{-j}\right) = \basintneg{j}.
\]
Finally, with the use of all these sequences and objects we can define
our functions $\gams{j}\evalat{\SI\setminus \Zstar_{j}}.$

Observe that we are using the intervals of the form
$\basintabs{\ell},$ $\basintabs[\delta]{\ell}$ and also $\BSG{\ell}{\all-1}$
when $\ell$ is negative.
To ease the use of these intervals we introduce the following notation:
\[
 \wbasint{\ell} := \begin{cases}
                      \basint{\ell} & \text{if $\ell \ge 0$, or}\\
                      \BSG{\ell}{\abs{\ell+1}} & \text{if $\ell < 0$,}
                   \end{cases}
\andq
 \wobasint{\ell} := \begin{cases}
                      \obasint{\ell} & \text{if $\ell \ge 0$, or}\\
                      \OBG{\ell}{\abs{\ell+1}} & \text{if $\ell < 0$.}
                   \end{cases}
\]
Notice that the ball $\wbasint{\ell}$ has diameter $\alpha_j$ for
$\ell \in \{j, -(j+1)\}.$

\begin{remark}\label{Rotation-intervals-formulas}
With the above notation $\basintabs{\ell} \varsubsetneq \wobasint{\ell}$ for every $\ell < 0.$
Moreover, for $\ell \in \Z$ and $j \in \Z^+,$
\begin{align*}
 R_\omega\left(\BSG{\ell}{j}\right)  &= \BSG{\ell+1}{j},\text{ and}\\
 R_\omega\left(\wbasint{\ell}\right) &= \begin{cases}
                      \BSG{\ell+1}{\ell} & \text{if $\ell \ge 0$, or}\\
                      \basintabs{\ell+1} & \text{if $\ell < 0$.}
                   \end{cases}
\end{align*}

Also, the same formulae holds with $\alpha$ replaced by $\delta$
and for open balls.
\end{remark}

The next crucial definition fixes in detail all quantities and objects mentioned above.

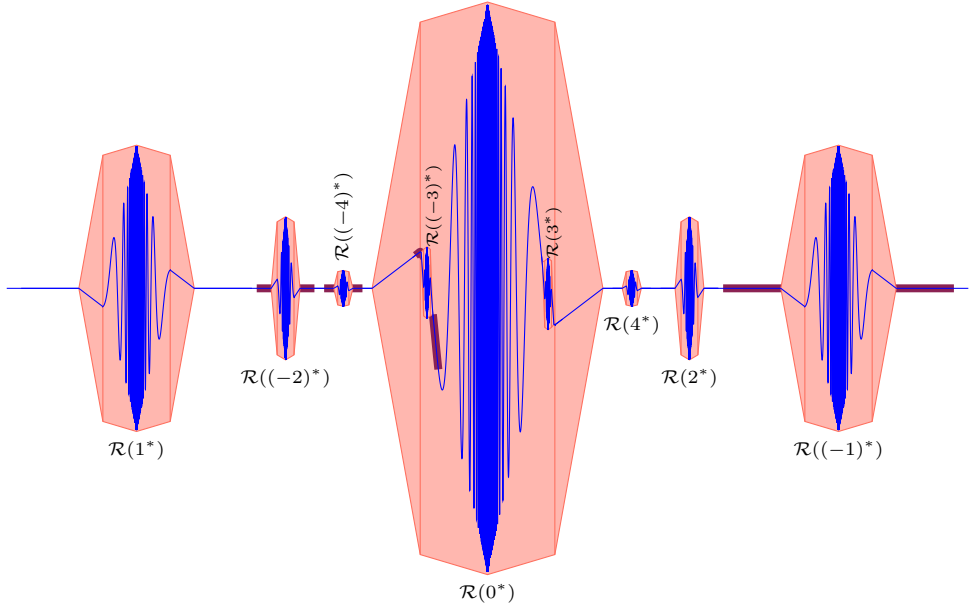
\begin{figure}[hbt]
\begin{tikzpicture}[domain=-0.5:0.5, x=\textwidth, y=0.3\textwidth]
\def\alphazero{0.12}\def\CentRegBound{0.07}
\def\BasicBox#1{\pgfmathsetmacro{\vertor}{(-1)^#1}
 \filldraw[draw=Bittersweet, fill=Bittersweet!50] 
                (-\alphazero,0) -- (-\CentRegBound, 1-\CentRegBound) --
                (0,1) --  (\CentRegBound, 1-\CentRegBound) --
                (\alphazero, 0) -- (\CentRegBound, \CentRegBound-1) --
                (0,-1) -- (-\CentRegBound, \CentRegBound-1)  -- cycle; 
 \foreach \side in { 1, -1 }{ 
     \draw[Bittersweet] (\side*\CentRegBound, 1-\CentRegBound)  -- (\side*\CentRegBound, \CentRegBound-1); 
     \begin{scope} \clip (-\CentRegBound-0.0005,-1) rectangle (\CentRegBound+0.0005,1);
         \draw[blue, xscale=\side*\alphazero, yscale=\side*\vertor] plot file {PseudoCurve_AlsManMor-base.table};
     \end{scope}
     \draw[blue,xscale=\side] (\CentRegBound, -\vertor*\side*0.12791) -- (\alphazero,0) -- (2*\alphazero,0); 
 };}

\draw[blue] (-0.5,0) -- (0.5,0);  
\foreach \j/\xcent/\xs in { 0/0/1, 1/-0.365/0.5, 2/0.21/0.125, 4/0.15/0.083 }{ \pgfmathsetmacro{\ys}{1.0/2^\j}
   \begin{scope}[shift={(\xcent,0)}, xscale=\xs, yscale=\ys]
      \BasicBox{\j}
      \node[below] at (0, -1) {\scriptsize$\win{\sstar{\j}}$};
   \end{scope}
   \ifnum\j>0
      \begin{scope}[shift={(-\xcent,0)}, xscale=\xs, yscale=\ys]
         \draw[garnet, line width=3pt] (-2*\alphazero,0) -- (2*\alphazero,0);\BasicBox{\j}
         \ifnum\j< 4 \node[below] at (0, -1) {\scriptsize$\win{\sstar{-\j}}$};
         \else       \node[above] at (0, 1) {\rotatebox{90}{\scriptsize$\win{\sstar{-\j}}$}};
         \fi
      \end{scope}
   \fi
};
%
\draw[garnet, line width=3pt, join=round]  (-0.063-\alphazero/10,0.115 ) -- (-\CentRegBound, 0.12791) --
                         plot[domain=-\CentRegBound:-\CentRegBound+0.004] (\x,{-((1 + 25*\x/3)^2 - 0.01)*sin(deg(-0.3769911185/\x))});
\draw[garnet, line width=3pt]  plot[domain=-0.05685:-0.063+\alphazero/10] (\x,{-((1 + 25*\x/3)^2 - 0.01)*sin(deg(-0.3769911185/\x))});
\foreach \signel in { 1 , -1 }{
    \begin{scope}[shift={(\signel*0.063,-\signel*0.02)}, xscale=\signel*0.05, yscale=\signel*0.125]
       \filldraw[draw=Bittersweet, fill=Bittersweet!50] 
                (-\alphazero,0.89) -- (-\CentRegBound, 1-\CentRegBound) --
                (0,1) --  (\CentRegBound, 1-\CentRegBound) --
                (\alphazero, -0.85) -- (\CentRegBound, \CentRegBound-1) --
                (0,-1) -- (-\CentRegBound, \CentRegBound-1)  -- cycle; 
       \foreach \side in { 1, -1 }{ 
                \draw[Bittersweet] (\side*\CentRegBound, 1-\CentRegBound)  -- (\side*\CentRegBound, \CentRegBound-1); 
                \begin{scope} \clip (-\CentRegBound-0.0005,-1) rectangle (\CentRegBound+0.0005,1);
                     \draw[blue, xscale=\side*\alphazero, yscale=-\side] plot file {PseudoCurve_AlsManMor-base.table};
                \end{scope}
       };
       \draw[blue, join=round] (\CentRegBound, 0.12791) -- (\alphazero,-0.85) -- (0.14,-0.86) -- (0.3,-0.7); 
       \draw[blue, join=round] (-\CentRegBound, -0.12791) -- (-\alphazero-0.003,0.89) -- (-\alphazero-0.015,1.05); 
    \end{scope}
};
\draw[blue] plot[domain=-0.05685:-0.063+\alphazero/10] (\x,{-((1 + 25*\x/3)^2 - 0.01)*sin(deg(-0.3769911185/\x))}); 
\node[above right] at (0.05, 0.07) {\rotatebox{90}{\tiny$\win{\sstar{3}}$}};
\node[above right] at (-0.075, 0.11) {\rotatebox{90}{\tiny$\win{\sstar{-3}}$}};
\end{tikzpicture}
\caption{The boxes $\win{\lstar}$ for $\ell \in \{-4, -3, -2, -1, 0, 1, 2, 3, 4\}$ and the graph of $\gams{4}.$
The wings are represented as a thick \textcolor{garnet}{garnet} curve surrounding the graph of $\gams{4}.$
For clarity the scale and separation between boxes is not preserved.
The circle $\SI$ is parametrized as $[-\tfrac{1}{2}, \tfrac{1}{2}).$
}\label{IterativeBoxes}
\end{figure}

\begin{definition}\label{PCgenerators}
We start by defining
$\basicbox{0} := \win{\sstar{0},n_{0},\alpha_{0},\delta_{0},0,0,0}$ and
$\varphi_{_{\sstar{0}}} := \varphi_{_{(\sstar{0},n_0,\alpha_0,\delta_0,0,0,0)}}$
by choosing (Definition~\ref{GenericBoxes})
$n_{0}=1,$ $\alpha_{0} < \frac{1}{2} = 2 ^{-n_0}$ and $\delta_0 < \alpha_0$
small enough so that the intervals
$\wbasint{0} = \basint{0},$ $\BSG{1}{0}$ and $\wbasint{-1} = \BSG{-1}{0}$ are pairwise disjoint;
and $\sstar{-2}, \sstar{2} \notin \wbasint{-1}$
and, additionally, $\Bd\left(\basint{0}\right) \cap \Orbom  = \emptyset.$

We also set $a_{0}^+ = a_{0}^- = a_{0} = 0,$
and we define the map $\map{\gams{0}}{\SI\setminus \{0\}}[\I]$ by
\[
\gams{0}(\theta) =
  \begin{cases}
      \varphi_{_{\sstar{0}}}(\theta) & \text{if $\theta \in \basint{0}\setminus\{0\}$,}\\
      0                              & \text{if $\theta \notin \basint{0}$.}
  \end{cases}
\]

For consistency with the definition of $\gams{j}$ in the case $j \ge 1,$ we define
the map $\map{\gams{-1}}{\SI\setminus \{0\}}[\I]$ by
$\gams{-1}(\theta) = 0$ for every $\theta \in \SI.$
Then, notice that, $a_{0} = \gams{-1}(\sstar{0}),$
$
 a^{\pm}_{0} = \varphi_{_{\sstar{0}}}(\sstar{0}\pm\alpha_{0}) = \gams{-1}(\sstar{0}\pm\alpha_{0}),
$
and $\gams{0}(\theta) = \gams{-1}(\theta)$ for every $\theta \notin \basint{0}.$

Next, for every $j\in \N$ we define $\basicbox{j},$ $\basicbox{-j}$
and $(\gams{j},\SI\setminus \Zstar_{j})$
from the corresponding boxes $\basicbox{i}$ and $\basintabs{i} \subset \wbasint{i}$
for $i \in Z_{j-1},$ and $(\gams{j-1},\SI\setminus \Zstar_{j-1})$ as follows.
We take $n_j,$ $\delta_{j}$ and $\alpha_{j}$ such that
(see Figure~\ref{IterativeBoxes} to fix ideas):
\begin{enumerate}[({\tsfR.}1)]
\item $n_j > n_{j-1}$, $\delta_{j} < \alpha_j < 2^{-n_j} < \delta_{j-1} < \alpha_{j-1}$ and
\[
\Bigl( \Bd\left(\basintneg{j}\right) \cup \Bd\left(\basint{j}\right) \Bigr)
    \cap \Orbom  = \emptyset.
\]
\item The intervals
\begin{align*}
   &\wbasint{j} = \basint{j},\\
   &R_\omega\left(\basint{j}\right) = \BSG{j+1}{j},\\
   &\wbasint{-j} = \BSG{-j}{j-1}\text{ and}\\
   &\wbasint{-(j+1)} = \BSG{-(j+1)}{j}
\end{align*}
are pairwise disjoint,
\[
  \gams{j-1} \left(\BSG{\ell}{j}\right) \subset \left[
         \gams{j-1}(\lstar) - 2^{-n_{j}},
         \gams{j-1}(\lstar) + 2^{-n_{j}}
   \right]
\]
for every $\ell \in \{j+1, -(j+1)\},$
\begin{align*}
 & \wbasint{\ell} \cap \Zstar_{j+1} = \lstarset
   \text{ for }
   \ell \in \{j, -(j+1)\} \text{ and}\\
 & \BSG{j+1}{j} \cap \Zstar_{j+1} = \sstarset{j+1},
\end{align*}
and $\sstar{-(j + 2)}, \sstar{j+2} \notin \wbasint{-(j+1)} = \BSG{-(j+1)}{j}.$
\item $\Bd\left(\BSG{k+1}{\ak}\right) \cap \left(\basint{j}\cup \basintneg{j}\right) = \emptyset$
      for every  $k\in Z_{j-1}$.
\item Assume that there exists $k\in Z_{j-1}$ such that
      $\BSG{j+1}{j} \cap \wbasint{k} \ne \emptyset$ and
      $\ak$ is maximal verifying these conditions. Then,
      $\BSG{j+1}{j}$ is contained in one of the two connected components
      of $\obasintabs{k} \setminus \kstarset$ when
      $\BSG{j+1}{j} \cap \basintabs{k} \ne \emptyset$, and
      $\BSG{j+1}{j}$ is contained in one of the two connected components
      of $\wobasint{k}\setminus \basintabs{k}$ if
      $\BSG{j+1}{j} \cap \basintabs{k} = \emptyset$
      (note that, in this case, $k$ must be negative).
\item Let $\ell \in \{j, -(j+1)\}$
(recall that the ball $\wbasint{\ell}$ has diameter $\alpha_j$
 for these two values of $\ell$ and only for them).
      \begin{enumerate}[{(\tsfR.5.}i)]
           \item If $\lstar \notin \bigcup_{i\in Z_{j-1}} \wbasint{i}$ then,
                 $\wbasint{\ell} \cap \wbasint{i} = \emptyset$
                 for every $i \in Z_{j-1}.$
           \item If $\lstar \in \wbasint{m}$ for some $m \in Z_{j-1}$
                 such that $\am$ is maximal with these properties, then
                 \begin{enumerate}[({\tsfR.5.ii.}1)]
                      \item $\wbasint{\ell} \cap \wbasint{i} = \emptyset$ for every
                            $i \in Z_{j-1}$ such that $\ai \ge \am,\ i \ne m,$ and
                      \item $\wbasint{\ell}$ is contained in (a connected component of)\\[-1ex]
                             \begin{minipage}{25em}
                             \begin{multline*}
                             \wobasint{m} \setminus \left(\Bd\left(\basintabs{m}\right) \cup \mstarset\right) =\\
                               \left(\mstar - \alpha_{_{\am - 1}},  \mstar - \alpha_{_{\am}}\right) \cup
                               \left(\mstar - \alpha_{_{\am}},  \mstar\right) \cup\\
                               \left(\mstar,  \mstar + \alpha_{_{\am}}\right) \cup
                               \left(\mstar + \alpha_{_{\am}},  \mstar + \alpha_{_{\am - 1}}\right)
                             \end{multline*}
                             \end{minipage}\\[1ex]
                             (observe that $\wbasint{\ell} \subset \wobasint{m} \setminus \basintabs{m}$
                              can only happen when $m < 0$ since
                              $\wbasint{m} = \basintabs{m}$ for $m \ge 0$).
                 \end{enumerate}
      \end{enumerate}
\item Let $\ell \in \{j, -j\}.$
      If $\wbasint{\ell} \cap \wbasint{m} = \emptyset$
      for every $m\in Z_j,$ $m \ne \ell$ then,
      to define $\basicbox{\ell}$ and the map $\varphi_{_{\lstar}},$ we set
      \[ a_{\ell} = \gams{j-1}(\lstar) = a_\ell^{\pm} = \gams{j-1}(\lstar \pm \alpha_{j}) = 0.\]
      Otherwise, there exists $m\in Z_{j-1}$ such that
      $\wbasint{\ell}$ is contained in a connected component of
      $\wobasint{m} \setminus \left( \Bd\left(\basintabs{m}\right) \cup \mstarset \right)$
      and $\am$ is maximal with these properties. Then,
      to define $\basicbox{\ell}$ and the map $\varphi_{_{\lstar}},$ we set
      \begin{enumerate}[{(\tsfR.6.}i)]
        \item $a_{\ell} := \gams{\am}(\lstar),$
              $a_\ell^{\pm} := \gams{\am}(\lstar \pm \alpha_{j})$ and
              $
              \Graph\Bigl(\gams{\am}\evalat{\BSG{\ell}{j}}\Bigr) \subset \basicbox{\ell}.
              $

        \item Assume that there exists $k \in Z_{\am} \subset Z_{j-1}$ such that
              $\wbasint{\ell} \subset \obasintabs{k} \setminus \kstarset.$
              Then, $\basicbox{\ell}$ is contained in one of the
              two connected components of
              $\Int\left(\basicbox{k} \setminus \setsilift{\sstar{k}}\right).$
      \end{enumerate}
\end{enumerate}
Finally we define $\map{\gams{j}}{\SI\setminus \Zstar_{j}}[\I]$ by
\[
\gams{j}(\theta) =
  \begin{cases}
    \varphi_{_{\jstar}}(\theta)    & \text{if $\theta \in \basint{j} \setminus \sstarset{j},$} \\
    \varphi_{_{\sstar{-j}}}(\theta) & \text{if $\theta \in \basintneg{j} \setminus \sstarset{-j},$} \\
    \gams{j-1}(\theta)                & \text{if $\theta \notin \left( \basint{j} \cup \basintneg{j} \cup \Zstar_{j-1}\right).$}
  \end{cases}
\]
(notice that $\Zstar_{j} = \Zstar_{j-1} \cup \{\jstar,\sstar{-j}\}$).
\end{definition}

For every $\ell \in \Z$ we define the \emph{winged region associated to $\ell$} as
\[
 \wbasicbox{\ell} := \begin{cases}
                      \basicbox{\ell} & \text{if $\ell \ge 0$, or}\\
                      \basicbox{\ell} \cup \Graph\Bigl(\gams{\all}\evalat{\wbasint{\ell} \setminus \obasintabs{\ell}}\Bigr) & \text{if $\ell < 0$.}
                   \end{cases}
\]

The next technical lemma shows that the objects from Definition~\ref{PCgenerators}
exist (that is, they are well defined), and studies some of the
basic properties of the family of pseudo-curve generators
$\{(\gams{i},\SI\setminus \Zstar_{i})\}_{i=0}^{\infty}.$

\begin{remark}[Explicit consequences of Definition~\ref{PCgenerators}]\label{PCgeneratorsExplicitConsequences}
The following statements are easy consequences of Definition~\ref{PCgenerators}.
They are stated explicitly for easiness of usage.
\begin{enumerate}[({\tsfR.}1)]
\item $n_j > j.$
      This follows from Definition~\ref{PCgenerators}(\tsfR.1)
      and the fact that we have set $n_0 = 1$ and $n_j > n_{j-1}$ for $j\in\N.$
\item For every $j \in \N,$
      \[
          \wbasint{-j} \cap \Zstar_{j+1} = \sstarset{-j}.
      \]
      This follows from Definition~\ref{PCgenerators}(\tsfR.2)
      for $j-1.$  We get\\[-1ex]
      \hspace*{1em}\begin{minipage}{31.5em}
            \[
               \wbasint{-j} \cap \Zstar_{j} = \sstarset{-j}
               \andq
               \sstar{-(j + 1)}, \sstar{j+1} \notin \wbasint{-j}.
            \]
      \end{minipage}\\[1ex]
      which shows the statement.
\setcounter{enumi}{5}
\item Let $j \in \N$ and $\ell \in \{j, -j\},$
      and assume that $\wbasint{\ell} \cap \wbasint{m} = \emptyset$
      for every $m\in Z_j,$ $m \ne \ell.$
      Then, $\gams{r}\evalat{\wbasint{\ell}} = \gams{0}\evalat{\wbasint{\ell}} \equiv 0$
      for $r=1,2, \dots, j-1.$
      \begin{enumerate}[{(\tsfR.6.}i)]
        \item Assume that here exists $m\in Z_{j-1}$ such that
              $\wbasint{\ell}$ is contained in a connected component of
              $\wobasint{m} \setminus \left( \Bd\left(\basintabs{m}\right) \cup \mstarset \right)$
              and $\am$ is maximal with these properties. Then,
              $\gams{r}\evalat{\wbasint{\ell}} = \gams{\am}\evalat{\wbasint{\ell}}$
              for $r=\am + 1, \am + 2, \dots, j-1.$

        \item Assume that there exists $k \in Z_{\am} \subset Z_{j-1}$ such that
              $\wbasint{\ell} \subset \obasintabs{k} \setminus \kstarset$ and
              $\ak$ is maximal with these properties. Then,
              $\gams{r}\evalat{\wbasint{\ell}} = \gams{\ak}\evalat{\wbasint{\ell}}$
              for $r=\ak + 1, \ak + 2, \dots, \am.$
      \end{enumerate}
\end{enumerate}
To prove (\tsfR.6) notice that when
$
\wbasint{\ell} \cap \basintabs{m}
      \subset \wbasint{\ell} \cap \wbasint{m}
      = \emptyset
$
for every $m\in Z_j,$ $m \ne \ell,$ from the definition of
$\gams{r}$ for $0 \le r < j$ we get that
$\gams{r}\evalat{\wbasint{\ell}} = \gams{0}\evalat{\wbasint{\ell}} \equiv 0$
for $r=1,2, \dots, j-1.$

\inidemopart{\tsfR.6.i}
The maximality of $\am,$ together with Definition~\ref{PCgenerators}(\tsfR.2),
imply that
$
\wbasint{\ell} \cap \basintabs{i}
        \subset \wbasint{\ell} \cap \wbasint{i}
        = \emptyset
$ for every $i \in Z_{j-1},$ $\ai \ge \am,$ $i \ne m.$
So, by the definition of the functions $\gams{r}$,
\[
\gams{r}\evalat{\BSG{\ell}{j}} = \gams{\am}\evalat{\BSG{\ell}{j}}
\andq[for]
r=\am + 1, \am + 2, \dots, j-1.
\]

\inidemopart{\tsfR.6.ii}
When $\ak = \am$ (\tsfR.6.ii) holds trivially.
So, assume that $\ak < \am.$
As in the case (\tsfR.6.i),
the maximality of $\ak$ and Definition~\ref{PCgenerators}(\tsfR.2)
imply that
$
\wbasint{\ell} \cap \basintabs{r} = \emptyset
$
for every $r \in Z_{j-1},$ $\abs{r} \ge \ak,$ $r \ne k.$
So, (\tsfR.6.ii) follows from the definition of
the functions $\gams{r}$.
\end{remark}

\begin{lemma}\label{Propertiesvarphi}
For every $j\in \Z^+$ the regions
$\basicbox{j}$ and $\basicbox{-j}$ (and hence $\wbasicbox{j}$ and $\wbasicbox{-j}$),
and the maps $(\gams{j},\SI\setminus \Zstar_{j})$ are well defined.
Moreover, the following statements hold:
\begin{enumerate}[(a)]
\item $(\gams{j}, \SI\setminus\Zstar_{j}) \in \C.$
Furthermore, for every $\ell \in \{j+1, -(j+1)\},$
\[
  \gams{j} \left(\BSG{\ell}{j}\right) \subset \left[
         \gams{j}(\lstar) - 2^{-n_{j}},
         \gams{j}(\lstar) + 2^{-n_{j}}
   \right].
\]

\item $\bigcup_{\ell \in \Z} \wbasicbox{\ell} \subset \SI \times [-1,1]$ and
      $\Graph\left(\gams{j}\evalat{\SI\setminus\Zstar_{j}}\right) \subset \SI \times [-1,1].$

\item For $\ell \in \{j, -j\}$ we have
      $\Graph\Bigl(\gams{j-1}\evalat{\BSG{\ell}{j}}\Bigr) \subset \basicbox{\ell},$
      $a_{\ell} = \gams{j-1}(\lstar),$ and
      $
       a^{\pm}_{\ell} =
       \varphi_{_{\lstar}}(\lstar\pm\alpha_{j}) =
       \gams{j-1}(\lstar\pm\alpha_{j}).
      $

\item $\Graph\Bigl(\gams{n}\evalat{\BSG{\ell}{j} \setminus \Zstar_n}\Bigr) \subset \basicbox{\ell}$
      for every $n \ge j$ and $\ell \in \{j, -j\}.$

\item For every $\ell \in \{j, -j\},$
      \[\hspace*{1.5em}
        \gams{j}\evalat{\left(\wbasint{\ell} \setminus \OBG{\ell}{j}\right) \cup R_\omega\left(\wbasint{\ell} \setminus \OBG{\ell}{j}\right)} =
        \gams{j-1}\evalat{\left(\wbasint{\ell} \setminus \OBG{\ell}{j}\right) \cup R_\omega\left(\wbasint{\ell} \setminus \OBG{\ell}{j}\right)}.
      \]
      Moreover, for every
      $\theta \in \Bd(\wbasint{\ell} \setminus \OBG{\ell}{j}) = \Bd(\BSG{\ell}{j}) \cup \Bd(\wbasint{\ell}),$
      we have
      $\theta \notin \wbasint{n} \cup \wbasint{-n}$ and
      $\gams{n}(\theta) = \gams{j}(\theta) = \gams{j-1}(\theta)$
      for every $n > j,$ and
      $R_\omega(\theta) \notin \basint{n} \cup \basintneg{n}$ and
      $\gams{n}\bigl(R_\omega(\theta)\bigr) = \gams{j-1}\bigl(R_\omega(\theta)\bigr)$
      for every $n \ge j.$

\item For every $\ell \in \Z,$ $\wbasicbox{\ell}$ is a compact connected set
      such that $\pi\left(\wbasicbox{\ell}\right) = \wbasint{\ell},$
      $\gams{\all}\evalat{\wbasint{\ell} \setminus \obasintabs{\ell}}$ is continuous
      and
      \[
        \diam\left(\wbasicbox{\ell}\right) = \begin{cases}
            \diam\left(\basicbox{\ell}\right) = \diam\left(\basicbox{-\ell}\right) =
            2 \cdot 2^{-n_{\ell}} \le 2^{-\ell} & \text{if $\ell \ge 0$,}\\
            2 \cdot 2^{-n_{\abs{\ell+1}}} \le 2\cdot 2^{-{\all}} & \text{if $\ell < 0$.}
        \end{cases}
      \]

\item Given $\ell,m\in \Z$ such that
$\all \ge \am,$ $\ell \ne m$ and $\wbasint{\ell} \cap \wbasint{m} \ne \emptyset,$
it follows that $\all > \am,$ and either
$\wbasint{\ell} \subset \obasintabs{m} \setminus \mstarset$
and the region $\wbasicbox{\ell}$ is contained in one of the
two connected components of
$\Int\left(\basicbox{m} \setminus \setsilift{\mstar}\right),$
or $m < 0$ and $\wbasint{\ell}$ is contained in one of
the two connected components of $\wobasint{m} \setminus \basintabs{m}$.
\end{enumerate}
\end{lemma}

\begin{proof}
We start by proving the first statement of the lemma and (a) by induction.

Observe that $n_0 = 1$, $\alpha_0,$ $\delta_0$ and $\gams{0}$ are defined so that
Definition~\ref{PCgenerators}(\tsfR.1--2) for $j=0$
and $(\gams{0}, \SI\setminus\Zstar_{0}) \in \C$
are verified except for the obvious fact that $\wbasint{-j} = \wbasint{j}.$
On the other hand, by construction, $\basint{0}$ is disjoint from
$\BSG{1}{0}$ and $\BSG{-1}{0}.$
Then, by the definition of $\gams{0},$
\[
  \gams{0} \left(\BSG{\ell}{0}\right) = \{0\} \subset [-\tfrac{1}{2}, \tfrac{1}{2}] =
  \left[
         \gams{0}(\lstar) - 2^{-n_{0}},
         \gams{0}(\lstar) + 2^{-n_{0}}
   \right]
\]
for $\ell \in \{1, -1\}.$
Hence, (a) holds.

Fix $j > 0$ and assume that we have defined
$n_{\ell}$, $\alpha_{\ell},$ $\delta_{\ell}$ and $\gams{\ell}$
such that all Definition~\ref{PCgenerators}(\tsfR.1--6) above and (a)
hold for $\ell = 0, 1, \dots, j-1.$

Since the elements of $\Zstar_{j + 2}$ are pairwise different,
we can choose an integer $n_{j} > n_{j-1}$ and $\delta_{j}$ and $\alpha_j$
small enough so that
\begin{itemize}
\item $0 < \delta_{j} < \alpha_j < 2^{-n_j} < \delta_{j-1},$
\item $\sstar{-(j + 2)}, \sstar{j+2} \notin \wbasint{-(j+1)} = \BSG{-(j+1)}{j},$
\item the three intervals
      $\wbasint{j} = \basint{j},$
      $R_\omega\left(\basint{j}\right) = \BSG{j+1}{j}$ and
      $\wbasint{-(j+1)}$
      are pairwise disjoint,
\item $\wbasint{\ell} \cap \Zstar_{j+1} = \lstarset$
      for $\ell \in \{j, -(j+1)\}$,\newline
      $\BSG{j+1}{j} \cap \Zstar_{j+1} = \sstarset{j+1}$
      and, additionally,
\item $
      \Bigl( \Bd\left(\basintneg{j}\right) \cup
           \Bd\left(\basint{j}\right) \Bigr) \cap \Orbom  = \emptyset.
      $
\end{itemize}
Then, Definition~\ref{PCgenerators}(\tsfR.1) is verified.
Moreover, from the above conditions it follows that
$\BSG{\ell}{j} \cap \Zstar_{j+1} = \lstarset$
for every $\ell \in \{j+1, -(j+1)\}.$
Thus, by statement~(a) for $j-1$,
$\gams{j-1}$ is defined and continuous on $\lstar \in \BSG{\ell}{j}$
because this interval is disjoint from $\Zstar_{j-1}.$
Hence, we can decrease the value of $\alpha_j$
(and, accordingly, the value of $0 < \delta_{j} < \alpha_j$),
if necessary, to get\smallskip
\begin{itemize}
\item \hfill$
  \gams{j-1} \left(\BSG{\ell}{j}\right) \subset \left[
         \gams{j-1}(\lstar) - 2^{-n_{j}},
         \gams{j-1}(\lstar) + 2^{-n_{j}}
   \right]
$\hfill\strut\par\par\medskip\par\par\noindent
for every $\ell \in \{j+1, -(j+1)\}.$
\end{itemize}

To see that Definition~\ref{PCgenerators}(\tsfR.2) is verified it remains
to show that the intervals $\wbasint{j},$ $\BSG{j+1}{j}$ and
$\wbasint{-(j+1)}$ are disjoint from $\wbasint{-j}.$
By induction, Definition~\ref{PCgenerators}(\tsfR.2) holds for $j-1$. Thus we see, that
$\sstar{-(j + 1)}, \sstar{j+1} \notin \wbasint{-j},$
and  $R_\omega\left(\basint{j-1}\right) = \BSG{j}{j-1}$
is disjoint from $\wbasint{-j}.$
Hence, we can decrease the value of $\alpha_j$
(and, accordingly, the value of $0 < \delta_{j} < \alpha_j$),
if necessary, until $\BSG{j+1}{j}$ and $\wbasint{-(j+1)} = \BSG{-(j+1)}{j}$
are disjoint from $\wbasint{-j}.$
On the other hand we have that
$\alpha_j < 2^{-n_j} < \delta_{j-1} < \alpha_{j-1}.$ So,
$\wbasint{j} = \basint{j} \subset \BSG{j}{j-1}$ is disjoint from $\wbasint{-j}.$

Up to now we have seen that we can choose
$n_{j},$ $\delta_{j}$ and $\alpha_j$ so that Definition~\ref{PCgenerators}(\tsfR.1--2)
hold for $j$.
Let us see that we can choose $\alpha_j$ such that Definition~\ref{PCgenerators}(\tsfR.3) also holds.
Observe that for every $\ell, i \in \Z$ and every $m\ge 0$ it follows that
$\Bd\left(\BSG{\ell}{m}\right) \cap \Orbom \ne \emptyset$ if and only if
$
  \Bd\left(R^i_\omega\left(\BSG{\ell}{m}\right)\right) \cap \Orbom =
  \Bd\left(\BSG{\ell + i}{m}\right) \cap \Orbom \ne \emptyset.
$
Therefore, by using Definition~\ref{PCgenerators}(\tsfR.1) inductively, we obtain
\[
 \bigcup_{k\in Z_{j-1}} \Bd\left(\BSG{k+1}{\ak}\right) \cap \{\sstar{-j}, \jstar\} \subset
 \bigcup_{k\in Z_{j-1}} \Bd\left(\BSG{k+1}{\ak}\right) \cap \Orbom  = \emptyset.
\]
Consequently, since $\bigcup_{k\in Z_{j-1}} \Bd\left(\BSG{k+1}{\ak}\right)$ is a finite set,
by decreasing again the value of $\alpha_j$, if necessary,
we can achieve that Definition~\ref{PCgenerators}(\tsfR.3) holds for $j$
and Definition~\ref{PCgenerators}(\tsfR.1--2) are still verified.

Next we will take care of Definition~\ref{PCgenerators}(\tsfR.4).
If $\sstar{j+1} \notin \bigcup_{i\in Z_{j-1}} \wbasint{i},$
by decreasing again the value of $\alpha_j$ (and $\delta_j$), if necessary,
we can achieve that
$\BSG{j+1}{j} \cap \left( \bigcup_{i\in Z_{j-1}} \wbasint{i} \right) = \emptyset$
while preserving that Definition~\ref{PCgenerators}(\tsfR.1--3) are verified for $j$.
In this case Definition~\ref{PCgenerators}(\tsfR.4) holds trivially.

Conversely, assume that there exists $k\in Z_{j-1}$ such that
$\sstar{j+1} \in \wbasint{k}$ and $\ak$ is maximal verifying these conditions.
By Definition~\ref{PCgenerators}(\tsfR.2), $k$ is unique (that is, the condition cannot be verified
by $k$ and $-k$ simultaneously).
On the other hand, by the Definition~\ref{PCgenerators}(\tsfR.1) for $\ak$ and $\ak - 1$
and the comment above,
$\sstar{j+1} \notin \Bd\left(\wbasint{k}\right) \cup \Bd\left(\basintabs{k}\right).$
Since $k \in Z_{j-1},$ $\ak \le j-1$ and, hence,
$\sstar{j+1} \notin \Zstar_{\ak}$ (in particular $\jstar \ne \kstar$).
Consequently, $\sstar{j+1}$ is contained in one of the connected components of
$\wobasint{k} \setminus \left( \Bd\left(\basintabs{k}\right) \cup \Zstar_{\ak} \right).$
Then, by decreasing again the value of $\alpha_j$, if necessary,
we can get that $\BSG{j+1}{j}$ is contained in the connected component of
$\wobasint{k} \setminus \left( \Bd\left(\basintabs{k}\right) \cup \Zstar_{\ak} \right)$
where $\sstar{j+1}$ lies, while preserving that
Definition~\ref{PCgenerators}(\tsfR.1--3) are verified for $j$.
Consequently, Definition~\ref{PCgenerators}(\tsfR.1--4) hold for $j$.

Now we will deal with Definition~\ref{PCgenerators}(\tsfR.5).
If $\lstar \notin \bigcup_{i \in Z_{j-1}} \wbasint{i},$
by decreasing again the value of $\alpha_j$, if necessary,
we can get Definition~\ref{PCgenerators}(\tsfR.5.i) while preserving that
Definition~\ref{PCgenerators}(\tsfR.1--4) are verified for $j$.

Assume that there exists $m \in Z_{j-1}$ such that
$\lstar \in \wbasint{m}$ and $\am$ is maximal with these properties.
As in the above construction, by Definition~\ref{PCgenerators}(\tsfR.1--2),
$\lstar \in \wobasint{m} \setminus \left(\Bd\left(\basintabs{m}\right) \cup \mstarset \right)$
and $m$ is unique (that is, the condition cannot be verified
simultaneously by $m$ and $-m$). Consequently,
$\lstar \notin \wbasint{i}$ for every
$i \in Z_{j-1}$ such that $\ai \ge \am,\ i \ne m.$
Thus, by decreasing again the value of $\alpha_j$, if necessary,
we can get that Definition~\ref{PCgenerators}(\tsfR.1--4) still hold,
Definition~\ref{PCgenerators}(\tsfR.5.ii.1) is verified and the interval
$\wbasint{\ell}$ is contained in the connected component of
$\wobasint{m} \setminus \left( \Bd\left(\basintabs{m}\right) \cup \mstarset \right)$
where $\lstar$ lies.
So, Definition~\ref{PCgenerators}(\tsfR.5.ii.2) also holds.

We claim that\\\noindent{\itshape
for every $\ell,m\in \Z$ such that
$\am \le \all \le j,$ $\ell \ne m,$ either
$\wbasint{\ell} \cap \wbasint{m} = \emptyset$
or $\am < \all$ and
$\wbasint{\ell}$ is contained in a connected component of
\[
\wobasint{m} \setminus \left( \Bd\left(\basintabs{m}\right) \cup \mstarset \right).
\]}\\[-1ex]
We prove the claim by induction.
Observe that the claim holds trivially for $\am \le \all \le 1$ because
$\wbasint{0},$ $\wbasint{1} = \basint{1} \subset \BSG{1}{0}$ and
$\wbasint{-1}$ are pairwise disjoint by construction.

Assume that the claim holds for every $\am \le \all < j.$
So, to prove the claim, we may assume that
$\ell \in \{j, -j\},$ $m \in Z_{j-1} \cup \{-\ell\}$
and $\wbasint{\ell} \cap \wbasint{m} \ne \emptyset$.
By Definition~\ref{PCgenerators}(\tsfR.2),
$\wbasint{j} \cap \wbasint{-j} = \emptyset.$
Consequently, $m \ne -\ell$
(that is, $m \in Z_{j-1}$ and $\all = j > \am$).
On the other hand, if $\ell = -j,$ Definition~\ref{PCgenerators}(\tsfR.2) for $j-1$ shows that
$\wbasint{j-1},$ $\wbasint{-(j-1)}$ and  $\wbasint{-j}$
are pairwise disjoint. Thus, $m \in Z_{j-2}$ in this case.

Hence, by the Definition~\ref{PCgenerators}(\tsfR.5) for $j$ when $\ell = j$ and
for $j-1$ when $\ell = -j,$ there exists $k \in Z_{j-1}$
(in fact when $\ell = -j,$ $k \in Z_{j-2}$) such that $\wbasint{\ell}$
is contained in a connected component of
$\wobasint{k} \setminus \left( \Bd\left(\basintabs{k}\right) \cup \kstarset \right)$
and $\all = j > \ak \ge \am.$

If $m = k$ then the claim holds.
Otherwise, $m \ne k$ and since $j=\all > \ak \ge \am$,
by the induction hypotheses, $\ak > \am,$ and  $\wbasint{k}$
is contained in a connected component of
$\wobasint{m} \setminus \left( \Bd\left(\basintabs{m}\right) \cup \mstarset \right).$
So, the claim holds also in this case. This ends the proof of the claim.

Finally, we consider Definition~\ref{PCgenerators}(\tsfR.6).
The fact that either $\wbasint{\ell} \cap \wbasint{m} = \emptyset$
for every $m\in Z_j,$ $m \ne \ell$ or
there exists $m\in Z_{j-1}$ such that
$\wbasint{\ell}$ is contained in a connected component of
$\wobasint{m} \setminus \left( \Bd\left(\basintabs{m}\right) \cup \mstarset \right)$
follows from the claim.

To show that Definition~\ref{PCgenerators}(\tsfR.6.i) can be guaranteed,
it is enough to decrease again the value of $\alpha_j$, if necessary,
until $\BSG{\ell}{j}$ is disjoint from $\Zstar_{\am}$ and
Definition~\ref{PCgenerators}(\tsfR.1--5) are still verified.
Thus by (a) for $\am,$ $\gams{\am}$ is well defined and continuous on
$\BSG{\ell}{j}.$
So, we can set $a_{\ell} := \gams{\am}(\lstar)$ and,
by decreasing again $\alpha_j$ (if necessary), we get
$
\Graph\Bigl(\gams{\am}\evalat{\BSG{\ell}{j}}\Bigr) \subset \basicbox{j}.
$

To show that Definition~\ref{PCgenerators}(\tsfR.6.ii) can be guaranteed we first assume that $k = m.$
As before, if necessary, we can increase the value of $n_{j}$ and, accordingly,
decrease the values of $\alpha_j < 2^{-n_{j}}$ and $0 < \delta_{j} < \alpha_j$
so that Definition~\ref{PCgenerators}(\tsfR.1--5) and (\tsfR.6.i)
are still verified for $j$ and in addition,
\[
(\lstar,a_{\ell}+2^{-n_{j}}), (\lstar,a_{\ell}-2^{-n_{j}}) \in \Int(\basicbox{k})
\]
and the region $\basicbox{\ell}$ is contained in one of the
two connected components of
$\Int\left(\basicbox{k} \setminus \setsilift{\kstar}\right).$

Assume now that $k \ne m$ (recall that $\ak \le \am < j$).
In this case we have $\wbasint{\ell} \subset \wobasint{m} \cap \obasintabs{k}$.
In particular, $\wobasint{m} \cap \obasintabs{k} \ne \emptyset$ and,
by the above claim, $\ak < \am$ and
$\wbasint{\ell} \subset \wbasint{m}$ is contained in a connected component of
$\wobasint{k} \setminus \left( \Bd\left(\basintabs{k}\right) \cup \kstarset \right).$
The fact that  $\wbasint{\ell} \subset \obasintabs{k} \setminus \kstarset$
implies that
$\wbasint{\ell} \subset \wbasint{m} \subset \obasintabs{k} \setminus \kstarset.$
Then, as above we can increase the value of $n_{j}$ and, accordingly,
decrease the values of $\alpha_j < 2^{-n_{j}}$ and $0 < \delta_{j} < \alpha_j$
so that Definition~\ref{PCgenerators}(\tsfR.1--5) and (\tsfR.6.i) are still verified,
\[
(\lstar,a_{\ell}+2^{-n_{j}}), (\lstar,a_{\ell}-2^{-n_{j}}) \in \Int(\basicbox{k})
\]
and the region $\basicbox{\ell}$ is contained in one of the
two connected components of
$\Int\left(\basicbox{k} \setminus \setsilift{\kstar}\right).$

Now assume that $\ak$ is not maximal verifying the assumptions.
Then, there exists $r \in Z_{\am} \subset Z_{j-1}$ such that
$\wbasint{\ell} \subset \obasintabs{r} \setminus \sstarset{r}$
and $\abs{r}$ is maximal with these properties.

We have $\ak \le \abs{r} \le \am < j$ and
\[
\wbasint{r} \cap \wbasint{k} \supset \obasintabs{r} \cap \obasintabs{k} \ne \emptyset
\]
because
$\wbasint{\ell} \subset \obasintabs{r} \cap \obasintabs{k}.$
Then, by the claim, $\ak < \abs{r}$ and $\wbasint{r}$
is contained in a connected component of
$\wobasint{k} \setminus \left( \Bd\left(\basintabs{k}\right) \cup \kstarset \right).$
The fact that  $\wbasint{\ell} \subset \obasintabs{k} \setminus \kstarset$
implies that
$\wbasint{r} \subset \obasintabs{k} \setminus \kstarset.$
By the part already proven and Definition~\ref{PCgenerators}(\tsfR.6.ii) for $\abs{r} < j$ we get
that $\basicbox{\ell}$ is contained in one of the
two connected components of
$\Int\left(\basicbox{r} \setminus \setsilift{\sstar{r}}\right)$
and $\basicbox{r}$ is contained in one of the
two connected components of
$\Int\left(\basicbox{k} \setminus \setsilift{\sstar{k}}\right).$
This shows that Definition~\ref{PCgenerators}(\tsfR.6.ii) can be guaranteed.

Let us prove that (a) holds for $j$.
Since the set $\SI\setminus\Zstar_{j}$ is residual,
to prove that $(\gams{j}, \SI\setminus\Zstar_{j}) \in \C$
we have to show that
$\gams{j}\evalat{\SI\setminus\Zstar_{j}}$
is continuous.
Note that, from Definition~\ref{PCgenerators}(\tsfR.6.ii),
$
  a^{\pm}_{\ell} = \varphi_{_{\lstar}}(\lstar\pm\alpha_{j}) =
  \gams{j-1}(\lstar\pm\alpha_{j}).
$
Hence, the continuity of $\gams{j}\evalat{\SI\setminus\Zstar_{j}}$
follows from the fact that $\gams{j-1}$ is continuous on
$\SI\setminus\Zstar_{j-1} \supset \SI\setminus\Zstar_{j}$ and
the continuity of $\varphi_{_{\jstar}}$ and $\varphi_{_{\sstar{-j}}}$
(Definition~\ref{GenericBoxes}).

This ends the proof of the first statement of the lemma and
the first statement of (a).
For every $\ell \in \{j+1, -(j+1)\},$
from By Definition~\ref{PCgenerators}(\tsfR.1,2) we get:
\begin{align*}
\gams{j-1} \left(\BSG{\ell}{j}\right)
   &\subset \left[
         \gams{j-1}(\lstar) - 2^{-n_{j}},
         \gams{j-1}(\lstar) + 2^{-n_{j}}
   \right]\\
\BSG{\ell}{j} &\text{ is disjoint from
     $\basint{j}$ and $\BSG{-j}{j-1} \supset \BSG{-j)}{j}$, and}\\
\lstarset & \notin \BSG{\ell}{j} \cap \Zstar_{j-1} \subset
     \BSG{\ell}{j} \cap \Zstar_{j+1} = \lstarset.
\end{align*}
So, from the definition of $\gams{j}$
it follows that
\[
   \gams{j}\evalat{\BSG{\ell}{j}} = \gams{j-1}\evalat{\BSG{\ell}{j}}
\]
and, thus, (a) holds.

Statement~(c) follows immediately
from Definition~\ref{PCgenerators}(\tsfR.6)
and Remark~\ref{PCgeneratorsExplicitConsequences}(\tsfR.6).

Next we prove (b,d,e,f,g).

\inidemopart{d} When $n = j,$ we get
$
  \BSG{\ell}{j} \setminus \Zstar_j = \BSG{\ell}{j} \setminus \lstarset
$
from Definition~\ref{PCgenerators}(\tsfR.2). Hence,
$
  \Graph\Bigl(\gams{j}\evalat{\BSG{\ell}{j} \setminus \Zstar_j}\Bigr) \subset \basicbox{\ell}
$
by the definition of $\gams{j}$ (Definition~\ref{PCgenerators}) and the definition of
$\varphi_{_{\lstar}}$ (Definition~\ref{GenericBoxes}).

Now assume that $n > j$ and fix $\theta \in \BSG{\ell}{j} \setminus \Zstar_n.$
We have to show that the point $(\theta,\gams{n}(\theta)) \in \basicbox{\ell}.$
If $\theta \notin \basintabs{m} $ for every $m$ such that $j < \am \le n$
then, by the iterative use of the definition of
$\gams{i}$ for $i=j+1,j+2,\dots, n$ (Definition~\ref{PCgenerators})
and Definition~\ref{GenericBoxes},
\[
   (\theta,\gams{n}(\theta)) =  (\theta,\gams{n-1}(\theta)) = \dots =
    (\theta,\gams{j+1}(\theta)) = (\theta,\gams{j}(\theta)) =
   (\theta,\varphi_{_{\lstar}}(\theta)) \in \basicbox{\ell}.
\]
Otherwise, by Definition~\ref{PCgenerators}(\tsfR.2),
there exists $m \in \Z$ such that $\all < \am \le n,$
$
\theta \in \basintabs{m} \setminus \Zstar_n,
$
and
$\theta \notin \basintabs{s}$ for every $s$ such that $\am < \abs{s} \le n.$
This implies that
$\wbasint{\ell} \cap \wbasint{m} \supset
  \BSG{\ell}{j} \cap \basintabs{m} \ne
  \emptyset
$ and $\am$ is maximal with these properties.
So, by the claim for $j = \am,$
$\wbasint{m}$ is contained in a connected component of
$\wobasint{\ell} \setminus \left( \Bd\left(\basintabs{\ell}\right) \cup \lstarset \right)$.
Moreover, since $\theta \in \wobasint{m} \cap \BSG{\ell}{j} \ne \emptyset,$
$\wbasint{m} \subset \basintabs{\ell} \setminus \lstarset$.
Thus, by Definition~\ref{PCgenerators}(\tsfR.6.ii) and
Remark~\ref{PCgeneratorsExplicitConsequences}(\tsfR.6.ii) for $j = \am,$
$\ell$ replaced by $m$ and $k$ replaced by $\ell$,
$\basicbox{m} \subset \basicbox{\ell}$
and (d) follows from the part already proven
by replacing $\ell$ by $m$ and $j$ by $\am.$

\inidemopart{g} By the claim we have that for every $\ell,m\in \Z$ such that
$\all \ge \am,$ $\ell \ne m$ and $\wbasint{\ell} \cap \wbasint{m} \ne \emptyset,$
it follows that $\all > \am,$ and
$\wbasint{\ell}$ is contained in a connected component of
$\wobasint{m} \setminus \left( \Bd\left(\basintabs{m}\right) \cup \mstarset \right)$.
Only it remains to show that if
$\wbasint{\ell} \subset \obasintabs{m} \setminus \mstarset,$
then the region $\wbasicbox{\ell}$ is contained in one of the
two connected components of
$\Int\left(\basicbox{m} \setminus \setsilift{\mstar}\right).$
By Definition~\ref{PCgenerators}(\tsfR.6.ii) we know that this holds for
$\basicbox{\ell}$ instead of $\wbasicbox{\ell}.$
Hence, if $\ell \ge 0,$ (g) holds because $\wbasicbox{\ell} = \basicbox{\ell}.$
Assume now that $\ell < 0.$ Since
$\wbasicbox{\ell} = \basicbox{\ell} \cup \Graph\Bigl(\gams{\all}\evalat{\wbasint{\ell} \setminus \obasintabs{\ell}}\Bigr)$
is connected,
$\basicbox{\ell} \subset \basicbox{m},$ and
$\Int\left(\basicbox{m} \setminus \setsilift{\mstar}\right)$
has two connected components, it is enough to show that
\[
\Graph\Bigl(\gams{\all}\evalat{\wbasint{\ell} \setminus \obasintabs{\ell}}\Bigr)
\subset \basicbox{m}.
\]
Since
$
  \wbasint{\ell} \setminus \obasintabs{\ell} \subset
  \wbasint{\ell} \subset
  \obasintabs{m} \setminus \mstarset,
$
statement (g) follows from (d) with $\ell$ replaced by $m$, $j$ by $\am$ and
$n$ replaced by $\all.$

\inidemopart{b}
With (g) in mind we set
\[
 \mathsf{D} := \set{\ell \in \Z}{\wbasicbox{\ell} \not\subset \basicbox{i} \text{ for every }i \in \Z\setminus\{\ell\}}.
\]
Clearly,
\begin{align*}
\bigcup_{\ell \in \Z} \wbasicbox{\ell} &=
 \left(\bigcup_{i \in \Z\setminus \mathsf{D}} \wbasicbox{i}\right) \cup
 \left(\bigcup_{\ell \in \mathsf{D}} \wbasicbox{\ell}\right) \\
&\subset
 \left(\bigcup_{i \in \mathsf{D}} \basicbox{i}\right) \cup
 \left(\bigcup_{\ell \in \mathsf{D}} \wbasicbox{\ell}\right)
 = \bigcup_{\ell \in \mathsf{D}} \wbasicbox{\ell}
\end{align*}
\begin{case}{Claim:}
For every $\ell \in \mathsf{D},$
$\gams{\all-1}\evalat{\wbasint{\ell} \setminus \obasintabs{\ell}} \equiv 0.$
\end{case}

First we prove statement (b) from the above claim and then we will prove the claim.
To this end we start by pointing out few elementary facts.

From the definition of $\wbasicbox{\ell}$ we see that
$\wbasicbox{\ell} \setminus \basicbox{\ell} = \emptyset$ for every $\ell \ge 0$
and
$
  \wbasicbox{\ell} \setminus \basicbox{\ell} \subset
  \Graph\Bigl(\gams{\all}\evalat{\wbasint{\ell} \setminus \obasintabs{\ell}}\Bigr)
$
for every $\ell < 0.$
So, in any case,
\[
  \wbasicbox{\ell} \setminus \basicbox{\ell} \subset
  \Graph\Bigl(\gams{\all}\evalat{\wbasint{\ell} \setminus \obasintabs{\ell}}\Bigr)
  \andq[for every]
  \ell \in \Z.
\]

On the other hand, the arc
$\wbasint{\ell} \supset \wbasint{\ell} \setminus \obasintabs{\ell}$
is disjoint from the arc
$\wbasint{-\ell} \supset \BSG{-\ell}{\all}$
by Definition~\ref{PCgenerators}(\tsfR.2). Thus,
by Definition~\ref{PCgenerators} and (a),
\[
  \gams{\all-1}\evalat{\wbasint{\ell} \setminus \obasintabs{\ell}} =
  \gams{\all}\evalat{\wbasint{\ell} \setminus \obasintabs{\ell}}.
\]

Furthermore, by the Claim and Definition~\ref{PCgenerators}(\tsfR.6),
$a_{\ell}^+ = a_{\ell}^- = a_{\ell} = 0$ for every $\ell \in D.$
So, by Remark~\ref{propiedadesR}(1),
\[
  \basicbox{\ell} \subset
      \basintabs{\ell} \times [-2^{-n_{\all}}, 2^{-n_{\all}}] \subset
      \basintabs{\ell} \times [-2^{-\all}, 2^{-\all}] \subset
      \basintabs{\ell} \times [-1,1].
\]

Therefore, summarizing and using again by the Claim,
\begin{align*}
\bigcup_{\ell \in \Z} \wbasicbox{\ell} &\subset
\bigcup_{\ell \in \mathsf{D}} \wbasicbox{\ell} \subset
 \bigcup_{\ell \in \mathsf{D}} \left( \basicbox{\ell} \cup
      \Graph\Bigl(\gams{\all}\evalat{\wbasint{\ell} \setminus \obasintabs{\ell}}\Bigr)
 \right)\\
&=
 \left(\bigcup_{\ell \in \mathsf{D}} \basicbox{\ell} \right) \cup \left(\bigcup_{\ell \in \mathsf{D}}
    \Graph\Bigl(\gams{\all-1}\evalat{\wbasint{\ell} \setminus \obasintabs{\ell}}\Bigr)\right)\\
&
 \subset \left(\bigcup_{\ell \in \mathsf{D}} \basintabs{\ell} \right) \times [-1,1] \cup \SI \times \{0\}
 \subset \SI \times [-1,1].
\end{align*}
So, the first part of (b) is proved, provided that the claim holds.
Let us prove the second statement of (b).
Observe that, since
\[
 \left(\bigcup_{\ell \in \Z} \basicbox{\ell} \right) \cup \SI \times \{0\} \subset
 \left(\bigcup_{\ell \in \Z} \wbasicbox{\ell} \right) \cup \SI \times \{0\}
 \subset \SI \times [-1,1],
\]
it is enough to show that
\[
 \Graph\left(\gams{j}\evalat{\SI\setminus\Zstar_{j}}\right) \subset
 \left(\bigcup_{\ell \in \Z} \basicbox{\ell} \right) \cup \SI \times \{0\}
\]
for every $j \in \Z^+.$
We will prove this statement by induction on $j.$

By construction we have
\[
 \Graph\left(\gams{0}\evalat{\SI\setminus\{\sstar{0}\}}\right) \subset
 \basicbox{0} \cup \SI \times \{0\} \subset
 \left(\bigcup_{\ell \in \Z} \basicbox{\ell} \right) \cup \SI \times \{0\}.
\]
So, the statement holds for $j=0$. Now assume that it holds for some $j \ge 0,$
and prove it for $j+1.$
By Definition~\ref{PCgenerators} and (d),
\begin{align*}
 \Graph\left(\gams{j+1}\evalat{\SI\setminus\Zstar_{j+1}}\right) &\subset
      \basicbox{j} \cup \basicbox{-j} \cup
           \Graph\left(\gams{j}\evalat{\SI\setminus\Zstar_{j}}\right) \\
 &\subset \basicbox{j} \cup \basicbox{-j} \cup
           \left(\bigcup_{\ell \in \Z} \basicbox{\ell} \right) \cup \SI \times \{0\}\\
 &\subset \left(\bigcup_{\ell \in \Z} \basicbox{\ell} \right) \cup \SI \times \{0\}.
\end{align*}

To end the proof of (b) it remains to show the Claim.

Let $\ell \in \mathsf{D}$ and $m \in Z_{\all},$ $m \ne \ell.$
Then, either
\begin{equation}\label{DichotomyforD}
\begin{cases}
 \wbasint{\ell} \cap \wbasint{m} = \emptyset\text{ or}\\
 \all > \am,\ m < 0 \text{ and } \wbasint{\ell} \subset \wobasint{m} \setminus \basintabs{m}.
\end{cases}
\end{equation}
To see this, observe that if $\wbasint{\ell} \cap \wbasint{m} \ne \emptyset$
then, by (g), $\all > \am$ and either $\wbasicbox{\ell} \subset \basicbox{m}$ or
$m < 0$ and $\wbasint{\ell} \subset \wobasint{m} \setminus \basintabs{m},$
and the first possibility is ruled out because $\ell \in \mathsf{D}.$

By using iteratively the dichotomy \eqref{DichotomyforD} we get that, for every $\ell \in \mathsf{D},$
there exists a sequence $m_0,m_1,\dots,m_k = \ell \in \Z$ with $k \ge 0$ such that
$\wbasint{m_0} \cap \wbasint{q} = \emptyset$ for every $q \in Z_{\abs{m_0}},\ q \ne m_0$
and, in the case $k > 0,$
$
\abs{m_0} < \abs{m_1} < \dots < \abs{m_k} = \all
$
and, for every $p=0,1,\dots,k-1,$
\begin{itemize}
   \item $m_{p} < 0,$
   \item $\wbasint{m_{p+1}} \subset \wobasint{m_p} \setminus \basintabs{m_p}$ and
   \item $\wbasint{m_{p+1}} \cap \wbasint{q} = \emptyset$
for every $q\in Z_{\abs{m_{p+1}}},$ $q \ne m_p, m_{p+1}$ and $\abs{m_p} \le \aq.$
\end{itemize}

The condition
$\wbasint{m_0} \cap \wbasint{q} = \emptyset$ for every $q \in Z_{\abs{m_0}},\ q \ne m_0$
implies
\[
 \gams{\abs{m_0}-1}\evalat{\wbasint{m_0}} = \gams{\abs{m_0}-2}\evalat{\wbasint{m_0}} =
 \dots = \gams{0}\evalat{\wbasint{m_0}} \equiv 0
\]
by Definition~\ref{PCgenerators}(\tsfR.6) and
Remark~\ref{PCgeneratorsExplicitConsequences}(\tsfR.6)
(with $\ell = m_0$). This ends the proof of the Claim when $k = 0.$

Assume now that $k > 0.$
As before we have
\[
  \gams{\abs{m_0}-1}\evalat{\wbasint{m_0} \setminus \obasintabs{m_0}} =
  \gams{\abs{m_0}}\evalat{\wbasint{m_0} \setminus \obasintabs{m_0}}.
\]
This, together with the inclusion,
\[
 \wbasint{m_{1}} \subset \wobasint{m_0} \setminus \basintabs{m_0}
\]
implies that
\[
 \gams{\abs{m_0}}\evalat{\wbasint{m_1}} \equiv 0.
\]
Then, by Definition~\ref{PCgenerators}(\tsfR.6.i) and
Remark~\ref{PCgeneratorsExplicitConsequences}(\tsfR.6.i) with $\ell = m_1$,
\[
 0 \equiv \gams{\abs{m_0}}\evalat{\wbasint{m_1}} = \gams{\abs{m_0}+1}\evalat{\wbasint{m_1}} =
 \dots = \gams{\abs{m_1}-1}\evalat{\wbasint{m_1}}.
\]

If $k = 1$ we are done. Otherwise, $k \ge 2$ and, as above,
\[
 \gams{\abs{m_1}}\evalat{\wbasint{m_2}} \equiv 0.
\]

By iterating the above arguments at most $k$ times the Claim holds.
This ends the proof of (b).

\inidemopart{e}
By Definition~\ref{PCgenerators}(\tsfR.2) and
Remark~\ref{PCgeneratorsExplicitConsequences}(\tsfR.2)
it follows that
\[
\theta \notin \Zstar_{j+1} \cup \OBG{\ell}{j} \cup \wobasint{-\ell}
\andq[for every]
\theta \in \wbasint{\ell} \setminus \OBG{\ell}{j}.
\]
So, by (a), $\gams{j-1}(\theta)$ is well defined and
$\gams{j-1}$ is continuous at $\theta.$
Thus, by the definition of $\gams{j}$ (Definition~\ref{PCgenerators})
and the continuity of $\gams{j-1}$ at $\theta,$
$\gams{j}(\theta) = \gams{j-1}(\theta).$

Now assume that
$\theta \in \Bd(\wbasint{\ell} \setminus \OBG{\ell}{j}) = \Bd(\BSG{\ell}{j}) \cup \Bd(\wbasint{\ell}).$
By (g), $\theta \notin \wbasint{n} \cup \wbasint{-n}$
for every $n > j.$
So, by the iterative use of the definition of
$\gams{i}$ for $i=j+1,j+2,\dots, n$ (Definition~\ref{PCgenerators})
we get
\[
  \gams{j}(\theta) = \gams{j+1}(\theta) = \dots =
  \gams{n-1}(\theta) = \gams{n}(\theta).
\]

Now we prove the part of (e) concerning $R_\omega(\wbasint{\ell} \setminus \OBG{\ell}{j}).$
We first assume that $\ell = j \ge 0.$ Then,
\[
   \wbasint{j} = \basint{j},\
   \theta \in \Bd(\basint{j})
   \andq
   R_\omega(\theta) \in \Bd(\BSG{j + 1}{j}).
\]
Again by Definition~\ref{PCgenerators}(\tsfR.2),
$R_\omega(\theta) \notin \Zstar_{j+1} \cup \basint{j} \cup \wbasint{-j}.$
So, by (a) and the definition of $\gams{j}$ (Definition~\ref{PCgenerators}),
$\gams{j-1}\bigl(R_\omega(\theta)\bigr)$ is well defined and
$\gams{j}\bigl(R_\omega(\theta)\bigr) = \gams{j-1}\bigl(R_\omega(\theta)\bigr).$
By Definition~\ref{PCgenerators}(\tsfR.3)
(with $j = n$ and $k = \ell = j$),
$R_\omega(\theta) \notin \basint{n} \cup \basintneg{n}$
for every $n > j.$
So,
$\gams{n}\bigl(R_\omega(\theta)\bigr) = \gams{j}\bigl(R_\omega(\theta)\bigr)$
as above.

Assume now that $\ell = -j < 0.$
In this case we have $\wbasint{\ell} = \BSG{\ell}{\abs{\ell+1}}$
and, hence,
$
R_\omega(\theta) \in \basintabs{\ell + 1} \setminus \OBG{\ell + 1}{j}.
$
By Definition~\ref{PCgenerators}(\tsfR.1) we have
\[
  \BSG{\ell + 1}{j} \subset
  \BSG{\ell+1}{\abs{\ell+1}} \subset
  \wbasint{\ell+1}.
\]
Thus,
$R_\omega(\theta) \in \wbasint{\ell+1} \setminus \sstarset{\ell+1}.$
Again by Definition~\ref{PCgenerators}(\tsfR.2) and
Remark~\ref{PCgeneratorsExplicitConsequences}(\tsfR.2)
(with $j$ replaced by $-(\ell + 1)$),
\[
  R_\omega(\theta) \notin
     \Zstar_{\ell} \cup \BSG{-\ell}{-(\ell+1)} \cup \wbasint{\ell} \supset
     \Zstar_{j} \cup \basint{j} \cup \wbasint{-j}.
\]
So, by (a) and the definition of $\gams{j}$ (Definition~\ref{PCgenerators}),
$\gams{j-1}\bigl(R_\omega(\theta)\bigr)$ is well defined and
$\gams{j}\bigl(R_\omega(\theta)\bigr) = \gams{j-1}\bigl(R_\omega(\theta)\bigr).$

To end the proof of (e), assume as above that
$\theta \in \Bd(\BSG{\ell}{j}) \cup \Bd(\wbasint{\ell})$
and, hence,
$R_\omega(\theta) \in \Bd(\BSG{\ell + 1}{j}) \cup \Bd\left(\basintabs{\ell + 1}\right).$
We have to show that, in this case,
$R_\omega(\theta) \notin \basint{n} \cup \basintneg{n}$
for every $n > j$
(the fact that $\gams{n}\bigl(R_\omega(\theta)\bigr) = \gams{j}\bigl(R_\omega(\theta)\bigr)$
follows as above).
When $R_\omega(\theta) \in \Bd(\BSG{\ell + 1}{j})$ this follows from
Definition~\ref{PCgenerators}(\tsfR.3) as before.
Assume now that
$R_\omega(\theta) \in \Bd\left(\basintabs{\ell + 1}\right).$
Then, by (g),
$R_\omega(\theta) \notin \wbasint{n} \cup \wbasint{-n}$ for every $n > j.$

\inidemopart{f} If $\ell \ge 0$ then the first two statements of (f)
follow directly from the definitions.
Moreover, by Remarks~\ref{propiedadesR}(2) and
\ref{PCgeneratorsExplicitConsequences}(\tsfR.1),
\[
\diam\left(\wbasicbox{\ell}\right) =
       \diam\left(\basicbox{\ell}\right) =
       \diam\left(\basicbox{-\ell}\right) =
           2 \cdot 2^{-n_{\ell}} \le 2 \cdot 2^{-(\ell+1)} = 2^{-\ell}.
\]

Assume that $\ell < 0.$
From Definition~\ref{PCgenerators}(\tsfR.2) and
Remark~\ref{PCgeneratorsExplicitConsequences}(\tsfR.2)
we get
$\left(\wbasint{\ell} \setminus \obasintabs{\ell}\right) \cap \Zstar_{\all} = \emptyset$
and, hence, $\gams{\all}$
is continuous in an open neighbourhood of
$\wbasint{\ell} \setminus \obasintabs{\ell}$
by (a). On the other hand, by (d),
$\left(\theta,\gams{\all}(\theta)\right) \in \basicbox{\ell}$
for every
$
 \theta \in \Bd\left(\basintabs{\ell}\right)\subset
    \wbasint{\ell} \setminus \obasintabs{\ell}.
$
Thus,
\[
\wbasicbox{\ell} = \basicbox{\ell} \cup \Graph\Bigl(\gams{\all}\evalat{\wbasint{\ell} \setminus \obasintabs{\ell}}\Bigr)
\]
is closed, connected and projects onto the whole $\wbasint{\ell}.$

On the other hand, by (e) and (a)
(since $\ell < 0$, $\abs{\ell+1} = \all -1$),
\begin{align*}
\gams{\all}\left(\wbasint{\ell} \setminus \obasintabs{\ell}\right)
   & = \gams{\all-1}\left(\BSG{\ell}{\abs{\ell+1}} \setminus \obasintabs{\ell}\right)\\
   & \subset \left[
         \gams{\all - 1}(\lstar) - 2^{-n_{\all - 1}},
         \gams{\all - 1}(\lstar) + 2^{-n_{\all - 1}}
   \right].
\end{align*}
Thus, by Remark~\ref{propiedadesR}(1), (c) and Definition~\ref{PCgenerators}(\tsfR.1),
\begin{align*}
\wbasicbox{\ell}
   & = \basicbox{\ell} \cup \Graph\Bigl(\gams{\all}\evalat{\wbasint{\ell} \setminus \obasintabs{\ell}}\Bigr)\\
   &\subset \basintabs{\ell} \times \left[\gams{\all - 1}(\lstar) - 2^{-n_{\all}}, \gams{\all - 1}(\lstar) + 2^{-n_{\all}}\right] \cup\\
   & \hspace*{2em}\left(\BSG{\ell}{\abs{\ell+1}} \setminus \obasintabs{\ell}\right) \times \left[
         \gams{\all - 1}(\lstar) - 2^{-n_{\all - 1}},
         \gams{\all - 1}(\lstar) + 2^{-n_{\all - 1}}
     \right] \\
   &\subset \BSG{\ell}{\abs{\ell+1}} \times \left[
         \gams{\all - 1}(\lstar) - 2^{-n_{\all - 1}},
         \gams{\all - 1}(\lstar) + 2^{-n_{\all - 1}}
     \right].
\end{align*}
Hence, by Definition~\ref{PCgenerators}(\tsfR.1) and
Remark~\ref{PCgeneratorsExplicitConsequences}(\tsfR.1),
\[
 \diam\left(\wbasicbox{\ell}\right) \le
     2\cdot \max\{\alpha_{\abs{\ell+1}}, 2^{-n_{\all - 1}}\}
                 = 2\cdot 2^{-n_{\all - 1}} \le 2 \cdot 2^{-{\all}}.
\]
\end{proof}

The next results allow us to define the limit pseudo-curve
generated by the sequence $\{(\gams{i},\SI\setminus \Zstar_{i})\}_{i=0}^{\infty}.$

\begin{lemma}\label{convergencia}
The sequence
$\{(\gams{i},\SI\setminus \Zstar_{i})\}_{i=0}^{\infty} \subset \C$
is convergent in $\C.$
\end{lemma}

\begin{proof}
By Proposition~\ref{Ccompleto} it suffices to show that
$\{(\gams{i},\SI\setminus \Zstar_{i})\}_{i=0}^{\infty}$
is a Cauchy sequence in $\C.$
By the definition of $\gams{i}$ (Definition~\ref{PCgenerators}) we have
\begin{align*}
\dinf\left(\gams{i-1},\gams{i}\right)
  & = \leftlimits{\sup}{\theta\in \SI\setminus \Zstar_{i}} \abs{\gams{i-1}(\theta)-\gams{i}(\theta)} \\
  & = \leftlimits{\sup}{\theta \in \left(\basint{i} \setminus \istarset\right) \cup
                     \left(\basintneg{i}\setminus \sstarset{-i}\right)
    } \abs{\gams{i-1}(\theta)-\gams{i}(\theta)}.
\end{align*}
By Lemmas~\ref{Propertiesvarphi}(c,d),
and Definition~\ref{PCgenerators}(\tsfR.2) and
Remark~\ref{PCgeneratorsExplicitConsequences}(\tsfR.2),
\[
(\theta,\gams{i-1}(\theta)), (\theta, \gams{i}(\theta)) \in \basicbox{\ell}
\andq[for]
\text{$\theta \in \BSG{\ell}{i} \setminus \lstarset$ and $\ell \in \{i,-i\}.$}
\]
Hence, by Lemma~\ref{Propertiesvarphi}(f),
\[
  \dinf(\gams{i-1},\gams{i}) \le \diam(\basicbox{i}) =
     \diam(\basicbox{-i}) \le 2^{-i}.
\]
Since $n_{i}$ is a strictly increasing sequence, for every $m \ge 0,$
\[
 \dinf(\gams{i+m},\gams{i}) \le
  \sum^{i+m}_{k=i+1} 2^{-k} <
  2^{-(i+1)} \sum_{k=0}^{\infty} \tfrac{1}{2^k} = 2\cdot  2^{-(i+1)},
\]
and consequently
$
\{(\gams{i},\SI\setminus \Zstar_{i})\}_{i=0}^{\infty}
$
is a Cauchy sequence in $\C.$
\end{proof}

Lemma~\ref{convergencia} allows us to define the following limit pseudo-curve generator
of the sequence $\{(\gams{i},\SI\setminus \Zstar_{i})\}_{i=0}^{\infty}.$

\begin{definition}\label{gammalimit}
There exists $(\gamma,\SI\setminus \Orbom) \in \C$ such that
\[
 (\gamma,\SI\setminus \Orbom)  = \lim_{i\to\infty} (\gams{i},\SI\setminus \Zstar_{i})
\]
(that is,
$\gamma(\theta) = \lim_{i\to\infty} \gams{i}(\theta)$  for every $\theta \in \SI\setminus \Orbom$).
Observe that
\[
  \SI\setminus \Orbom = \bigcap_{i=1}^{\infty} \left(\SI \setminus \Zstar_{i} \right)
\]
is a residual set in $\SI.$
\end{definition}

Now, we are ready to define the sequence of pseudo-curves associated to
the sequence
$\{(\gams{i},\SI\setminus \Zstar_{i})\}_{i=0}^{\infty},$
and to the limit pseudo-curve generator $(\gamma,\SI\setminus \Orbom).$
This will finally define the pseudo-curve $\A$ that we want to construct.

\begin{definition}\label{PCAtLast}
We denote by
\[
\A_{j} := \pc[\gams{j},\SI\setminus\Zstar_{j}] =
     \overline{\Graph(\gams{j},\SI\setminus \Zstar_{j})}
\]
the pseudo-curve defined by
$(\gams{j},\SI\setminus \Zstar_{j}) \in \C$, and
\[
 \A = \pc[\gamma,\SI\setminus \Orbom] := \overline{\Graph(\gamma,\SI\setminus \Orbom)}.
\]
By Definition~\ref{gammalimit} and Proposition~\ref{Hdimdinfty},
$
\A = \lim_{j\to\infty} \pc[\gams{j},\SI\setminus\Zstar_{j}].
$
\end{definition}

The next lemmas study the properties the pseudo-curves $\A_{j}$ and $\A.$

\begin{lemma}\label{propiedadesA}
The following statements hold for every $\ell\in \Z$:
\begin{enumerate}[(a)]
\item $\setfibth{\A_n} \subset \setfibth{\basicbox{\ell}}$
      for every $n \ge \all-1$ and $\theta \in \basintabs{\ell}.$
\item $\setfibpt{\A_n}{\lstar} = \setfibpt{\A_{\all}}{\lstar} \subset \setfibbb{\ell}$
      for every $n \ge \all.$
      Moreover, $\setfibpt{\A_{\all}}{\lstar} = \setfibbb{\ell}$
      is a non-degenerate interval.
\item $\setfibth{\A_{\ell}} = \{(\theta, \gams{\ell}(\theta)\}$
      for every $\theta \in \SI\setminus\Zstar_{\ell}.$
\item $\A_{\all} \subset \SI \times [-1,1].$
\end{enumerate}
\end{lemma}

\begin{proof}
\inidemopart{a}
By Lemma~\ref{Propertiesvarphi}(c,d),
$
\Graph\Bigl(\gams{n}\evalat{\basintabs{\ell} \setminus \Zstar_n}\Bigr) \subset \basicbox{\ell}.
$
Then, the statement follows from the compacity of $\basicbox{\ell}$.

\inidemopart{b}
From the definition of $\gams{i}$ (Definition~\ref{PCgenerators}) and
Definition~\ref{PCgenerators}(\tsfR.2),
for every $n > \all$ there exists an $\varepsilon(n) > 0$ such that
$\gams{n}(\theta) = \gams{\all}(\theta)$ for every
$\theta \in \ball{\lstar}{\varepsilon(n)} \setminus \lstarset.$
Hence $\setfibpt{\A_n}{\lstar} = \setfibpt{\A_{\all}}{\lstar}.$
Moreover,
$\gams{\all}$ coincides with $\varphi_{_{\lstar}}$ in a neighbourhood of $\lstar.$
Thus, $\setfibpt{\A_{\all}}{\lstar} = \setfibbb{\ell}$ and it is an interval
by Definition~\ref{GenericBoxes} and Remark~\ref{propiedadesR}(4).

Finally statement (c) follows from
Lemma~\ref{PC-properties}(a) and Definition~\ref{PCAtLast},
and (d) from Lemma~\ref{Propertiesvarphi}(b).
\end{proof}

\begin{lemma}\label{propiedadesPCA}
The following statements hold.
\begin{enumerate}[(a)]
\item $\setfibth{\A} \subset \setfibth{\basicbox{\ell}}$
      for every $\ell\in\Z$ and $\theta \in \basintabs{\ell}.$
\item $\setfibpt{\A}{\lstar} = \setfibpt{\A_{\all}}{\lstar}$
      for every $\ell\in\Z.$
      In particular $\setfibpt{\A}{\lstar}$ is a non-degenerate interval.
\item If $\theta\notin \Orbom,$ then
      $\setfibth{\A} = \{(\theta,\gamma(\theta))\}.$
\item $\A \subset \SI \times [-1,1].$
\end{enumerate}
\end{lemma}

\begin{proof}
Statement (c) follows directly from Lemma~\ref{PC-properties}(a).

Now we prove (a).
From Lemma~\ref{propiedadesA}(a),
$\setfibth{\A_n} \subset \basicbox{\ell}$
for every $\ell \in \Z$ and $n \ge \all.$
On the other hand,
by Definition~\ref{gammalimit} and Proposition~\ref{Hdimdinfty},
$\setfibth{\A} = \lim_{n\to\infty} \setfibth{\A_n}.$
Hence the result follows from the compacity of $\basicbox{\ell}$.

By Lemma~\ref{propiedadesA}(b)
and the part of the lemma already proved we have
\[
 \setfibpt{\A}{\lstar} = \lim_{n\to\infty} \setfibpt{\A_n}{\lstar} = \setfibpt{\A_{\all}}{\lstar}.
\]
Statement (d) follows from Lemma~\ref{propiedadesA}(d),
the compacity of $\SI \times [-1,1]$ and the fact that
$\A = \lim_{j\to\infty} \A_j.$
\end{proof}

The next proposition, summarizes the main properties of the set $\A.$

\begin{proposition}\label{teoremacentral}
The set $\A$ is a connected, does not contain any arc of curve and
$\Omega\setminus \A$ has two connected components.
\end{proposition}

\begin{proof}
From statements (b) and (c) of the previous lemma, we know that
$\setfibth{\A}$ is connected for every $\theta \in\SI.$

If $\A$ is not connected there exist closed (in $\A$) sets $U$ and $V$
such that $U \cap V = \emptyset$ and $U \cup V = \A$.
Observe that $\pi(U) \cup \pi(V) = \pi(\A) = \SI$
because every pseudo-curve is a circular set.
Moreover, since $\A$ is compact,
$U$ and $V$ are also compact sets of $\Omega$.
Hence, $\pi(U)$ are $\pi(V)$ compact in $\SI.$
Since {\SI} is connected, $\pi(U)\cap\pi(V) \ne \emptyset.$
For every $\theta \in \pi(U) \cap \pi(V)$ we have,
\[
  \setfibth{\A} = \setfibth{(U \cup V)} = \setfibth{U} \cup \setfibth{V}.
\]
The sets $\setfibth{U}$ and $\setfibth{V}$ are closed, non-empty and disjoint.
Consequently, $\setfibth{\A}$ is not connected; a contradiction.
This proves that $\A$ is connected.

By Lemma~\ref{propiedadesPCA}(b),
$\setfibpt{\A}{\lstar}$ is a non-degenerate interval for every $\ell\in \Orbom.$
Then, since $\Orbom$ is dense in $\SI,$
$\A$ does not contain any arc of curve by
Lemma~\ref{PC-properties-invariant}(b).

To prove that $\Omega\setminus \A$ has two connected components we define
\begin{align*}
\Omega_{-} &:= \set{(\theta,y) \in\Omega}{y < \min \set{x \in \I}{(\theta,x) \in \A}}\text{, and}\\
\Omega_{+} &:= \set{(\theta,y) \in\Omega}{y > \max \set{x \in \I}{(\theta,x) \in \A}}.
\end{align*}
By Lemma~\ref{propiedadesPCA}(d) we know that
\[
-1 \le \min \set{x \in \I}{(\theta,x) \in \A} \le \max \set{x \in \I}{(\theta,x) \in \A} \le 1.
\]
Hence,
$\Omega \setminus \A = \Omega_{-} \cup \Omega_{+},$
$\Omega_{+}$ and $\Omega_{-}$ are disjoint open circular subsets of $\Omega$
and $\Omega_{-} \supset \SI \times [-2, -1]$ and $\Omega_{+} \supset \SI \times [1, 2]$
(in particular, for every $\theta\in \SI,$ $\setfibth{\Omega_{+}}$ and $\setfibth{\Omega_{-}}$ are
non-degenerate intervals).
Thus, $\Omega_{+}$ and $\Omega_{-}$ are arc-wise connected and, hence, connected.
\end{proof}

\section{A collection of auxiliary functions $G_i$ defined on the boxes $\protect\wbasicbox{i}$}\label{FunctionsGi}

In this section we define a family of auxiliary functions
{\map{G_i}{\basicbox{i}}[\Omega]} with $i \in \Z$ and study their properties.

In  what follows we consider the supremum metric {\dinf} on the class of all functions
{\map{F}{A}[\Omega]} with $A \subset \Omega$. That is, given {\map{F,G}{A}[\Omega]}
we set
\[
 \dinf(F,G) := \leftlimits{\sup}{(\theta,x) \in A} \dom(F(\theta,x), G(\theta,x)).
\]
In the special case when $F$ and $G$ are skew products with the same base,
that is when $F(\theta,x) = (R(\theta),f(\theta,x))$ and
$G(\theta,x) = (R(\theta),g(\theta,x)),$ then
\[
 \dinf(F,G) :=
        \leftlimits{\sup}{(\theta,x) \in A}  \abs{f(\theta,x) - g(\theta,x)}.
\]
Observe that $(\cSO, \dinf)$ is a complete metric space.


Before defining the maps $G_i$ we need to introduce the necessary notation, and
recall and collect some basic facts that we will use in this definition
and to study their properties.

For every $i \in \Z,$ we define
\begin{align*}
& {\map{M_i}{\wbasint{i}}[\I]}
& \text{by}& \qquad
M_i(\theta) := \max \set{x \in \I}{(\theta,x) \in \wbasicbox{i}}, \text{ and}\\
& {\map{m_i}{\wbasint{i}}[\I]}
& \text{by} & \qquad
m_i(\theta) := \min \set{x \in \I}{(\theta,x) \in \wbasicbox{i}}.
\end{align*}

The next simple lemma states the basic properties of the maps $m_i$ and $M_i.$

\begin{lemma}\label{voresdelescaixesalesvores}
The following statements hold for every $i \in \Z$
\begin{enumerate}[(a)]
\item $-1 \le m_i(\theta) \le M_i(\theta) \le 1$ for every $\theta \in \wbasint{i}$.
\item $m_i $ and $M_i$ are continuous.
\item $m_i\evalat{\basintabs{i}}$ and $M_i\evalat{\basintabs{i}}$ are piecewise linear.
\item $m_i(\theta) = M_i(\theta) = \gams{\ai}(\theta)$ if and only if $\theta \in \wbasint{i} \setminus \obasintabs{i}.$
\end{enumerate}
\end{lemma}

\begin{proof}
It follows easily from Definition~\ref{GenericBoxes},
the definition of a winged region and
Lemma~\ref{Propertiesvarphi}(b,f).
\end{proof}

Notice that, for every $i \in \Z,$
\[
\wbasicbox{i} = \LSleftlimits{\bigcup}{\theta \in \wbasint{i}} \setfibth{\wbasicbox{i}}
             = \LSleftlimits{\bigcup}{\theta \in \wbasint{i}} \{\theta\} \times [m_i(\theta), M_i(\theta)].
\]

In what follows the interval $[m_i(\theta), M_i(\theta)] \subset \I,$
defined for every $\theta \in \wbasint{i},$
will be denoted by $\I_{i,\theta}.$
Clearly, for every $\theta \in \wbasint{i},$
$\setfibth{\wbasicbox{i}} = \{\theta\} \times \I_{i,\theta}.$

By Definition~\ref{PCgenerators}(\tsfR.2) and
Remark~\ref{PCgeneratorsExplicitConsequences}(\tsfR.2),
\[
\wbasint{i} \setminus \istarset
\andq[is disjoint from]
\Zstar_{\ai}.
\]

Hence, Lemmas~\ref{Propertiesvarphi}(a,d) and \ref{propiedadesA}(c)
can be summarized as:
\begin{equation}\label{gammathetaproperties}
\begin{cases}
\gams{\all}\evalat{\wbasint{\ell} \setminus \lstarset}
  \quad\text{is continuous,}\\
\gams{\all}(\theta) \in \I_{\ell,\theta}
  \andq[for every]
  \theta \in \wbasint{\ell} \setminus \lstarset,\text{ and}\\
\setfibth{\A_{\all}} = \{(\theta, \gams{\all}(\theta)\}
  \andq[for every]
  \theta \in \wbasint{\ell} \setminus \lstarset
\end{cases}
\end{equation}
for $\ell \in \{i, i+1\}.$

Now we define a family of continuous maps
{\map{G_i}{\wbasicbox{i}}[\Omega]} with $i\in \Z,$ by
\[
G_i(\theta,x) = \bigl(R_\omega(\theta), g_i(\theta,x)\bigr)
\]
Also, for every $\theta\in \wbasint{i},$
we will denote the map
{\map{g_i(\theta, \cdot)}{\I_{i,\theta}}[\I]}
by $g_{_{i,\theta}}.$

To define the functions $g_{_{i,\theta}},$ for clarity, we will consider separately two different situations:
\begin{itemize}
\item $i \ge 0,$ when $\wbasicbox{i} = \basicbox{i}$, $\wbasint{i} = \basintabs{i}$ and
      $G_i(\basicbox{i})$ strictly contains the smaller box $\basicbox{i+1},$ and
\item $i \le -1,$ when $G_i(\wbasicbox{i})$ is strictly contained in the bigger box $\basicbox{i+1}.$
\end{itemize}
We start by defining $g_{_{i,\theta}}$ for $i \ge 0$
in three different ways, depending on the base point $\theta \in \basint{i}$.
In this definition, for simplicity we will use
$\basicbox{i}$ instead of $\wbasicbox{i}$ and
$\basintabs{i}$ instead of $\wbasint{i}.$

Notice that, by Definition~\ref{PCgenerators}(\tsfR.1) and Lemma~\ref{Propertiesvarphi}(c),
\begin{equation}\label{alphadeltaintervals}
\begin{split}
& \text{for every $i \ge 0$}\\
& \BSG[\delta]{i}{i+1} \subset \OBG{i}{i+1}
  \andq
  \BSG{i}{i+1} \subset \obasint[\delta]{i} \subset \obasint{i},\text{ and}\\
& \gams{i-1}(\istar) = a_i
 \andq
 \gams{i}(\sstar{i+1}) = a_{i+1}.
\end{split}
\end{equation}
\begin{definition}[\bfseries Definition of $\boldsymbol{g_{i}}$ for $\boldsymbol{i\ge 0}$]\label{defi-gi-positiva}\strut
\begin{labeledlist}{$\boldsymbol{\theta \in \cball{\istar}{\alpha_{i+1}} \setminus \ball{\istar}{\delta_{i+1}}}$}
 \item [$\boldsymbol{\theta \in \cball{\istar}{\delta_{i+1}}}$]
       $
        g_{_{i,\theta}}(x) := \gams{i}(\sstar{i+1}) + \frac{2^{n_{i}}}{2^{n_{i+1}}}
          \left(\gams{i-1}(\istar) - x\right).
       $

 \item [$\boldsymbol{\theta \in \cball{\istar}{\alpha_{i+1}} \setminus \ball{\istar}{\delta_{i+1}}}$]
       we define $g_{_{i,\theta}}$ to be the unique piecewise affine map
       with two affine pieces, defined on $\I_{i,\theta},$
       whose graph joins
       $\left(m_{i}(\theta), M_{i+1}\bigl(R_\omega(\theta)\bigr)\right)$
       with
       $\left(\gams{i}(\theta),\gams{i+1}\bigl(R_\omega(\theta)\bigr)\right),$
       and this with the point
       $\left(M_{i}(\theta), m_{i+1}\bigl(R_\omega(\theta)\bigr)\right)$
       (in particular,
        $g_{_{i,\theta}}\bigl(\gams{i}(\theta)\bigr)= \gams{i+1}\bigl(R_\omega(\theta)\bigr)$),

 \item [$\boldsymbol{\theta \in \cball{\istar}{\alpha_{i}} \setminus \ball{\istar}{\alpha_{i+1}}}$]
       $g_{_{i,\theta}}(x) := \gams{i+1}\bigl(R_\omega(\theta)\bigr)$
       (that is, $g_{_{i,\theta}}$ is constant).
\end{labeledlist}
\end{definition}

The next lemma states the basic properties of the functions $G_i$ for $i \ge 0.$

\begin{lemma}\label{gpositiva}
The following statements hold for every $i \ge 0:$
\begin{enumerate}[(a)]
\item The map $g_{_{i,\theta}}$ is well defined and non-increasing for every $\theta \in \basint{i}.$
      Moreover, $-1 \le g_{_{i,\theta}}(x) \le 1$ for every $\theta \in \basint{i}$ and $x \in \I_{i,\theta}.$
      Furthermore, the function $G_i$ is continuous.
\item $G_i\evalat{\setfibth{\basicbox{i}}}$ is affine and
      $G_i\bigl(\setfibth{\basicbox{i}}\bigr) = \setfibpt{\basicbox{i+1}}{R_\omega(\theta)}$
      for every $\theta \in \BSG[\delta]{i}{i+1};$
      $G_i\evalat{\setfibth{\basicbox{i}}}$ is piecewise affine with two pieces and
      $G_i\bigl(\setfibth{\basicbox{i}}\bigr) = \setfibpt{\basicbox{i+1}}{R_\omega(\theta)}$
      for every $\theta \in \BSG{i}{i+1} \setminus \OBG[\delta]{i}{i+1};$ and\\
      $G_i\bigl(\setfibth{\basicbox{i}}\bigr) = \setfibpt{\A_{i+1}}{R_\omega(\theta)}$
      for every
      $\theta \in \basint{i} \setminus \OBG{i}{i+1}.$
\item $G_i(\setfibth{\A_{i}}) = \setfibpt{\A_{i+1}}{R_\omega(\theta)}$ for every $\theta \in \basint{i}.$
\end{enumerate}
\end{lemma}

\begin{proof}
We will prove all statements of the lemma simultaneously and according to the regions
in the definition of the map $g_i.$

\inidemopartfree{\textbullet} We start with the region $\setfibball[\delta]{i}{i+1}.$\\[\medskipamount]
Let $z \in [-\delta_{i}, \delta_{i}] \subset \R$ and let
$\theta = \istar + z \in \basint[\delta]{i}.$
From Definition~\ref{GenericBoxes} and \eqref{alphadeltaintervals} we get
\begin{equation}\label{mMformula}
\begin{split}
m_i(\theta) &= a_i - 2^{-{n_i}}(1-z) =
               \gams{i-1}(\istar) - 2^{-{n_i}}(1-z), \text{ and}\\
M_i(\theta) &= a_i + 2^{-{n_i}}(1-z) =
               \gams{i-1}(\istar) + 2^{-{n_i}}(1-z).
\end{split}
\end{equation}
In a similar way, for every $\theta \in \BSG[\delta]{i}{i+1}$
(that is, $z \in [-\delta_{i+1}, \delta_{i+1}]$), we have
$R_\omega(\theta) = \sstar{i+1} + z \in \basint[\delta]{i+1},$
and
\begin{equation}\label{mMRotformula}
\begin{split}
m_{i+1}(R_\omega(\theta)) &= a_{i+1} - 2^{-{n_{i+1}}}(1-z) =
               \gams{i}(\sstar{i+1}) - 2^{-{n_{i+1}}}(1-z), \text{ and}\\
M_{i+1}(R_\omega(\theta)) &= a_{i+1} + 2^{-{n_{i+1}}}(1-z) =
               \gams{i}(\sstar{i+1}) + 2^{-{n_{i+1}}}(1-z).
\end{split}
\end{equation}
Hence, for every $\theta \in \BSG[\delta]{i}{i+1},$
\begin{equation}\label{mappingendpoints}
\begin{split}
g_{_{i,\theta}}(m_i(\theta))
    &= \gams{i}(\sstar{i+1}) +
       \tfrac{2^{n_{i}}}{2^{n_{i+1}}} 2^{-{n_i}} (1-z)
     = \gams{i}(\sstar{i+1}) + 2^{-{n_{i+1}}} (1-z)\\
    &=  M_{i+1}(R_\omega(\theta)),\\
g_{_{i,\theta}}(M_i(\theta))
    &= \gams{i}(\sstar{i+1}) -
       \tfrac{2^{n_{i}}}{2^{n_{i+1}}} 2^{-{n_i}} (1-z)
     = \gams{i}(\sstar{i+1}) - 2^{-{n_{i+1}}} (1-z)\\
    &=  m_{i+1}(R_\omega(\theta)).
\end{split}
\end{equation}
So,
$g_{_{i,\theta}}\evalat{\I_{i, \theta}}$
is the affine map whose graph joins the point
$\bigl(m_{i}(\theta), M_{i+1}(R_\omega(\theta))\bigr)$ with
$\bigl(M_{i}(\theta), m_{i+1}(R_\omega(\theta))\bigr).$
In particular, $g_{_{i,\theta}}$ sends the interval
$\I_{i, \theta}$ affinely onto
$\I_{i+1, R_\omega(\theta)}$ or, equivalently,
$G_i$ sends the interval $\setfibth{\basicbox{i}}$ affinely onto
$\setfibpt{\basicbox{i+1}}{R_\omega(\theta)}.$
Then, by Lemma~\ref{Propertiesvarphi}(b), this implies that
$-1 \le g_{_{i,\theta}}(x) \le 1$ for every $x \in \I_{i,\theta}.$
Moreover, the continuity of the maps $m_i,\ M_i,\ m_{i+1} \circ R_\omega$
and $M_{i+1} \circ R_\omega$ imply that $g_{i}$ is well defined
and continuous on
$\setfibball[\delta]{i}{i+1}$

Next we will prove that
$
G_i\bigl(\setfibth{\A_{i}}\bigr) = \setfibpt{\A_{i+1}}{R_\omega(\theta)}
$
for every $\theta \in \BSG[\delta]{i}{i+1}.$
We take
$\theta = \istar + z \in \BSG[\delta]{i}{i+1}\setminus\istarset.$
Then, clearly, $z \in [-\delta_{i+1}, \delta_{i+1}]\setminus\{0\} \subset \R.$
By Definitions~\ref{PCgenerators} and \ref{GenericBoxes} and statement \eqref{alphadeltaintervals},
\begin{align*}
\gams{i}(\theta)
&= \varphi_{_{\istar}}(\theta) =
     a_i + 2^{-n_i} d =
     \gams{i-1}(\istar) + 2^{-n_i} d \in \I_{i, \theta}, \text{ and}\\
\gams{i+1}(R_\omega(\theta))
&= \varphi_{_{\sstar{i+1}}}(\theta) =
     a_{i+1} - 2^{-n_{i+1}} d =
     \gams{i-1}(\istar) - 2^{-n_{i+1}} d \in \I_{i+1, R_\omega(\theta)},
\end{align*}
where $d = (-1)^{i} \phi(z).$
So, for every
$\theta \in \BSG[\delta]{i}{i+1}\setminus\istarset,$
\begin{equation}\label{Formula1}
g_{_{i,\theta}}(\gams{i}(\theta))
= \gams{i}(\sstar{i+1}) - \frac{2^{n_{i}}}{2^{n_{i+1}}} 2^{-n_i} d
= \gams{i+1}(R_\omega(\theta)).
\end{equation}
Thus, from \eqref{alphadeltaintervals} and \eqref{gammathetaproperties} we get
\begin{align*}
G_i\bigl(\setfibth{\A_{i}}\bigr)
&= G_i\bigl(\{(\theta, \gams{i}(\theta))\}\bigr)
=\{(R_\omega(\theta), g_{_{i,\theta}}(\gams{i}(\theta)))\}\\
&= \{(R_\omega(\theta), \gams{i+1}(R_\omega(\theta)))\}
= \setfibpt{\A_{i+1}}{R_\omega(\theta)}
\end{align*}
for every
$\theta \in \BSG[\delta]{i}{i+1}\setminus\istarset.$
On the other hand, by the part already proven,
$g_{_{i,\istar}}$ sends the interval
$\I_{i,\istar}$ affinely to
$\I_{i+1,\sstar{i+1}}$
or, equivalently, $G_i$ sends the interval
$
  \setfibbb{i} = \istarset \times \I_{i,\istar}
$
affinely onto
$
  \setfibbb{i+1} = \iistarset \times \I_{i,\sstar{i+1}}.
$
This implies that
$
G_i\bigl(\setfibpt{\A_{i}}{\istar}\bigr) = \setfibpt{\A_{i+1}}{\sstar{i+1}}
$
by Lemma~\ref{propiedadesA}(b).
Hence,
$
G_i\bigl(\setfibth{\A_{i}}\bigr) = \setfibpt{\A_{i+1}}{R_\omega(\theta)}
$
for every $\theta \in \BSG[\delta]{i}{i+1}.$

\inidemopartfree{\textbullet} Now we study
$
\setfib{\basicbox{i}}{\left(\BSG{i}{i+1} \setminus
    \OBG[\delta]{i}{i+1}\right)}.
$\\[\medskipamount]
Observe that
$
R_\omega(\cball{\istar}{\alpha} \setminus \istarset) =
   \cball{\sstar{i+1}}{\alpha} \setminus \iistarset
$
for $\alpha \in \{\alpha_i, \alpha_{i+1}\}.$
Then, by \eqref{gammathetaproperties}
\begin{equation}\label{gammaRotthetaontoca}
\begin{split}
& \gams{i+1} \circ R_\omega\evalat{\basint{i} \setminus \istarset}
  \quad\text{is continuous, and}\\
& \gams{i+1}(R_\omega(\theta)) \in \I_{i+1, R_\omega(\theta)}
  \andq[for every]
  \theta \in \BSG{i}{i+1} \setminus \istarset.
\end{split}
\end{equation}
So, the continuity of the maps
$m_i,\ M_i,\ m_{i+1} \circ R_\omega$ and $M_{i+1} \circ R_\omega$
imply that $g_{i}$ is well defined and continuous on
$
\setfib{\basicbox{i}}{\left(\BSG{i}{i+1} \setminus
    \OBG[\delta]{i}{i+1}\right)},
$
and
\[
 \bigl(\gams{i}(\theta),\gams{i+1}(R_\omega(\theta))\bigr) \in
 \I_{i, \theta}
 \times
 \I_{i+1, R_\omega(\theta)}
\]
for every
$\theta \in \BSG{i}{i+1} \setminus \OBG[\delta]{i}{i+1}.$
Consequently, $g_{_{i,\theta}}$ maps
$\I_{i, \theta}$ piecewise affinely with two pieces onto
$\I_{i+1, R_\omega(\theta)}$ or, equivalently,
$G_i$ sends the interval $\setfibth{\basicbox{i}}$
piecewise affinely with two pieces onto
$\setfibpt{\basicbox{i+1}}{R_\omega(\theta)}.$
Again, by Lemma~\ref{Propertiesvarphi}(b), this implies that
$-1 \le g_{_{i,\theta}}(x) \le 1$ for every $x \in \I_{i,\theta}.$
On the other hand,
from \eqref{alphadeltaintervals} and \eqref{gammathetaproperties} we have
\begin{align*}
G_i\bigl(\setfibth{\A_{i}}\bigr)
&= G_i\bigl(\{(\theta, \gams{i}(\theta))\}\bigr)
=\{(R_\omega(\theta), g_{_{i,\theta}}(\gams{i}(\theta)))\}\\
&= \{(R_\omega(\theta), \gams{i+1}(R_\omega(\theta)))\}
= \setfibpt{\A_{i+1}}{R_\omega(\theta)}
\end{align*}
for every $\theta \in \BSG{i}{i+1} \setminus \OBG[\delta]{i}{i+1}.$

\inidemopartfree{\textbullet} Finally, we study the region
$
\setfib{\basicbox{i}}{\left(\basint{i} \setminus
   \OBG{i}{i+1}\right)}.
$\\[\medskipamount]
In this case, by definition and Lemma~\ref{Propertiesvarphi}(b) we have
$-1 \le g_{_{i,\theta}}(x) \le 1$ for every $x \in \I_{i,\theta}.$
By \eqref{gammaRotthetaontoca},
$g_{i}(\cdot, x) = \gams{i+1} \circ R_\omega$ is well defined and continuous in
both variables on
$
\setfib{\basicbox{i}}{\left(\basint{i} \setminus
   \OBG{i}{i+1}\right)}
$
because $m_i$ and $M_i$ are continuous.
Moreover, for every
$\theta \in \basint{i} \setminus \OBG{i}{i+1}$
and $x$ such that $(\theta,x) \in \setfibth{\basicbox{i}},$
we have
\[
\{G_i(\theta, x)\}
 = \{(R_\omega(\theta), g_{i}(\theta, x))\} =
   \{(R_\omega(\theta), \gams{i+1}(R_\omega(\theta))\} =
   \setfibpt{\A_{i+1}}{R_\omega(\theta)}
\]
by Definition~\ref{PCAtLast} and Lemma~\ref{PC-properties}(a).
Thus, by Lemma~\ref{propiedadesA}(a),
\[
 G_i\bigl(\setfibth{\A_{i}}\bigr) =
 G_i\bigl(\setfibth{\basicbox{i}}\bigr) =
 \setfibpt{\A_{i+1}}{R_\omega(\theta)}.
\]

From all the previous arguments (b) and (c) follow.
To end the proof of (a) we have to see that $G_i$ is well defined and globally continuous.
This amounts to show that it is well defined on the fibres
\begin{align*}
\setfibpt{\basicbox{i}}{(\istar \pm \delta_{i+1})}
&= \{\istar \pm \delta_{i+1}\} \times \I_{i, \istar \pm \delta_{i+1}} \text{ and}\\
\setfibpt{\basicbox{i}}{(\istar \pm \alpha_{i+1})}
&= \{\istar \pm \alpha_{i+1}\} \times \I_{i, \istar \pm \alpha_{i+1}}.
\end{align*}
We will only show that the two definitions of $g_{i}$ coincide on
$\{\theta\} \times \I_{i, \theta}$
with $\theta \in \{\istar + \delta_{i+1},  \istar + \alpha_{i+1}\}.$
The case $\theta \in \{\istar - \delta_{i+1},  \istar - \alpha_{i+1}\}$
follows analogously.

We start with $\theta = \istar + \alpha_{i+1} \in \obasint[\delta]{i}.$
In this case,
$R_\omega(\theta) = \sstar{i+1} + \alpha_{i+1} \in \Bd(\basint{i+1})$ and,
by Definition~\ref{GenericBoxes} and Lemma~\ref{Propertiesvarphi}(c),
\[
M_{i+1}(R_\omega(\theta)) =
m_{i+1}(R_\omega(\theta)) = a^+_{i+1} =
\gams{i+1}(R_\omega(\theta)).
\]
Thus, the piecewise affine map whose graph joins the points
$\bigl(m_{i}(\theta), M_{i+1}(R_\omega(\theta))\bigr),$
$\bigl(\gams{i}(\theta),\gams{i+1}(R_\omega(\theta))\bigr),$ and
$\bigl(M_{i}(\theta), m_{i+1}(R_\omega(\theta))\bigr)$
is the constant map $\gams{i+1}(R_\omega(\theta)).$
Hence, $g_{_{i,\theta}}$ is well defined
for $\theta = \istar + \alpha_{i+1}.$

Now we deal with the case
$\theta = \istar + \delta_{i+1} \in \basint[\delta]{i}.$
By \eqref{mappingendpoints} and \eqref{Formula1} we know that the points
$\bigl(m_{i}(\theta), M_{i+1}(R_\omega(\theta))\bigr),$
$\bigl(\gams{i}(\theta),\gams{i+1}(R_\omega(\theta))\bigr)$ and
$\bigl(M_{i}(\theta), m_{i+1}(R_\omega(\theta))\bigr)$
belong to
$
\Graph\left(x \mapsto \gams{i}(\sstar{i+1}) +
     \frac{2^{n_{i}}}{2^{n_{i+1}}} \left(\gams{i-1}(\istar) - x\right)
\right).
$
Consequently, the map
$
\gams{i}(\sstar{i+1}) +
\frac{2^{n_{i}}}{2^{n_{i+1}}} \left(\gams{i-1}(\istar) - x\right)
$
coincides with the piecewise affine map whose graph joins
$\bigl(m_{i}(\theta), M_{i+1}(R_\omega(\theta))\bigr),$
$\bigl(\gams{i}(\theta),\gams{i+1}(R_\omega(\theta))\bigr)$ and
$\bigl(M_{i}(\theta), m_{i+1}(R_\omega(\theta))\bigr).$
This ends the proof of (a).
\end{proof}

Now we define $g_{_{i,\theta}}$ for $i < 0.$
In this case, since we are going from a smaller box $\wbasicbox{i}$
to a bigger one, we only need to define $g_{_{i,\theta}}$ in two
different ways, depending on the base point $\theta \in \wbasint{i}$.

As in the previous case we need to fix some facts about the elements
that we will use in the definition.

By Definition~\ref{PCgenerators}(\tsfR.1) and Lemma~\ref{Propertiesvarphi}(c),
\begin{equation}\label{alphadeltaintervalsineg}
\begin{split}
   & \text{for every $i < 0$\hspace*{1.5em}}\\
   & \BSG[\delta]{i+1}{\ai} \subset\
      \BSG{i+1}{\ai} \subset
      \OBG[\delta]{i+1}{\aii} \subset
      \obasintabs{i+1},\\
   & R_\omega\left(\wbasint{i}\right) = \basintabs{i+1},\quad
     \basintabs[\delta]{i} \subset \obasintabs{i}, \text{ and}\\
& \gams{\aii}(\istar) = a_i \quad
       \text{and}\quad
       \gams{\abs{i+2}}(\sstar{i+1}) = a_{i+1}.
\end{split}
\end{equation}
Consequently, from \eqref{gammathetaproperties}
and Definitions~\ref{GenericBoxes} and \ref{PCgenerators} we get
\begin{align*}
& m_{i}(\theta) < \gams{\ai}(\theta) < M_{i}(\theta) \text{ and}\\
& m_{i+1}(R_\omega(\theta)) < \gams{\aii}(R_\omega(\theta)) < M_{i+1}(R_\omega(\theta))
\end{align*}
for every $\theta \in \obasintabs{i} \setminus \istarset$
(and $R_\omega(\theta) \in \OBG{i+1}{\ai} \setminus \iistarset$).
Then,
\[
\widetilde{\kappa}_i(\theta) = \min \left\{1,
   \frac{m_{i+1}(R_\omega(\theta)) - \gams{\aii}(R_\omega(\theta))
   }{\tfrac{2^{n_{\ai}}}{2^{n_{_{\aii}}}} (\gams{\ai} (\theta) - M_{i}(\theta))},
   \frac{ M_{i+1}(R_\omega(\theta)) - \gams{\aii}(R_\omega(\theta))
   }{\tfrac{2^{n_{\ai}}}{2^{n_{_{\aii}}}} (\gams{\ai}(\theta) - m_{i}(\theta))}
\right\} > 0
\]
defines a continuous function
$
\map{\widetilde{\kappa}_i}{\obasintabs{i} \setminus \obasintabs[\delta]{i}}[{(0,1]}].
$
To define the map $g_{i}$ we need an auxiliary function
\[
\map{\kappa_{i}}{\basintabs{i} \setminus \obasintabs[\delta]{i}}[{[0,1]}]
\]
such that
$\kappa_{i}$ is non-decreasing and continuous,
$\kappa_{i}(\istar\pm \delta_{\ai}) = \widetilde{\kappa}_{i}(\istar\pm \delta_{\ai}),$ and
$\kappa_{i}(\theta) \le \widetilde{\kappa}_{i}(\theta)$ for every
$\theta \in \obasintabs{i} \setminus \obasintabs[\delta]{i}.$
In principle any such function would do, but for definiteness, and to show that such function exists,
we note that we can take, for instance,
\[
 \kappa_{i}(\theta) = \begin{cases}
   \inf_{t \in [\theta, \istar - \delta_{\ai}] \cap \obasintabs{i}} \widetilde{\kappa}_i(t)
      & \text{if $\theta \le \istar - \delta_{\ai}$},\\
   \inf_{t \in [\istar + \delta_{\ai}, \theta] \cap \obasintabs{i}} \widetilde{\kappa}_i(t)
      & \text{if $\theta \ge \istar + \delta_{\ai}$}.
\end{cases}
\]
It is easy to check that this map verifies the desired properties.

\begin{definition}[\bfseries Definition of $\boldsymbol{g_{i}}$ for $\boldsymbol{i< 0}$]\label{defi-gi-negativa}
For every $(\theta,x) \in \wbasicbox{i}$ we set
\[
g_{_{i,\theta}}(x) := \begin{cases}
  \frac{2^{n_{\ai}}}{2^{n_{_{\aii}}}}
         \left(\gams{\aii}(\istar) - x\right)
      + \gams{\abs{i+2}}(\sstar{i+1})
  & \text{if $\theta \in \basintabs[\delta]{i},$}\\
  \frac{2^{n_{\ai}}}{2^{n_{_{\aii}}}} \kappa_{i}(\theta)
         \left(\gams{\ai}(\theta) - x\right)
      + \gams{\aii}(R_\omega(\theta))
  &\text{if $\theta \in \basintabs{i} \setminus \obasintabs[\delta]{i}$}\\
  \gams{\aii}(R_\omega(\theta))
  &\text{if $\theta \in \wbasint{i} \setminus \obasintabs{i}.$}
\end{cases}
\]
\end{definition}

The next lemma states the basic properties of the functions $G_i$ for $i < 0.$

\begin{lemma}\label{gnegativa}
The following statements hold for every $i < 0:$
\begin{enumerate}[(a)]
\item The map $g_{_{i,\theta}}$ is well defined and non-increasing for every $\theta \in \basint{i}.$
      Moreover, $-1 \le g_{_{i,\theta}}(x) \le 1$ for every $\theta \in \basint{i}$ and $x \in \I_{i,\theta}.$
      Furthermore, the function $G_i$ is continuous.
\item $G_i\evalat{\setfibth{\wbasicbox{i}}}$ is affine,
      $G_i\bigl(\setfibth{\wbasicbox{i}}\bigr) \subset \setfibpt{\basicbox{i+1}}{R_\omega(\theta)}$
      for every $\theta \in \wbasint{i}$ and
      $G_i\bigl(\setfibth{\wbasicbox{i}}\bigr) = \setfibpt{\basicbox{i+1}}{R_\omega(\theta)}$
      for every $\theta \in \basintabs[\delta]{i}.$
\item $G_i(\setfibth{\A_{\ai}}) = \setfibpt{\A_{\aii}}{R_\omega(\theta)}$ for every $\theta \in \wbasint{i}.$
\end{enumerate}
\end{lemma}

\begin{proof}
First we will prove that the map $G_i$ is continuous and that
$G_i\evalat{\setfibth{\wbasicbox{i}}}$ is affine, according to the
three regions in the definition.

\inidemopartfree{\textbullet}
As in the previous lemma we start with
$\setfib{\wbasicbox{i}}{\basintabs[\delta]{i}} =
 \setfibball[\delta]{i}{\ai}.$\\[\medskipamount]
As in the same case of Lemma~\ref{gpositiva}, by using
\eqref{alphadeltaintervalsineg} instead of \eqref{alphadeltaintervals},
it follows that
$g_{_{i,\theta}}\evalat{\I_{i, \theta}}$
is the affine map whose graph joins the points
$\bigl(m_{i}(\theta), M_{i+1}(R_\omega(\theta))\bigr)$ and
$\bigl(M_{i}(\theta), m_{i+1}(R_\omega(\theta))\bigr),$
$g_{i}$ is well defined and continuous on
$\setfibball[\delta]{i}{\ai},$
\begin{align*}
& g_{_{i,\theta}}(\gams{\ai}(\theta))
   = \gams{\aii}(R_\omega(\theta))
      \text{ for every $\theta \in \basintabs[\delta]{i}\setminus\istarset,$}\\
& \text{$G_i$ sends the interval $\setfibth{\basicbox{i}}$ affinely onto
      $\setfibth{\basicbox{i+1}}$, and}\\
& G_i\bigl(\setfibth{\A_{\ai}}\bigr)
   = \setfibpt{\A_{\aii}}{R_\omega(\theta)}\quad
      \text{for every $\theta \in \basintabs[\delta]{i}$}.
\end{align*}

\inidemopartfree{\textbullet}
$
\setfib{\wbasicbox{i}}{\bigl(\basintabs{i} \setminus \obasintabs[\delta]{i}\bigr)} =
\setfib{\basicbox{i}}{\bigl(\basintabs{i} \setminus \obasintabs[\delta]{i}\bigr)}.
$\\[\medskipamount]
From \eqref{gammathetaproperties} we know that the maps
$\gams{\ai}$ and $\gams{\aii} \circ R_\omega$
are continuous on the domain
$\basintabs{i} \setminus \obasintabs[\delta]{i}.$
Hence, the continuity of $g_{i}$ follows from the continuity of the maps
$\kappa_i$, $m_i,\ M_i,\ m_{i+1} \circ R_\omega$ and $M_{i+1} \circ R_\omega.$

Notice that, from the definition of $g_{i}$ in this region we clearly have that
\[
 g_{_{i,\theta}}(\gams{\ai}(\theta)) = \gams{\aii}(R_\omega(\theta)), \text{ and}
\]
$G_i\evalat{\setfibth{\wbasicbox{i}}} = g_i(\theta, \cdot)$ is affine.

\inidemopartfree{\textbullet}
$
\setfib{\wbasicbox{i}}{\bigl(\wbasint{i} \setminus \obasintabs{i}\bigr)}.
$\\[\medskipamount]
In this case we have $m_{i}(\theta) = \gams{\ai}(\theta) = M_{i}(\theta)$
by definition. Then, the map
$G_i\evalat{\setfibth{\wbasicbox{i}}} = g_i(\theta, \cdot)$ is affine because
it is constant, and $g_i$ is continuous because
$\gams{\ai}$ and $\gams{\aii} \circ R_\omega$
are continuous on the domain $\wbasint{i} \setminus \istarset$
by \eqref{gammathetaproperties}.

To end the proof of (a) we have to see that $G_i$ is well defined
and globally continuous.
This amounts to show that it is well defined on the fibres
\[
\setfibpt{\basicbox{i}}{(\istar \pm \delta_{\ai})}
\andq
\setfibpt{\basicbox{i}}{(\istar \pm \alpha_{\ai})}
\]
We start by showing that the two definitions of $g_{i}$ coincide on
the fibres
$\setfibth{\basicbox{i}}$ for $\theta \in \{\istar \pm \alpha_{\ai}\}.$
In this case we have $m_{i}(\theta) = \gams{\ai}(\theta) = M_{i}(\theta).$
Consequently, $\I_{i, \theta} = \{\gams{\ai}(\theta)\}$ and
\[
  \frac{2^{n_{\ai}}}{2^{n_{_{\aii}}}} \kappa_{i}(\theta)
         \left(\gams{\ai}(\theta) - x\right)
      + \gams{\aii}(R_\omega(\theta)) =
  \gams{\aii}(R_\omega(\theta))
\]
for $x \in \I_{i, \theta}.$

Next we consider
$\setfibth{\basicbox{i}} = \{\theta\} \times \I_{i,\theta}$ with
$\theta = \istar + \delta_{\ai}.$
We will show that the two definitions of $g_{i}$ coincide on this set.
The case $\theta = \istar - \delta_{\ai}$ follows analogously.

For simplicity we will denote
\begin{align*}
g^{\delta_{\ai}}_{_{i,\theta}}(x) &:= \frac{2^{n_{\ai}}}{2^{n_{_{\aii}}}}
             \left(\gams{\aii}(\istar) - x\right)
        + \gams{\abs{i+2}}(\sstar{i+1}),\text{ and}\\
\xi_{_{i,\theta}}(x) &:= \frac{2^{n_{\ai}}}{2^{n_{_{\aii}}}}
         \left(\gams{\ai}(\theta) - x\right)
      + \gams{\aii}(R_\omega(\theta)).
\end{align*}
Notice that $g^{\delta_{\ai}}_{_{i,\theta}}$ is the map
$g_{_{i,\theta}}$ as defined in the first region while
\[
\kappa_i(\theta)\left( \xi_{_{i,\theta}} - \gams{\aii}(R_\omega(\theta)) \right)
      + \gams{\aii}(R_\omega(\theta))
\]
is the map  $g_{_{i,\theta}}$ as defined in the second region.
In a similar way to the previous lemma we have that
$
  \bigl(\gams{\ai}(\theta),\gams{\aii}(R_\omega(\theta))\bigr)
  \in
  \Graph\left(g^{\delta_{\ai}}_{_{i,\theta}}\right).
$
Hence, since $g^{\delta_{\ai}}_{_{i,\theta}}$ is affine with slope
$-\tfrac{2^{n_{\ai}}}{2^{n_{_{\aii}}}},$
it follows that
$g^{\delta_{\ai}}_{_{i,\theta}} = \xi_{_{i,\theta}}.$
So, to end the proof of the lemma, we only have to see that
$
  \kappa_{i}(\istar + \delta_{\ai}) =
  \widetilde{\kappa}_{i}(\istar + \delta_{\ai}) = 1.
$

Since the points
$\bigl(m_{i}(\theta), M_{i+1}(R_\omega(\theta))\bigr)$ and
$\bigl(M_{i}(\theta), m_{i+1}(R_\omega(\theta))\bigr)$
also belong to
$
\Graph\left(g^{\delta_{\ai}}_{_{i,\theta}}\right) =
\Graph\left(\xi_{_{i,\theta}}\right),
$
it follows that
\begin{align*}
m_{i+1}(R_\omega(\theta))
  &= \xi_{_{i,\theta}}(M_{i}(\theta))
   = \frac{2^{n_{\ai}}}{2^{n_{_{\aii}}}}
         \left(\gams{\ai}(\theta) - M_{i}(\theta)\right)
      + \gams{\aii}(R_\omega(\theta)), \text{ and}\\
M_{i+1}(R_\omega(\theta))
  &= \xi_{_{i,\theta}}(m_{i}(\theta))
   = \frac{2^{n_{\ai}}}{2^{n_{_{\aii}}}}
         \left(\gams{\ai}(\theta) - m_{i}(\theta)\right)
      + \gams{\aii}(R_\omega(\theta)).
\end{align*}
This shows that
$
  \widetilde{\kappa}_{i}(\istar + \delta_{\ai}) =
  \widetilde{\kappa}_{i}(\theta) = 1
$
and ends the proof of (a).

Now we prove (b) according to the three regions in the definition.
From the part of the lemma already proven we already know that
$G_i\evalat{\setfibth{\wbasicbox{i}}}$ is affine, and
$
  G_i\bigl(\setfibth{\wbasicbox{i}}\bigr) =
  \setfibpt{\basicbox{i+1}}{R_\omega(\theta)}
$
for every $\theta \in \basintabs[\delta]{i}.$
So, to end the proof of (b) we have to see that
\begin{equation}\label{fibinclu}
g_{_{i,\theta}}(\I_{i, \theta}) \subset \I_{i+1, R_\omega(\theta)}
\end{equation}
for every $\theta \in \wbasint{i} \setminus \basintabs[\delta]{i}$
(by definition, since $i < 0,$ $\wbasint{i} = \BSG{i}{\aii};$
therefore,
$R_\omega(\theta) \in \basintabs{i+1}$ and
$\I_{i+1, R_\omega(\theta)} = \setfibpt{\basicbox{i+1}}{R_\omega(\theta)}$).

For $\theta \in \wbasint{i} \setminus \obasintabs{i},$
by \eqref{gammathetaproperties}, we have
\[
  g_{_{i,\theta}}(\I_{i, \theta})
   = \{\gams{\aii}(R_\omega(\theta))\}
   \subset \I_{i+1, R_\omega(\theta)}.
\]
Now we consider
$\theta \in \obasintabs{i} \setminus \basintabs[\delta]{i}.$
Since
\[
\kappa_{i}(\theta) \le \widetilde{\kappa}_{i}(\theta) \le
 \frac{ M_{i+1}(R_\omega(\theta))-\gams{\aii}(R_\omega(\theta))
   }{ \tfrac{2^{n_{\ai}}}{2^{n_{_{\aii}}}} (\gams{\ai}(\theta)-m_{i}(\theta))},
\]
we have
\begin{align*}
g_{_{i,\theta}}(m_{i}(\theta))
 &\le \frac{2^{n_{\ai}}}{2^{n_{_{\aii}}}}
       \frac{M_{i+1}(R_\omega(\theta))-\gams{\aii}(R_\omega(\theta))
             }{\tfrac{2^{n_{\ai}}}{2^{n_{_{\aii}}}} (\gams{\ai}(\theta)-m_{i}(\theta))}
       \left(\gams{\ai}(\theta) - m_{i}(\theta)\right)
          + \gams{\aii}(R_\omega(\theta))\\
  &= M_{i+1}(R_\omega(\theta)).
\end{align*}
An analogous computation shows that
$g_{_{i,\theta}}(M_{i}(\theta)) \ge m_{i+1}(R_\omega(\theta)).$
Hence, \eqref{fibinclu} holds because $g_{_{i,\theta}}$ is affine.
This ends the proof of (b).

Then, by Lemma~\ref{Propertiesvarphi}(b), Statement~(b) of the lemma
implies that $-1 \le g_{_{i,\theta}}(x) \le 1$ for every $x \in \I_{i,\theta}.$

By the part of the lemma already proved we know that
$
  G_i\bigl(\setfibth{\A_{\ai}}\bigr)
     = \setfibpt{\A_{\aii}}{R_\omega(\theta)}
$
for every $\theta \in \basintabs[\delta]{i}$.
On the other hand, as in the previous lemma,
from \eqref{alphadeltaintervalsineg} and \eqref{gammathetaproperties}
we get
\begin{align*}
G_i\bigl(\setfibth{\A_{\ai}}\bigr)
&= G_i\bigl(\{(\theta, \gams{\ai}(\theta))\}\bigr)
 =\{(R_\omega(\theta), g_{_{i,\theta}}(\gams{\ai}(\theta)))\}\\
&= \{(R_\omega(\theta), \gams{\aii}(R_\omega(\theta)))\}
 = \setfibpt{\A_{\aii}}{R_\omega(\theta)}
\end{align*}
for every $\theta \in \wbasint{i} \setminus \basintabs[\delta]{i}.$
So, (c) holds.
\end{proof}

Up to now we have defined the family of auxiliary functions
{\map{G_i}{\wbasicbox{i}}[\Omega]} with $i \in \Z.$
The next step before being able to define
the family $\{T_m\} \subset \cSO$ is to fix some
stratification in the set of boxes $\wbasicbox{i}.$

\section{A stratification in the set of boxes $\protect\wbasicbox{i}$}\label{stratification}

In this section we introduce a notion of \emph{depth}
in the set of arcs $\wbasint{i}$ defined earlier.
This notion introduce a stratification in the set of boxes
$\wbasicbox{i}$ that we study below.

\begin{definition}\label{depth}
For every $\ell\in\Z$ we define the \emph{depth of $\ell$},
which will be denoted by $\dep(\ell)$,
as the cardinality of the set (see Lemma~\ref{Propertiesvarphi}(g))
\begin{align*}
 \set{i\in\Z}{\wbasint{\ell} \varsubsetneq \wbasint{i}} &=
 \set{i\in\Z}{\wbasint{\ell} \cap \wbasint{i} \ne \emptyset} =\\
 \set{i\in\Z}{\wbasicbox{\ell} \varsubsetneq \wbasicbox{i}} &=
 \set{i\in\Z}{\wbasicbox{\ell} \cap \wbasicbox{i} \ne \emptyset}.
\end{align*}

Also, for every $m \in \Z^+,$ we denote
\begin{align*}
  \DS &:= \set{\ell \in \Z}{\dep(\ell) = m},\\
  \sstar{\DS} &:= \set{\istar}{i\in\DS},\text{ and}\\
  \mu_m &:= \min\set{\ai}{i \in \DS}.
\end{align*}
\end{definition}

The next lemma studies the stratification on $\Z$ created
by the notion of \emph{depth}.

\begin{lemma}\label{Dsets}
	The following statements hold:
	\begin{enumerate}[(a)]
	   \item $\DS[m+1] \subset \set{\ell \in \Z}{\exists\; i \in \DS
		        	\text{ such that }
					\wbasint{\ell} \varsubsetneq \wbasint{i}}.$
		\item For every $\ell, k \in \DS$ it follows that
		$\wbasint{\ell} \cap \wbasint{k} = \emptyset.$
	\end{enumerate}
\end{lemma}

\begin{proof}
Observe that if $\wbasint{\ell} \varsubsetneq \wbasint{i}$ then $\dep(\ell)\ge \dep(i)+1.$
Hence, (a) holds.

Statement (b) follows from Lemma~\ref{Propertiesvarphi}(g).
\end{proof}

In what follows, for every $m \in \Z^+$ we set
\[
	\wIBD := \bigcup_{i \in \DS} \wbasint{i}\supset \sstar{\DS}.
\]
Note that, by Lemma~\ref{Dsets}(b), $\wIBD$ is a disjoint union of closed arcs.
Therefore, for every $\theta \in \wIBD,$ there exists a unique
$i \in \DS$ such that $\theta \in \wbasint{i}.$
We will denote such integer $i$ by $\bt{\theta} \in \DS.$

The next two lemmas study the properties of the winged boxes
$\wbasint{i}$ and $\wbasicbox{i}$ according to the depth stratification.
Lemma~\ref{QuePassaALesAles} is the real motivation to introduce the winged boxes.

\begin{lemma}\label{denso}
The following statements hold:
\begin{enumerate}[(a)]
\item The sequence $\{\mu_m\}_{m=0}^{\infty}$ is strictly increasing.
In particular $\lim_{m\to\infty} \mu_m = \infty.$

\item For every $m \in \Z^+,$ $\wIBD$ is dense in $\SI,$
      $\wIBD[m+1]\subset \wIBD$
      and $\sstar{\DS} \cap \wIBD[m+1] = \emptyset.$

\item $\Orbom \subset \wIBD[0],$ and
      $\setfibth{\A} = \{(\theta,0)\}$
      for every $\theta \in \SI\setminus \wIBD[0].$

\item Let $i\in\Z$ and $\theta \in \wbasint{i} \setminus \wIBD[\dep(i)+1].$
Then, $\theta \notin \Orbom$ unless $\theta = \istar,$ and
$\setfibth{\A_n} = \setfibth{\A_{\ai}}$
for every $n \ge \ai.$
In particular $\setfibth{\A} = \setfibth{\A_{\ai}}.$
\end{enumerate}
\end{lemma}

\begin{proof}
By Lemmas~\ref{Dsets}(a) and \ref{Propertiesvarphi}(g) it follows that
for every $m \in \Z+$ and $\ell \in \DS[m+1]$ there exists $i\in \DS$
such that $\wbasint{\ell} \varsubsetneq \wbasint{i}$ and $\ai < \all.$
Thus, $\wIBD[m+1]\subset \wIBD$ and $\mu_m < \mu_{m+1}.$
This proves (a) and the second statement of (b).

Next we will show that $\istar \notin\wIBD[m+1]$
for every $i\in \DS.$
Assume by way of contradiction that there exists
$i\in \DS$ such that $\istar \in\wIBD[m+1].$
Let $k = \bt[m+1]{\istar} \in \DS[m+1].$
Clearly, $i \ne k$ and $\istar \in \wbasint{k}.$
Then, by Lemma~\ref{Propertiesvarphi}(g),
$\ak < \ai$  and $\wbasint{i} \varsubsetneq \wbasint{k}.$
Thus,
\[
  m = \dep(i) \ge \dep(k) + 1 = m+2;
\]
a contradiction.

Now we prove the first statement of (c).
From the definitions and the part of (b) already proven we have
\[
 \Orbom \subset \bigcup_{i \in \Z} \wbasint{i}
        \subset \bigcup_{m=0}^\infty \wIBD = \wIBD[0].
\]

To end the proof of (b) it remains to show the density of $\wIBD.$
We will do it by induction on $m.$
Clearly $\wIBD[0] \supset \Orbom$ is dense in $\SI$ because so is $\Orbom.$
Suppose that (b) holds for $\wIBD.$
We will show that (b) also holds for $\wIBD[m+1].$
Choose $\theta\in \wIBD$ and set $i = \bt{\theta}.$
Since $\Orbom$ is dense in $\SI,$ there exists a sequence
$\{s_{n}\}_{n=0}^{\infty} \subset \Z$ such that
$\sstar{s}_{n} \in \wobasint{i}$ and
$\lim_{n\to\infty} \sstar{s}_n = \theta.$
As above, we get that $\dep(s_{n}) \ge \dep(i) + 1 = m + 1.$
Moreover,
$\sstar{s}_n \in \wIBD[\dep(s_{n})] \subset \wIBD[m+1]$
for every $n.$
Consequently, $\wIBD \subset \overline{\wIBD[m+1]},$ and
the density of $\wIBD[m+1]$ follows from the density of $\wIBD.$

Next we prove the second statement of (c).
From above it follows that
\[
  \bigcup_{i \in \Z} \basintabs{i} \subset
  \bigcup_{i \in \Z} \wbasint{i} \subset \wIBD[0].
\]
Hence, by the definition of the maps $\gams{m}$ (Definition~\ref{PCgenerators})
it follows that $\gams{m}(\theta) = \gams{0}(\theta) = 0$
for every $\theta \notin \wIBD[0]$ and $m \in \Z^+.$
So, $\gamma(\theta) = \lim_{m\to\infty} \gams{m}(\theta) = 0,$
and
$\setfibth{\A} = \{(\theta,\gamma(\theta))\} = \{(\theta,0)\}$
by Lemma~\ref{propiedadesPCA}(c).
This ends the proof of (c).

\inidemopart{d}
If $\theta = \istar$ then the statement follows from
Lemmas~\ref{propiedadesA}(b) and \ref{propiedadesPCA}(b).
So, we assume that $\theta \ne \istar.$

By Definition~\ref{PCgenerators}(\tsfR.2) and
Remark~\ref{PCgeneratorsExplicitConsequences}(\tsfR.2)
we get that
$\theta \notin \Zstar_{\ai + 1}.$
Hence, if $\theta \in \Orbom,$ it follows that
$\theta = \kstar \in \wIBD[\dep(k)]$ with $\ak > \ai + 1$ and
$\wbasint{k} \cap \wbasint{i} \ne \emptyset.$
Thus, by Lemma~\ref{Propertiesvarphi}(g), $\dep(k) \ge \dep(i)+1$.
By (b), this implies that
$\theta = \kstar \in \wIBD[\dep(i)+1];$
a contradiction.
Therefore, $\theta \notin \Orbom.$
On the other hand,
$\theta \notin \wbasint{-i}$
by Definition~\ref{PCgenerators}(\tsfR.2).

If $\theta \notin \basintabs{k}$
for every $k \in \Z$ such that $\ak > \ai,$ then
$\gams{n}(\theta) = \gams{\ai}(\theta)$ and
$\setfibth{\A_n} = \setfibth{\A_{\ai}}$ for every $n \ge \ai,$
by Definition~\ref{PCgenerators} and Lemma~\ref{propiedadesA}(c).

Now assume that $\theta \in \basintabs{k}$ for some $k \in \Z$
such that $\ak > \ai$ and $\ak$ is minimal with these properties.
If $\theta \in \wobasint{k},$
as above we get that $\dep(k) \ge \dep(i)+1$ and
$\theta \in \wIBD[\dep(k)] \subset \wIBD[\dep(i)+1].$
Thus, $\theta \in \Bd(\wbasint{k}) = \Bd(\basintabs{k})$
and $k \ge 0.$
So, by Lemma~\ref{Propertiesvarphi}(c) and
the definition of the maps $\gams{j}$ (Definition~\ref{PCgenerators}),
$\gams{\ak}(\theta) = \gams{\ak - 1}(\theta).$
Moreover, by Lemma~\ref{Propertiesvarphi}(e),
$\gams{j}(\theta) = \gams{\ak}(\theta)$ for every $j > \ak.$
On the other hand, the minimality of $\ak$ implies that
$\theta \notin \basintabs{\ell}$ for every
$\ell \in \Z$ such that $\ak > \all > \ai.$
Hence, by the definition of the maps $\gams{j}$ (Definition~\ref{PCgenerators}),
$\gams{j}(\theta) = \gams{\ai}(\theta)$ for every $\ak > j > \ai.$
In short, we have proved that
$\gams{j}(\theta) = \gams{\ai}(\theta)$ for every $j \ge \ai.$
Thus, as above, $\setfibth{\A_n} = \setfibth{\A_{\ai}}$ for every $n \ge \ai.$
This ends the proof of the lemma.
\end{proof}

\begin{lemma}\label{QuePassaALesAles}
Assume that $\wbasint{i} \subset \wbasint{k}$
for some $i \in \DS,\ k \in \DS[m-1]$ and $m \in \N.$
Then, $\ak < \ai$ and $\akk < \aii$
unless $k \ge 0$ and $i = -(k+2)$ (whence $\akk = \aii$).
Moreover, the following statements hold:
\begin{enumerate}[(a)]
\item For every $\theta \in \wbasint{i},$
\[
 \gams{\ak}(\theta) = \gams{\ak+1}(\theta) = \dots = \gams{\ai-1}(\theta) \in \I_{i, \theta}
\]
and, when $\akk < \aii,$
\[
  \gams{\akk}\left(R_\omega(\theta)\right) =
  \gams{\akk+1}\left(R_\omega(\theta)\right) = \dots =
  \gams{\aii-1}\left(R_\omega(\theta)\right)
\]
\item For every $\theta \in \wbasint{i} \setminus \obasintabs{i},$
\[
  \gams{\ai}(\theta) = \gams{\ai-1}(\theta)
  \andq
  \I_{i, \theta} = \{\gams{\ai}(\theta)\} = \{\gams{\ak}(\theta)\} \subset \I_{k, \theta}.
\]
\end{enumerate}
\end{lemma}

\begin{proof}
The fact that $\ak < \ai$ follows from Lemma~\ref{Propertiesvarphi}(g).
Therefore, either
$\akk < \aii$ or
$k \ge 0,$ $i = -(k+2)$ and $\akk = \aii$ or
$k \ge 0,$ $i = -(k+1)$ and $\akk > \aii.$
In the last case, $\wbasint{i} = \wbasint{-(k+1)}$ and $\wbasint{k}$
must be disjoint by Definition~\ref{PCgenerators}(\tsfR.2) (with $j = k$);
which is a contradiction.
Thus $\akk < \aii$ unless $k \ge 0$ and $i = -(k+2)$ ($\akk = \aii$).

By Definition~\ref{PCgenerators}(\tsfR.2) and
Remark~\ref{PCgeneratorsExplicitConsequences}(\tsfR.2),
$\wbasint{i} \cap \Zstar_{\ai-1} = \emptyset.$
Hence, from the definition of the maps $\gams{j}$ (Definition~\ref{PCgenerators}),
to prove that
\[
  \gams{\ak}\evalat{\wbasint{i}} =
  \gams{\ak + 1}\evalat{\wbasint{i}} = \dots =
  \gams{\ai-2}\evalat{\wbasint{i}} =
  \gams{\ai-1}\evalat{\wbasint{i}},
\]
it is enough to show that
$\basintabs{\ell} \cap \wbasint{i} = \emptyset$
for every $\ell$ such that $\ak < \all < \ai.$
Assume that $\basintabs{\ell} \cap \wbasint{i} \ne \emptyset$
for some $\ell$ such that $\ak < \all < \ai.$
Then,
\[
 \emptyset \ne \basintabs{\ell} \cap \wbasint{i} \subset
 \wbasint{\ell} \cap \wbasint{i} \subset
 \wbasint{\ell} \cap \wbasint{k}
\]
and, by Lemma~\ref{Propertiesvarphi}(g),
\[
 \wbasint{i} \varsubsetneq \wbasint{\ell} \varsubsetneq \wbasint{k}.
\]
So, in a similar way as before,
\[
  m = \dep(i) \ge \dep(\ell) + 1 \ge \dep(k) + 2 = m+1;
\]
a contradiction.
This ends the proof of the first statement of (a).

Now we show that if $\akk < \aii-1,$ then
\[
\gams{\akk}\left(R_\omega(\theta)\right) =
  \gams{\akk+1}\left(R_\omega(\theta)\right) = \dots =
  \gams{\aii-1}\left(R_\omega(\theta)\right),
\]
and are well defined.

First we prove that
$\gams{\ell}\left(R_\omega(\theta)\right)$
is well defined for every $\ell = 0,1,\dots,\aii-1.$
For  every $\theta \in \wbasint{i}$ we have
\[
R_\omega(\theta) \in R_\omega\left(\wbasint{i}\right) =
\begin{cases}
    \BSG{i+1}{i}    &\text{when $i \ge 0$, and}\\
    \basintabs{i+1} \subset \wbasint{i+1} &\text{when $i < 0$}.
\end{cases}
\]
In any case, by Definition~\ref{PCgenerators}(\tsfR.2) and
Remark~\ref{PCgeneratorsExplicitConsequences}(\tsfR.2)
with $j=i$ when $i \ge 0$ and $\ell = -(j+1) = i+1$ when $i < 0,$
and Lemma~\ref{Propertiesvarphi}(a),
\[
R_\omega(\theta) \notin \begin{cases}
    \Zstar_{i}      &\text{when $i \ge 0$, and}\\
    \Zstar_{\aii-1} &\text{when $i < 0$,}
\end{cases}
\]
and $\gams{\ell}\left(R_\omega(\theta)\right)$ is well defined
for $\ell = 0,1,\dots,\aii-1$
(recall that $\Zstar_{m} \subset \Zstar_{m+1}$ for every $m \ge 0$).

Now, assume by way of contradiction that
\[
\gams{\ell}\left(R_\omega(\theta)\right) \ne \gams{\ell-1}\left(R_\omega(\theta)\right)
\text{ for some }
\ell \in \{\akk+1, \akk+2, \dots, \aii-1\},
\]
and $\ell$ is minimal with this property (observe that $\ell \ge 1$).
By the definition of the map $\gams{\ell}$ (Definition~\ref{PCgenerators}),
\[ R_\omega(\theta) \in \OBG{q+1}{\ell} \andq[with] q \in \{\ell-1, -(\ell+1)\} \]
and, hence, $\theta \in \OBG{q}{\ell}.$

Since $\akk+1 \le \ell < \aii$, when $q = -(\ell+1) \le -2,$
\[
  \akk+2 \le -q \le \aii
  \text{ and }
  \OBG{q}{\ell} = \wobasint{-(\ell+1)} = \wobasint{q}.
\]
Otherwise, when $q = \ell-1 \ge 0,$ $\akk \le q \le \aii-2$ and
\[
  \OBG{q}{\ell} \subset \obasint{\ell-1} = \wobasint{\ell-1} = \wobasint{q},
\]
by Definition~\ref{PCgenerators}(\tsfR.1).

Next we want to use Lemma~\ref{Propertiesvarphi}(g) to show that
$
\wbasint{i} \varsubsetneq \wbasint{q} \varsubsetneq \wbasint{k}.
$
To this end we have to compare $\aq$ with $\ai$ and $\ak.$

Notice $\wbasint{q} \cap \wbasint{k} \ne \emptyset$ because
\[
 \theta \in \wobasint{q} \cap \wbasint{i} \subset \wobasint{q} \cap \wbasint{k}.
\]
If $k \ge 0,\ \aq \ge \akk > \ak.$
When $k,q < 0,\ \aq \ge \akk + 2 = \ak + 1 > \ak.$
If $k < 0$ and $q \ge 0,\ \aq = q \ge \akk = \ak - 1.$
If $q = \ak - 1$ (that is, $k = -(q+1)$), as above,
by Definition~\ref{PCgenerators}(\tsfR.2) with $j = q$ we get
$\wbasint{k} \cap \wbasint{q} = \emptyset;$ a contradiction.
So, $\aq > \ak$ unless $\aq = \ak$ and $k < 0 \le q.$
Summarizing, we have shown that $\aq \ge \ak$ and $q \ne k.$
Then, from Lemma~\ref{Propertiesvarphi}(g) we get that $\aq > \ak$
and $\wbasint{q} \varsubsetneq \wbasint{k}.$

Now we will study the relation of $\wbasint{q}$ with the box $\wbasint{i}.$
From above we get that $\wbasint{q} \cap \wbasint{i} \ne \emptyset.$
If $i < 0,\ \aq \le \aii = \ai-1.$
When $q,i \ge 0,$ we have $\aq = q \le \aii-2 = \ai - 1.$
If $i \ge 0$ and $q < 0,\ \aq \le \aii = \ai + 1.$

Assume that $i \ge 0$ and $q = -(i+1) < 0.$
In this case, additionally, $q = -(\ell+1)$ and, thus, $i = \ell \ge 1.$
Then,
\begin{align*}
R_\omega(\theta) & \in R_\omega\left(\wbasint{i}\right)
                   = R_\omega\left(\basint{i}\right)
                   = \BSG{i+1}{i}, \text{ and}\\
R_\omega(\theta) & \in \OBG{q+1}{\ell} = \OBG{-i}{i} \subset \wobasint{-i},
\end{align*}
which is a contradiction by Definition~\ref{PCgenerators}(\tsfR.2).
Summarizing, $\aq < \ai$ unless $\aq = \ai$ and $q < 0 \le i$
(that is, $\aq \le \ai$ and $q \ne i$).
Then, again by Lemma~\ref{Propertiesvarphi}(g), $\aq < \ai$
and $\wbasint{i} \varsubsetneq \wbasint{q} \varsubsetneq \wbasint{k}.$
So, as before,
\[
 m = \dep(i) \ge \dep(q) + 1 \ge \dep(k) + 2 = m+1;
\]
a contradiction.
This ends the proof of (a).

Now we assume that $\theta \in \wbasint{i} \setminus \obasintabs{i}.$
By Lemmas~\ref{Propertiesvarphi}(e) and \ref{voresdelescaixesalesvores}(d),
\[
  \gams{\ai}(\theta) = \gams{\ai-1}(\theta)
  \andq
  \I_{i, \theta} = \{\gams{\ai}(\theta)\} =
      \{\gams{\ai-1}(\theta)\} = \{\gams{\ak}(\theta)\}.
\]

On the other hand, by Lemma~\ref{denso}(b),
$\sstarplain{\DS[m-1]} \cap \wIBD = \emptyset$
which implies that $\theta \ne \kstar$
because $\kstar \in \sstarplain{\DS[m-1]}$ and
$\theta \in \wbasint{i} \setminus \obasintabs{i} \subset \wIBD.$
So, by \eqref{gammathetaproperties},
\[
 \I_{i, \theta} = \{\gams{\ai}(\theta)\} = \{\gams{\ak}(\theta)\} \subset \I_{k, \theta}.
\]

Now we prove that
$\gams{\ai-1}(\theta) \in \I_{i, \theta}$
for every $\theta \in \wbasint{i}.$
From above, we have $\I_{i, \theta} = \{\gams{\ai-1}(\theta)\}$
for every $\theta \in \wbasint{i} \setminus \obasintabs{i}.$
Moreover, when $\theta \in \obasintabs{i}$
the statement follows directly from Lemma~\ref{Propertiesvarphi}(c).
Thus, (b) is proved.
\end{proof}

\section{Boxes in the wings}\label{BoxesintheWings}

To prove Theorem~\ref{MainTh} we will inductively construct a Cauchy sequence
$\{T_m\}_{m=0}^{\infty} \subset \cSO$
that gives the function $T$ from Theorem~\ref{MainTh} as a limit.

This section is devoted to study the points in the wings of
boxes in the circle and its interaction with boxes of higher depth.
The resulting technology is necessary to be able to construct the
sequence $\{T_m\}_{m=0}^{\infty}$ so that it is Cauchy sequence.
Unfortunately this will complicate even more the definition of the
functions $T_m$ and the proof of its continuity.

We start by introducing some more notation.
For every $m \in \Z^+$ we set
\begin{align*}
 \IBD   &:= \bigcup_{i \in \DS} \basintabs{i} \subset \wIBD,\text{ and}\\
 \WDB  &:= \Bigl\{\theta \in \wIBD\setminus\IBD \,\colon
                  \theta \in \IBD[j] \text{ for some $j > m$}
 \Bigr\}.
\end{align*}
On the other hand, the smallest number $j$ from the above definition will be called the
\emph{least essential depth of $\theta$ below $m$,}
and will be denoted by $\led{\theta}.$
That is, $\led{\theta}$ denotes the positive integer larger than $m$
such that
\[
\theta \in \wIBD[j] \setminus \IBD[j]
\text{ for } j= m, m+1,\dots, \led{\theta}-1
\andq
\theta \in \IBD[\led{\theta}].
\]

The following simple lemmas are useful
to better understand and use the above definitions.
The next lemma establishes the relation between boxes in the wings of increasing depth.

\begin{lemma}\label{DepthintheWings}
Assume that $\theta \in \WDB$ for some $m \in \Z^+$
and set $\ell = \led{\theta}.$
Then, the following statements hold.
\begin{enumerate}[(a)]
\item For every $j = m, m+1,\dots, \ell$
the numbers $\is_j = \bt[j]{\theta} \in \DS[j]$ are well
defined and are all of them negative except, perhaps,
$\is_{\ell} = \bt[\led{\theta}]{\theta}.$

\item
\begin{align*}
& \hspace*{2em}\abs{\is_{m}} < \abs{\is_{m+1}} < \dots <
    \abs{\is_{\ell-1}} < \abs{\is_{\ell}},\text {and}\\
\theta &\in \basintabs{\is_{\ell}}
    \subset \wobasint{\is_{\ell-1}} \setminus \basintabs{\is_{\ell-1}}\\
   &\subset \wobasint{\is_{\ell-2}} \setminus \basintabs{\is_{\ell-2}}
    \subset \cdots
    \subset \wobasint{\is_{m}} \setminus \basintabs{\is_{m}}.
\end{align*}

\item For every $j = m, m+1,\dots,\ell-1,$
      $\basintabs{\is_{\ell}} \subset \WDB[j],$
      $\led[j]{\nu} = \led{\theta}$ and
      $\bt[\led[j]{\nu}]{\nu} = \bt[\led{\theta}]{\theta} = \is_{\ell}$
      for every $\nu \in \basintabs{\is_{\ell}}.$

\item $\I_{\is_{m}, \nu} = \{\gams{\abs{\is_{m}}}(\nu)\} \subset \I_{\is_{\ell}, \nu}$
      for every $\nu \in \obasintabs{\is_{\ell}}$ and
      \[
        \I_{\is_{m}, \nu} = \{\gams{\abs{\is_{m}}}(\nu)\} =
        \{m_{\is_{\ell}}(\nu)\} = \{M_{\is_{\ell}}(\nu)\} =
        \{\gams{\abs{\is_{\ell}}}(\nu)\} = \I_{\is_{\ell}, \nu}
      \]
      for every $\nu \in \Bd\left(\basintabs{\is_{\ell}}\right).$
\end{enumerate}
\end{lemma}

\begin{proof}
Since $\wbasint{i} = \basint{i}$ for every $i \ge 0,$
\begin{equation}\label{thewings}
\wIBD\setminus\IBD = \bigcup_{\substack{i \in \DS\\ i < 0}}
              \left(\wbasint{i} \setminus \basintabs{i}\right)
\end{equation}
for every $m \in \Z^+.$

Statement (a) follows from Lemma~\ref{denso}(b) and \eqref{thewings}.
Then, (b) follows from Lemma~\ref{Propertiesvarphi}(g).
Statement (c) is an easy consequence of (b) and the definitions.

Now we prove (d) iteratively.
Fix $\nu \in \obasintabs{\is_{\ell}}.$
By (b)
\[
\nu \in \wobasint{\is_{m+1}} \setminus \basintabs{\is_{m+1}}
       \subset \wobasint{\is_{m}} \setminus \basintabs{\is_{m}}
\]
provided that $\ell = \led{\theta} > m+1.$
Hence, by Lemmas~\ref{voresdelescaixesalesvores}(d) and \ref{QuePassaALesAles},
\begin{align*}
\gams{\abs{\is_{m}}}(\nu) &= \gams{\abs{\is_{m}}+1}(\nu) = \dots = \gams{\abs{\is_{m+1}}}(\nu),\text{ and}\\
\I_{\is_{m}, \nu} &= \{\gams{\abs{\is_{m}}}(\nu)\} = \{\gams{\abs{\is_{m+1}}}(\nu)\} = \I_{\is_{m+1}, \nu}.
\end{align*}
By iterating this argument we get,
\[
\gams{\abs{\is_{m}}}(\nu) = \gams{\abs{\is_{m}}+1}(\nu) = \dots = \gams{\abs{\is_{\ell-1}}}(\nu)
\andq
\I_{\is_{m}, \nu} = \I_{\is_{\ell-1}, \nu}.
\]
Again by (b) and Lemmas~\ref{voresdelescaixesalesvores}(d) and \ref{QuePassaALesAles},
\[
\gams{\abs{\is_{m}}}(\nu) = \gams{\abs{\is_{m}}+1}(\nu) = \dots = \gams{\abs{\is_{\ell}}}(\nu)
\andq
\I_{\is_{m}, \nu} = \I_{\is_{\ell}, \nu}
\]
when $\nu \in \Bd\left(\basintabs{\is_{\ell}}\right)$
and, otherwise,
\[
\gams{\abs{\is_{m}}}(\nu) = \gams{\abs{\is_{m}}+1}(\nu) = \dots = \gams{\abs{\is_{\ell}-1}}(\nu)
\andq
\I_{\is_{m}, \nu} \subset \I_{\is_{\ell}, \nu}.
\]
\end{proof}

Equipped with above results and definition we are going to define two maps,
analogous to the maps $m_i$ and $M_i,$ on the wings of the negative boxes.

\begin{definition}\label{curvesinthewings}
For every $m \in \Z^+$ we define
\begin{align*}
 \WDS   &:= \set{\bt[\led{\theta}]{\theta}}{\theta \in \WDB} \subset \Z,\\
 \WIB   &:= \Int(\WDB) = \LSleftlimits{\bigcup}{i \in \WDS} \obasintabs{i},\\
 \WB    &:= \bigcup_{\substack{i \in \DS\\ i < 0}} \left(\wbasint{i} \setminus \obasintabs{i}\right),\text{ and}\\
 \wEIBD & := \bigcup_{i \in \DS} \Bd\left(\wbasint{i}\right) \subset \wIBD.
\end{align*}

By Lemma~\ref{DepthintheWings}(a,c), $\WDS$ is well defined and
\[
\WIB \subset \WDB \subset \wIBD\setminus\IBD \subset \WB.
\]
Consequently,
\[ \wIBD = \IBD \cup \WB. \]

Then, we can define functions
{\map{\tau_m}{\WB}[\I]}
and
{\map{\lambda_m}{\WB}[\I]}
as follows:
\begin{align*}
\tau_m(\theta) &:= \begin{cases}
   M_{\bt[\led{\theta}]{\theta}}(\theta) & \text{if $\theta \in \WIB$,}\\
   \gams{\abs{\bt{\theta}}}(\theta) & \text{otherwise,}
\end{cases}\\
\lambda_m(\theta) &:= \begin{cases}
   m_{\bt[\led{\theta}]{\theta}}(\theta) & \text{if $\theta \in \WIB$,}\\
   \gams{\abs{\bt{\theta}}}(\theta) & \text{otherwise.}
\end{cases}
\end{align*}
Clearly, by Lemmas~\ref{voresdelescaixesalesvores}(a) and \ref{Propertiesvarphi}(b),
\[
 -1 \le \lambda_m(\theta) \le \tau_m(\theta) \le 1
\]
for every $\theta \in \WB.$
So, we can define
\[
\IW{\theta} := [\lambda_m(\theta), \tau_m(\theta)] \subset [0,1].
\]
\end{definition}

The next lemmas will help us in the definition and study of the maps $T_m$.

\begin{lemma}\label{VerticalIntervalsIntheWings}
The following statements hold for every $m \in \Z^+.$
\begin{enumerate}[(a)]
\item $\WIB \cap \IBD = \WIB \cap \wEIBD = \emptyset.$
\item Let $\theta \in \WB.$ Then,
$
\I_{\bt{\theta}, \theta} = \left\{\gams{\abs{\bt{\theta}}}(\theta)\right\},
$
\begin{align*}
\I_{\bt{\theta}, \theta} &= \IW{\theta} && \text{when $\theta \notin \WIB$, and}\\
\I_{\bt{\theta}, \theta} &\subset \IW{\theta} && \text{when $\theta \in \WIB$.}
\end{align*}
\item Assume that $m \in \N$ and let $U$
be a connected component of $\WB$ such that $U \subset \WB[m-1].$
Then,
$\WDB \cap U \subset \WDB[m-1],$
$\WIB \cap U = \WIB[m-1] \cap U$ and
$\IW{\theta} = \IW[m-1]{\theta}$ for every $\theta \in U.$
\end{enumerate}
\end{lemma}

\begin{proof}
\inidemopart{a}
By Lemma~\ref{DepthintheWings}(b),
\[
 \theta \in  \wobasint{\bt{\theta}} \setminus \basintabs{\bt{\theta}}
\]
and $\bt{\theta} < 0$ for every $\theta \in \WIB \subset \WDB.$
So, by Lemma~\ref{Dsets}(b), we get
$\theta \notin \IBD \cup \wEIBD.$

\inidemopart{b}
The fact that
$
\I_{\bt{\theta}, \theta} = \left\{\gams{\abs{\bt{\theta}}}(\theta)\right\}
$
follows from Lemma~\ref{voresdelescaixesalesvores}(d).
The other two statements follow from Definition~\ref{curvesinthewings}
and Lemma~\ref{DepthintheWings}(d).

\inidemopart{c}
The assumption that $U$ is a connected component of $\WB$
and $U \subset \WB[m-1]$ implies
by Lemmas~\ref{Dsets}(b) and \ref{Propertiesvarphi}(g)
that there exist $i\in \DS$ and $k \in \DS[m-1],$ $i,k < 0,$
such that $U$ is a connected component of
\[
   \wbasint{i}\setminus\obasintabs{i} \subset
      \wobasint{k}\setminus\basintabs{k} \subset \WB[m-1].
\]
Again by Lemma~\ref{Dsets}(b)
this implies that $U \subset \wIBD[m-1]\setminus\IBD[m-1].$
Moreover, by definition,
$\WDB \subset \wIBD\setminus\IBD.$
Consequently, $\WDB \cap U \subset \WDB[m-1].$

Let $\theta \in \WIB \cap U \subset \WDB \cap U \subset \WDB[m-1] \cap U.$
By  Definition~\ref{curvesinthewings} and Lemma~\ref{DepthintheWings}(a,b),
$i = \bt{\theta}$ and
there exists $\ell = \bt[\led{\theta}]{\theta} \in \WDS$ such that
\[
 \theta \in \obasintabs{\ell} \subset
      \wbasint{i}\setminus\obasintabs{i} \subset
      \wobasint{k}\setminus\basintabs{k}.
\]
Therefore, again by
Lemma~\ref{DepthintheWings}(a--c) and Definition~\ref{curvesinthewings},
$\led[m-1]{\theta} = \led{\theta},$
\[ \ell = \bt[\led{\theta}]{\theta} = \bt[\led[m-1]{\theta}]{\theta} \in \WDS[m-1] \]
and
$\theta \in \obasintabs{\ell} \subset \WIB[m-1].$
Hence, $\WIB \cap U \subset \WIB[m-1].$

Now assume that $\theta \in \WIB[m-1] \cap U.$
As above, there exist $r = \bt{\theta} \in \DS$
and $\ell = \bt[\led[m-1]{\theta}]{\theta} \in \WDS[m-1]$ such that
\[
 \theta \in \obasintabs{\ell} \subset
      \wobasint{r}\setminus\basintabs{r} \subset
      \wobasint{k}\setminus\basintabs{k}.
\]
Since $\theta \in U \subset \wbasint{i}$,
Lemma~\ref{Dsets}(b) gives
$i = r$ and $\theta \in \obasintabs{\ell} \subset U.$
Moreover, by Lemma~\ref{DepthintheWings}(c),
$\ell = \bt[\led[m-1]{\theta}]{\theta} = \bt[\led{\theta}]{\theta}\in \WDS$
and, so, $\theta \in \obasintabs{\ell} \subset \WIB.$
Thus, $\WIB \cap U = \WIB[m-1] \cap U.$

To end the proof of the lemma we have to show that
$\IW{\theta} = \IW[m-1]{\theta}$ for every $\theta \in U.$
Assume first that
$\theta \in U \setminus \WIB \subset \WB \setminus \WIB.$
Then,
\[
\theta \in U \setminus \WIB =
           U \setminus \WIB[m-1] \subset
           \WB[m-1] \setminus \WIB[m-1]
\]
and, by (b) and Lemmas~\ref{voresdelescaixesalesvores}(d) and \ref{QuePassaALesAles},
\[
\IW{\theta} = \I_{i, \theta} = \{\gams{\ai}(\theta)\} =
              \{\gams{\ak}(\theta)\} = \I_{k, \theta} = \IW[m-1]{\theta}.
\]
If $\theta \in U \cap \WIB = U \cap \WIB[m-1]$ then we get
\begin{align*}
 \IW{\theta}
   &= \left[m_{\bt[\led{\theta}]{\theta}}(\theta),
            M_{\bt[\led{\theta}]{\theta}}(\theta)\right] \\
   &= \left[m_{\bt[\led[m-1]{\theta}]{\theta}}(\theta),
            M_{\bt[\led[m-1]{\theta}]{\theta}}(\theta)\right]
    = \IW[m-1]{\theta}
\end{align*}
from Definition~\ref{curvesinthewings} and
Lemma~\ref{DepthintheWings}(c).
\end{proof}

\begin{lemma}\label{continuouscurvesinthewings}
Let $m \in \Z^+$ and let $U$ be a connected component of $\WB.$
Then, the functions $\lambda_m\evalat{U}$ and $\tau_m\evalat{U}$
are continuous.
\end{lemma}

\begin{proof}
We will prove only the continuity of $\lambda_m\evalat{U}.$
The proof of the continuity of $\tau_m\evalat{U}$ is analogous.

By Lemmas~\ref{DepthintheWings}(c) and \ref{voresdelescaixesalesvores}(b)
we get
\begin{equation}\label{contindepboxes}
\text{\parbox{0.9\textwidth}{
for every $\ell \in \WDS,$
$\ell = \bt[\led{\nu}]{\nu}$ for every $\nu \in \basintabs{\ell},$
and the function $m_{\ell}$ is continuous on $\basintabs{\ell}.$
}}\end{equation}

Let $\ell \in \WDS$ be such that $\obasintabs{\ell} \subset \WIB \cap U.$
Thus, by \eqref{contindepboxes}, the function
$\lambda_m = m_{\ell}$ is continuous on $\obasintabs{\ell}.$

So, we have to show that $\lambda_m$ is continuous at every
$\theta \in U \setminus \WIB.$ To show this we will use a simple usual
$\varepsilon$--$\delta$ game. Fix $\varepsilon > 0.$

By Lemma~\ref{Dsets}(b) it follows that $U$ is a connected component
of $\wbasint{i} \setminus \obasintabs{i}$ for some $i \in \DS,$ $i < 0,$
and
\begin{equation}\label{whichbox}
\bt{\nu} = i \andq[for every] \nu \in U.
\end{equation}
By Lemma~\ref{Propertiesvarphi}(a) and
Definition~\ref{PCgenerators}(\tsfR.2) and
Remark~\ref{PCgeneratorsExplicitConsequences}(\tsfR.2),
the function $\gams{\ai}\evalat{U}$ is continuous.
So,
\begin{equation}\label{gamscont}
\text{\parbox{0.9\textwidth}{
there exists $\overline{\delta}_{\ai} = \overline{\delta}_{\ai}(\theta) > 0$ such that
$\abs{\gams{\ai}(\theta), \gams{\ai}(\nu)} < \varepsilon/2$
provided that $\dSI(\theta, \nu) < \overline{\delta}_{\ai}.$
}}\end{equation}

On the other hand, by \eqref{contindepboxes},
\begin{equation}\label{contindepboxesepsdel}
\text{\parbox{0.9\textwidth}{
for every $\ell \in \WDS,$ there exists $\delta_{\ell} > 0$ such that
$\abs{m_{\ell}(\widetilde{\theta}), m_{\ell}(\nu)} < \varepsilon/2$
for every
$\widetilde{\theta} \in \Bd\left(\basintabs{\ell}\right)$ and
$\nu \in \Bd\basintabs{\ell}$
such that $\dSI(\theta, \nu) < \delta_{\ell}.$
}}\end{equation}

Now we will define $\delta.$
Note that there exists $N \in \N$ such that $2^{-N} < \varepsilon/2.$
Then we set:
\[
 \delta = \delta(\theta) := \min\left\{
    \overline{\delta}_{\ai}(\theta),
    \min \set{\delta_{\ell}}{\ell \in \WDS \text{ and } \all < N}
\right\}.
\]
Clearly, $\delta > 0$ because the set $\set{\ell \in \WDS}{\all < N}$ is finite.

To end the proof of the lemma we have to show that
\[
 \abs{\lambda_m(\theta) - \lambda_m(\nu)} < \varepsilon
\]
whenever $\nu \in U$ and $\dSI(\theta, \nu) < \delta.$

Assume that $\nu \in U$ and $\dSI(\theta, \nu) < \delta$
(recall that we have the assumption that $\theta \notin \WIB$).
If $\nu \notin \WIB,$ then
$\dSI(\theta, \nu) < \delta \le \overline{\delta}_{\ai}(\theta)$
and, by \eqref{whichbox} and \eqref{gamscont},
\[
 \abs{\lambda_m(\theta) - \lambda_m(\nu)} =
 \abs{\gams{\ai}(\theta) - \gams{\ai}(\nu)} < \varepsilon/2 < \varepsilon.
\]

Now assume that there exists $\ell \in \WDS$
such that $\nu \in \obasintabs{\ell} \subset \WIB.$
Clearly, there exists
$\widetilde{\theta} \in \Bd\left(\basintabs{\ell}\right)$
such that
\begin{align*}
& \dSI(\theta, \widetilde{\theta}) < \dSI(\theta, \nu) < \delta \le \overline{\delta}_{\ai}(\theta)\text{ and}\\
& \dSI(\widetilde{\theta}, \nu) < \dSI(\theta, \nu) < \delta.
\end{align*}
Observe that, by Lemma~\ref{Dsets}(b), $\widetilde{\theta} \notin \WIB.$
Hence, by \eqref{whichbox} and Lemma~\ref{DepthintheWings}(c,d),
\[
 \lambda_m(\widetilde{\theta}) = \gams{\ai}(\widetilde{\theta}) = m_{\ell}(\widetilde{\theta}).
\]

If $\all < N,$ then $\dSI(\widetilde{\theta}, \nu) < \delta \le \delta_{\ell}$ and,
by \eqref{contindepboxesepsdel},
$
\abs{m_{\ell}(\widetilde{\theta}) - m_{\ell}(\nu)} < \varepsilon/2.
$
Otherwise, by Lemma~\ref{Propertiesvarphi}(f),
\[
\abs{m_{\ell}(\widetilde{\theta}) - m_{\ell}(\nu)} <
  \diam\left(\basicbox{\ell}\right)
      \le  2^{-\all} \le 2^{-N} < \varepsilon/2.
\]
In any case,
$
\abs{m_{\ell}(\widetilde{\theta}) - m_{\ell}(\nu)} < \varepsilon/2.
$
Thus, again by \eqref{whichbox} and \eqref{gamscont},
\begin{align*}
\abs{\lambda_m(\theta) - \lambda_m(\nu)}
  & \le \abs{\lambda_m(\theta) - \lambda_m(\widetilde{\theta})} +
        \abs{\lambda_m(\widetilde{\theta}) - \lambda_m(\nu)}\\
  & = \abs{\gams{\ai}(\theta) - \gams{\ai}(\widetilde{\theta})} +
      \abs{m_{\ell}(\widetilde{\theta}) - m_{\ell}(\nu)} < \varepsilon.
\end{align*}
\end{proof}

\section{A Cauchy sequence of skew products. Proof of Theorem~\ref{MainTh}}\label{skew-product}

In this section prove Theorem~\ref{MainTh}.
To do this we inductively construct a Cauchy sequence
$\{T_m\}_{m=0}^{\infty} \subset \cSO$
that gives the function $T$ from Theorem~\ref{MainTh} as a limit.

The sequence $\{T_m\}_{m=0}^{\infty} \subset \cSO$ is defined so that
\[
  T_m(\theta,x) = (R_\omega(\theta), f_m(\theta, x))
\]
and $\map{f_m}{\Omega}[\I]$ is continuous in both variables.
To build these functions we will use the auxiliary functions
{\map{G_i}{\basicbox{i}}[\Omega]} with $i \in \Z$
from Section~\ref{FunctionsGi}.
The maps $f_m(\theta, \cdot)$ will also be denoted as $f_{m,\theta},$
and will be defined non-increasing, and such that
$f_{m,\theta}(2) = -2$ and $f_{m,\theta}(-2) = 2$
for every $\theta \in \SI.$

To make more evident the strategy of the construction of this sequence of maps
we will separate several cases, and we will state without proofs
the results that study these maps.
After establishing all the definitions and results related to the
construction of the sequence $\{T_m\}_{m=0}^{\infty}$
without having been distracted by the technicalities involving
the proofs, we will proceed to provide the missing proofs.
More precisely, we will start by defining the map
$T_0$ and stating without proof the  proposition that summarizes
the necessary properties of this map.
Next we will inductively define the maps
$\{T_m\}_{m=1}^{\infty} \subset \cSO$
and state without proof the proposition that establishes the
properties of the whole sequence
$\{T_m\}_{m=0}^{\infty}$.

Then, as we have said, we prove Theorem~\ref{MainTh}
and in the next three sections we will provide all pending proofs.

In what follows $\mathcal{C}(\I,\I)$
will denote the class of all continuous maps from $\I$ to itself.
We endow $\mathcal{C}(\I,\I)$ with the supremum metric denoted by
$\norm{\cdot}$ so that $(\mathcal{C}(\I,\I), \norm{\cdot})$
is a complete metric space.

Next we define the map $T_0$.

\begin{definition}[The map $T_0$]\label{T0mapDefi}
Assume first that $\theta \in \wIBD[0]$ and let $\is = \bt[0]{\theta}$
(that is $\theta \in \wbasint{\is}$).
In this case we set:
\[
f_{0,\theta}(x) = \begin{cases}
   g_{_{\is,\theta}}(x)
        & \text{if $x\in \I_{\is,\theta}$}, \\[0.75ex]
   \frac{g_{_{\is,\theta}}\left(m_{\is}(\theta)\right) - 2}{m_{\is}(\theta) + 2} (x + 2) + 2
        & \text{if $x \in [-2,m_{\is}(\theta)]$},\\[1ex]
      \frac{g_{_{\is,\theta}}\left(M_{\is}(\theta)\right) + 2}{M_{\is}(\theta) - 2} (x - 2) - 2
        & \text{if $x \in [M_{\is}(\theta),2]$}.
\end{cases}
\]
If $\theta \in \SI \setminus \wIBD[0]$ then we define
$f_{0,\theta}$ to be the unique piecewise affine map
with two affine pieces whose graph joins the point
$(-2,2)$ with $(0, \gamma(R_\omega(\theta))),$
and this with the point $(2, -2)$.
\end{definition}


Next we introduce some more notation to be able to define the maps
$\{T_m\}_{m=1}^{\infty}.$
For every $k \in \Z$ we set
\[
  \wbasband{k} := \setsilift{\wbasint{k}} = \wbasint{k} \times \I
\]
and, for every $m \in \Z^+,$
\[
 \wIVD := \setsilift{\wIBD} = \wIBD \times \I = \bigcup_{i \in \DS} \wbasband{i}.
\]

\begin{definition}[The maps $T_m$ with $m > 0$]\label{seqTmDefi}
Now we assume that we have defined the function
$T_{m-1}$ for some $m \ge 1$ and we define
\[
  T_m(\theta,x) = (R_\omega(\theta), f_m(\theta, x))
\]
as follows.
By Lemma~\ref{Dsets}(b),
for every $(\theta, x) \in \wIVD,$ we have
\[
 \theta \in \wbasint{\is} \subset \wIBD
 \andq[with]
 \is = \bt{\theta} \in \DS
\]
(and, of course, $x \in \I$).
Then we define:
\[
f_{m,\theta}(x) = \begin{cases}
   f_{m-1,\theta}(x)
        & \text{if $\theta \in \SI \setminus \wIBD;\ x\in \I$},\\
   g_{_{\is,\theta}}(x)
        & \text{if $\theta \in \IBD;\ x \in \I_{\is,\theta}$}, \\[0.75ex]
   \frac{2 - g_{_{\is,\theta}}\left(m_{\is}(\theta)\right)\hfill}{
         2 - f_{m-1,\theta}\left(m_{\is}(\theta)\right)
        } (f_{m-1,\theta}(x) - 2) + 2
        & \text{if $\theta \in \IBD;\ x \in [-2,m_{\is}(\theta)]$},\\[1ex]
   \frac{2 + g_{_{\is,\theta}}\left(M_{\is}(\theta)\right)\hfill}{
         2 + f_{m-1,\theta}\left(M_{\is}(\theta)\right)
        } (f_{m-1,\theta}(x) + 2) - 2
        & \text{if $\theta \in \IBD;\ x \in [M_{\is}(\theta),2]$},\\
   \gams{\abs{\is+1}}\left(R_\omega(\theta)\right)
        & \text{if $\theta \in \WB;\ x \in \IW{\theta}$}, \\[0.75ex]
   \frac{2 - \gams{\abs{\is+1}}\left(R_\omega(\theta)\right)\hfill}{
         2 - f_{m-1,\theta}\left(\lambda_m(\theta)\right)
        } (f_{m-1,\theta}(x) - 2) + 2
        & \text{if $\theta \in \WB;\ x \in [-2,\lambda_m(\theta)]$},\\[1ex]
   \frac{2 + \gams{\abs{\is+1}}\left(R_\omega(\theta)\right)\hfill}{
         2 + f_{m-1,\theta}\left(\tau_m(\theta)\right)
        } (f_{m-1,\theta}(x) + 2) - 2
        & \text{if $\theta \in \WB;\ x \in [\tau_m(\theta),2]$}.
\end{cases}
\]
Since $\wIVD \subset \wIVD[m-1],$ $f_{m-1,\theta}$ is defined on $\wIVD$.
Moreover, the above formula defines $f_{m,\theta}$ for every
$\theta \in \wIBD$ since, by Definition~\ref{curvesinthewings},
$\wIBD = \IBD \cup \WB.$
We also remark that $f_{m,\theta}$ formally is defined in two different
ways when $\theta \in \WB \cap \IBD.$ Later on we will show that
$f_{m,\theta}$ is well defined.
\end{definition}

The next proposition studies the maps $\{T_m\}_{m=0}^{\infty}$
and describes their properties.

\begin{proposition}\label{seqTmProperties}\label{T0mapProperties}
The following statements hold for every $m \in \Z^+.$
\begin{enumerate}[(a)]
\item The map $T_m$ is well defined, continuous
      and belongs to $\cSO$.
\item For every $\theta\in\SI,$
      $f_{m,\theta}$ is non-increasing, and
      $f_{m,\theta}(2) = -2,$ $f_{m,\theta}(-2) = 2.$
      Moreover,
      $-1 \le f_{0,\theta}\left(M_{\bt{\theta}}(\theta)\right) \le
              f_{0,\theta}\left(m_{\bt{\theta}}(\theta)\right) \le 1$
      for every $\theta \in \wIBD.$
\item For every $i \in \DS,$
      $T_m\evalat{\wbasicbox{i}} = G_i,$
      $T_m\left(\setfibpt{\A_{\ai}}{\istar}\right) = \setfibpt{\A_{\aii}}{\sstar{i+1}},$
      and\newline
      $T_k\evalat{\istarset \times \I} = T_m\evalat{\istarset \times \I}$
      (that is, $f_{k,\istar} = f_{m,\istar}$) for every $k > m.$
\end{enumerate}
\end{proposition}

The next result shows that the sequence $\{T_m\}_{m=0}^{\infty}$
has a limit in $\cSO$.

\begin{proposition}\label{distTmTm-1}
For every $m \ge 2$ and $\theta \in \SI,$
\begin{equation}\label{fitanorma}
 \norm{f_{m,\theta} - f_{m-1,\theta}} \le 2 \cdot 2^{-\abs{\bt[m-1]{\theta}}}.
\end{equation}
Moreover, the sequence $\{T_m\}_{k = 0}^{\infty}$ is a Cauchy sequence.
\end{proposition}

Finally we are ready to prove the main result of the paper.
It follows from the next result which gives a more concrete version
of Theorem~\ref{MainTh}.

\begin{theorem}
There exists a map $T \in \cSO$ with
$f(\theta,\cdot)$ non-increasing for every $\theta \in \SI,$
such that $T$ permutes the upper and lower circles of $\Omega$
(thus having a periodic orbit of period two of curves),
and there exists a connected pseudo-curve $\A \subset \Omega$
which does not contain any arc of a curve
such that $T(\A) = \A$ and there does not exist
any $T$-invariant curve.
\end{theorem}

\begin{proof}
By Propositions~\ref{seqTmProperties} and \ref{distTmTm-1},
there exists a map
\[
  T(\theta, x) = (R_\omega(\theta), f(\theta, x))
       = (R_\omega(\theta), \lim_{m\to\infty} f_m(\theta, x)) \in \cSO
\]
with $f(\theta,\cdot)$ non-increasing for every $\theta \in \SI$
such that $T$ permutes the upper and lower circles of $\Omega$
(that is, $f(\theta,2) = -2$ and $f(\theta,-2) = 2$).
As the connected set $\A$ we take the one given by
Proposition~\ref{teoremacentral} (and Definition~\ref{PCAtLast}).

To end the proof of the theorem we need to show that $T(\A) = \A,$
since this already implies that there does not exist any $T$-invariant curve.
To see it, assume by way of contradiction that there exists an invariant curve
and denote its graph by $B.$
Since $B$ is the graph of a (continuous) curve, it is
compact and connected.
On the other hand, let $\Omega_{+}$ and $\Omega_{-}$ be the two
connected components of $\Omega\setminus \A$
from the proof of Proposition~\ref{teoremacentral}.
The facts that
$T(\A) = \A,$
$f(\theta,\cdot)$ is decreasing for every $\theta \in \SI,$ and
$T$ permutes the upper and lower circles of $\Omega$
imply that
$T(\Omega_{+}) = \Omega_{-}$
and
$T(\Omega_{-}) = \Omega_{+}.$
Hence, by the invariance of $B,$
$B \nsubseteq \Omega_{+}$ and $B \nsubseteq \Omega_{-}.$
The connectivity of $\A$ and $B$ imply that there exists
$(\theta,x) \in \A \cap B.$
Consequently,
\[
 B = \overline{\set{T^n(\theta,x)}{n\in\Z^+}} \subset \A;
\]
a contradiction because $\A$ does not contain any arc of a curve.

So, only it remains to prove that $T(\A) = \A.$
By using Proposition~\ref{seqTmProperties}(c)
and Lemma~\ref{propiedadesPCA}(b) we get that
$T_m\left(\setfibpt{\A}{\istar}\right) = \setfibpt{\A}{\sstar{i+1}},$
and
$T_k\evalat{\setfibpt{\A}{\istar}} = T_m\evalat{\setfibpt{\A}{\istar}}$
for every $k,m \in \Z^+,\ k\ge m$ and $i \in \DS.$
Consequently, by the definition of the map $T$ we have,
$T(\setfibpt{\A}{\istar}) = \setfibpt{\A}{\sstar{i+1}}$
for every $i\in\Z$ or, equivalently,
$T\bigl(\setfib{\A}{\Orbom}\bigr) = \setfib{\A}{\Orbom}.$

Now we consider $\setfibth{\A}$ with
$\theta \in \SI\setminus\Orbom.$
Since $\Orbom$ is dense in $\SI,$ there exists a sequence
$\{(\theta_{n}, x_{n})\}_{n=0}^{\infty} \subset \setfib{\A}{\Orbom}$
such that
$
\lim_{n\to\infty} \theta_{n} = \theta.
$
By the compacity of $\A$ we can assume (by taking a convergent subsequence, if necessary)
that $\{(\theta_{n}, x_{n})\}_{n=0}^{\infty}$ is convergent to a point $(\theta,x) \in \A.$
By Lemma~\ref{propiedadesPCA}(c),
$\setfibth{\A} = (\theta,x)$ (and $x = \gamma(\theta)$).
On the other hand, by the part of the statement already proven,
$T(\theta_{n}, x_{n}) \in \A$ for every $n.$
Hence, by the continuity of $T$ and the compacity of $\A,$
\[
T(\theta,x) = (R_\omega(\theta), f(\theta,x))
            = \lim_{n\to\infty} T(\theta_{n}, x_{n})
            \in \setfibpt{\A}{R_\omega(\theta)}.
\]
Since $\theta \notin \Orbom$
we have that $R_\omega(\theta) \notin \Orbom$ and,
again by Lemma~\ref{propiedadesPCA}(c),
$\setfibpt{\A}{R_\omega(\theta)}$ consists of a unique point.
Hence,  $T(\setfibth{\A}) = \setfibpt{\A}{R_\omega(\theta)}$
for every $\theta \in \SI\setminus\Orbom.$
Equivalently,
$
T\Bigl(\setfib{\A}{\bigl(\SI\setminus\Orbom\bigr)}\Bigr) =
\setfib{\A}{\bigl(\SI\setminus\Orbom\bigr)}.
$
This ends the proof of the theorem.
\end{proof}

\section{Proof of Proposition~\ref{T0mapProperties} in the case $m=0$}\label{proofofT0mapProperties}

This section is devoted to prove Proposition~\ref{T0mapProperties} for $m=0$;
that is, to study the map $T_0$.
It is the first technical counterpart of Section~\ref{skew-product}.

To prove Proposition~\ref{T0mapProperties} for $T_0$ we will need some
more notation and a technical lemma.

Given a skew product
$F(\theta,x) = (R_\omega(\theta), \zeta(\theta,x)$
from $\Omega = \SI \times \I$ to itself
we define the \emph{fibre map function of $F,$}
{\map{\mathsf{fib}(F)}{\SI}[\mathcal{C}(\I,\I)]}
by $\mathsf{fib}(F)(\theta) := \zeta(\theta, \cdot).$
A simple exercise shows that $F$ is continuous if and only if
$\zeta(\theta, \cdot)$ is continuous
for every $\theta \in \SI,$ and $\mathsf{fib}(F)$ is continuous.

\begin{lemma}\label{f0alesvores}
Let $\theta \in \Bd\left(\wbasint{i}\right)$ for some $i \in \DS[0].$
Then, $m_i(\theta) = M_i(\theta) = 0,$
$g_i(\theta,m_i(\theta)) = \gamma(R_\omega(\theta)),$
and $f_{0,\theta}$ is the unique piecewise affine map
with two affine pieces whose graph joins the point
$(-2,2)$ with $(0, \gamma(R_\omega(\theta))),$
and this with the point $(2, -2)$.
\end{lemma}

\begin{proof}
By Lemma~\ref{voresdelescaixesalesvores}(d) and Definition~\ref{T0mapDefi},
we have $m_i(\theta) = M_i(\theta).$
Hence, $f_{0,\theta}$ is the piecewise affine map with two affine pieces
whose graph joins the point $(-2,2)$ with
$\left(m_{i}(\theta), g_{_{i,\theta}}\left(m_{i}(\theta)\right)\right),$
and this with the point $(2, -2)$.
So, we need to show that $m_{i}(\theta) = 0,$ and
$g_{_{i,\theta}}\left(m_{i}(\theta)\right) = \gamma(R_\omega(\theta)).$

Lemma~\ref{Propertiesvarphi}(g) and the fact that $\dep{i} = 0,$
$\wbasint{i} \cap \wbasint{\ell} = \emptyset$ for every $\ell\in Z_{\ai},$
$i \ne \ell.$
Consequently, by Definition~\ref{PCgenerators}(\tsfR.6),
$m_i(\theta) = M_i(\theta) = a_i^- = 0.$

Now we show that
$g_{_{i,\theta}}(m_i(\theta)) = \gamma(R_\omega(\theta)).$
From the definition of the map $g_i$
(Definitions~\ref{defi-gi-positiva} and \ref{defi-gi-negativa}),
Lemma~\ref{Propertiesvarphi}(e) and
Definitions~\ref{gammalimit} and \ref{PCgenerators}(\tsfR.1),
we get
\[
g_{_{i,\theta}}(m_i(\theta))
  = \gams{\aii}(R_\omega(\theta))
  = \gamma(R_\omega(\theta)).
\]
This ends the proof of the lemma.
\end{proof}

\begin{proof}[Proof of Proposition~\ref{T0mapProperties} for $m=0$]
By Lemma~\ref{Propertiesvarphi}(b),
\[ -1 \le m_{\bt[0]{\theta}}(\theta) \le M_{\bt[0]{\theta}}(\theta) \le 1 \]
for every $\theta \in \wIBD[0].$
So, $T_0$ is well defined.

\inidemopart{b}
If $\theta \in \SI\setminus\wIBD[0],$ then the statement follows directly from
Definition~\ref{T0mapDefi}.
Now assume that $\theta \in \wIBD[0]$ and let $i = \bt[0]{\theta}.$
From the definition of the maps $g_{i,\theta}$
(Definitions~\ref{defi-gi-positiva} and \ref{defi-gi-negativa})
and Definition~\ref{T0mapDefi}, it follows that
$f_{0, \theta}\evalat{\I_{i,\theta}}$ is piecewise affine and non-increasing.
On the other hand, again by Definition~\ref{T0mapDefi},
$f_{0, \theta}\evalat{[-2, m_i(\theta)]}$ and $f_{0, \theta}\evalat{[M_i(\theta), 2]}$
are affine with negative slope and $f_{0,\theta}(2) = -2$ and $f_{0,\theta}(-2) = 2.$
The fact that
\[
   -1 \le f_{0,\theta}\left(M_{\bt[0]{\theta}}(\theta)\right) \le
          f_{0,\theta}\left(m_{\bt[0]{\theta}}(\theta)\right) \le 1
\]
for every $\theta \in \wIBD[0]$ follows from Definition~\ref{T0mapDefi}
and Lemmas~\ref{gpositiva}(a) and \ref{gnegativa}(a).
This ends the proof of (b).

\inidemopart{c}
Recall that
\[
 \wbasicbox{i} =
    \LSleftlimits{\bigcup}{\theta \in \wbasint{i}} \{\theta\} \times \I_{i,\theta}.
\]
Hence, from Definition~\ref{T0mapDefi} and the definition of $G_i$
(Definitions~\ref{defi-gi-positiva} and \ref{defi-gi-negativa})
it follows that
\[
T_m(\theta,x) = \bigl(R_\omega(\theta), f_m(\theta, x)\bigr)
              = \bigl(R_\omega(\theta),  g_{_{i,\theta}}(x)\bigr)
              = G_i(\theta, x),
\]
for every $(\theta,x) \in \wbasicbox{i}.$ Thus,
$
T_0\left(\setfibpt{\A_{\ai}}{\istar}\right)
  = \setfibpt{\A_{\aii}}{\sstar{i+1}}
$
from Lemmas~\ref{propiedadesA}(b), \ref{gpositiva}(c) and \ref{gnegativa}(c).
On the other hand, Lemma~\ref{denso}(b) implies that
$\istar \in \wIBD[0]$ but $\istar \notin \wIBD[k]$
for every $k \in \N.$
Then, we get $f_{k,\istar} = f_{0,\istar}$
from Definition~\ref{seqTmDefi}.

\inidemopart{a}
Since $T_0$ is a skew product with base $R_\omega$
we only have to prove that $f_0$ is continuous.

By Definition~\ref{T0mapDefi}, for every $\theta \in \SI$,
the map $f_{0,\theta}$ is continuous.
So we have to prove that the map
$\mathsf{fib}(T_0)$ (that is, the map $s \mapsto f_{0,s}$)
is continuous.

In the rest of the proof we will denote
\[
\IndSetWWings{IB}{0} := \bigcup_{i \in \DS[0]} \wobasint{i} \subset \wIBD[0].
\]

Clearly, since for every $i \in \Z,$
the maps $m_i$ and $M_i$ are continuous
on $\wbasint{i},$ it follows that the map
$s \mapsto f_{0,s}$ is continuous on $\IndSetWWings{IB}{0}.$
Thus, we have to see that the fibre map function
is continuous at every
$\theta \in \SI \setminus \IndSetWWings{IB}{0};$
that is,
$\lim_{j\to\infty} f_{0,\theta_j} = f_{0,\theta}$
for every
$\{\theta_j\}_{j=1}^{\infty} \subset \SI$ converging to $\theta.$
Given $\alpha > 0,$
we can consider four sets associated to such a sequence:
\begin{align*}
& \set{j\in \N}{\theta_j \in \SI \setminus \IndSetWWings{IB}{0}},\quad
  \set{j\in \N}{\theta_j \in \IndSetWWings{IB}{0} \setminus \ball{\theta}{\alpha}},\\
& \set{j\in \N}{\theta_j \in (\theta, \theta+\alpha) \cap \IndSetWWings{IB}{0}} \andq
  \set{j\in \N}{\theta_j \in (\theta-\alpha, \theta) \cap \IndSetWWings{IB}{0}}.
\end{align*}
Observe that the second set
$\set{j\in \N}{\theta_j \in \IndSetWWings{IB}{0} \setminus \ball{\theta}{\alpha}}$
is always finite and that any of the other three sets gives rise to
a subsequence of $\{\theta_j\}_{j=1}^{\infty}$ converging to $\theta,$
when it is infinite.
Consequently, the continuity of the fibre map function $s \mapsto f_{0,s}$
at $\theta$ is equivalent to the fact that
$\lim_{j\to\infty} f_{0,\theta_j} = f_{0,\theta}$
for every
$\{\theta_j\}_{j=1}^{\infty}$ converging to $\theta$
and such that, for some $\alpha > 0,$
$\{\theta_j\}_{j=1}^{\infty}$ is contained either in
$\SI\setminus\IndSetWWings{IB}{0},$ or
$(\theta, \theta+\alpha) \cap \IndSetWWings{IB}{0},$ or
$(\theta-\alpha, \theta) \cap \IndSetWWings{IB}{0}.$
We will only deal with the first two cases
since the proof in the last case (for $(\theta-\alpha, \theta)$)
can be done symmetrically.

\begin{case}{Case 1:}
$\lim_{j\to\infty} \theta_j = \theta$ and
$\{\theta_j\}_{j=1}^{\infty} \subset \SI\setminus\IndSetWWings{IB}{0}.$
\end{case}
By Definition~\ref{T0mapDefi} and Lemma~\ref{f0alesvores},
$f_{0,\theta_j}$ (respectively $f_{0,\theta}$)
is the unique piecewise affine map with two affine pieces
whose graph joins the point
$(-2,2)$ with $(0, \gamma(R_\omega(\theta_j)))$
(respectively $(0, \gamma(R_\omega(\theta)))$),
and this with the point $(2, -2)$.
By Lemma~\ref{denso}(c) and Definition~\ref{gammalimit}
the function $\gamma$ is continuous at
$R_\omega(\theta) \notin \Orbom.$
Hence,
$\lim_{j\to\infty} \gamma(R_\omega(\theta_j)) = \gamma(R_\omega(\theta))$
and, thus,
$\lim_{j\to\infty} f_{0,\theta_j} = f_{0,\theta}.$

\begin{case}{Case 2:}
$\lim_{j\to\infty} \theta_j = \theta$ and
$\{\theta_j\}_{j=1}^{\infty} \subset (\theta, \theta+\alpha) \cap \IndSetWWings{IB}{0}.$
\end{case}

If there exists $i \in \DS[0]$ such that
$\theta$ is the left endpoint of $\wbasint{i} \subset \wIBD[0]$
then the result follows from Definition~\ref{T0mapDefi}, the
continuity of the maps $m_i$ and $M_i$ and
the continuity of the maps $g_i$
(Lemmas~\ref{gpositiva}(a) and \ref{gnegativa}(a)).

Assume now that $\theta$ is not the left endpoint of $\wobasint{i}$
for every $i \in \DS[0].$
For every $j \in \N$ we set $\is_j := \bt[0]{\theta_j} \in \DS[0]$
(that is, $\theta_j \in \wobasint{\is_j}$).

We claim that $\lim_{j\to\infty} \abs{\is_j} = \infty$
and consequently, by Definition~\ref{PCgenerators}(\tsfR.1),
\begin{equation}\label{limit-ns-zero}
   \lim_{j\to\infty} 2^{-n_{\abs{\is_j + 1}}} =
   \lim_{j\to\infty} 2^{-n_{\abs{\is_j}}} = 0.
\end{equation}
To prove this claim, assume by way of contradiction that
there exists $L$ such that
for every $k \in \N$ there exists $j_k \ge k$
such that $\abs{\is_{j_k}} \le L.$
Then,
\[
\{\theta_{j_k}\}_{k=1}^{\infty} \subset
    \bigcup_{k=1}^\infty \wobasint{\is_{j_k}}
\]
and, since $\set{\is_{j_k}}{k\in \N}$ is finite, it follows that
there exists $i \in \set{\is_{j_k}}{k\in \N} \subset \DS[0]$
and a subsequence of $\{\theta_{j_k}\}_{k=1}^{\infty},$
that by abuse of notation will also be called
$\{\theta_{j_k}\},$
such that
$\{\theta_{j_k}\}_{k=1}^{\infty} \subset \wobasint{i}.$
So,
\[
  \theta = \lim_{k\to\infty} \theta_{j_k} \in \wbasint{i};
\]
a contradiction. So, the claim (and hence \eqref{limit-ns-zero}) holds.

Next we claim that the conditions
\begin{align}
& \lim_{j\to\infty} M_{\is_j}(\theta_j)
     = \lim_{j\to\infty} m_{\is_j}(\theta_j) = 0,\text{ and}\label{ElqueCalDemostrar1}\\
&\text{\parbox{0.9\textwidth}{there exists a sequence $\{x_j\}_{j=1}^{\infty}$ with
     $x_j \in \I_{\is_j,\theta_j} = [m_{\is_j}(\theta_j), M_{\is_j}(\theta_j)]$
     for every $j,$ such that $\lim\limits_{j\to\infty} f_{0,\theta_j}(x_j) = \gamma(R_\omega(\theta))$}}\label{ElqueCalDemostrar2}
\end{align}
imply
\[ \lim_{j\to\infty} f_{0,\theta_j} = f_{0,\theta}. \]

To prove the claim notice that,
by Definition~\ref{T0mapDefi} and Lemma~\ref{f0alesvores},
$f_{0,\theta}$ is the unique piecewise affine map
with two affine pieces whose graph joins the point
$(-2,2)$ with $(0, \gamma(R_\omega(\theta))),$
and this with the point $(2, -2).$
On the other hand, for every $j,$
\begin{itemize}
\item $f_{0,\theta_j}\evalat{[-2, m_{\is_j}(\theta_j)]}$
is the affine map joining the point $(-2,2)$ with the point
$(m_{\is_j}(\theta_j), g_{\is_j}(\theta_j, m_{\is_j}(\theta_j))),$
and
\item $f_{0,\theta_j}\evalat{[M_{\is_j}(\theta_j),2]}$
is the affine map joining the point
$(M_{\is_j}(\theta_j), g_{\is_j}(\theta_j, M_{\is_j}(\theta_j)))$
with the point$(2,-2)$
\end{itemize}
(see Figure~\ref{fig:mapsF0}).
Moreover, from the part of the proposition already proven we know that
$f_{0,\theta_j}$ is non-increasing and continuous.
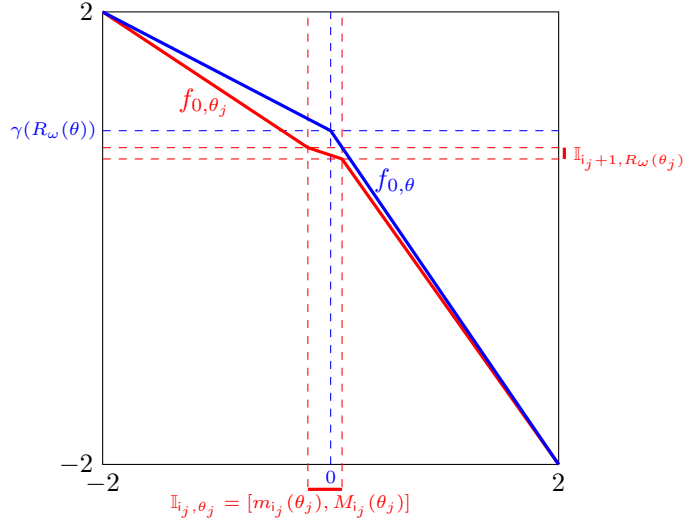
\begin{figure}
\begin{tikzpicture}[scale=1.5]
\draw (-2,-2) rectangle (2,2);
\foreach \c in {-2, 2} { \node[below] at (\c,-2) {$\c$}; \node[left] at (-2,\c) {$\c$}; }

\draw[dashed, color=blue] (0,-2) -- (0,2); \node[below, color=blue] at (0,-1.98) {\scriptsize$0$};
\draw[dashed, color=blue] (-2, 0.95) -- (2, 0.95);
\node[left, color=blue] at (-1.97,0.95) {\scriptsize$\gamma(R_\omega(\theta))$};

\draw[dashed, color=red] (-0.2,-2.2) -- (-0.2,2); \draw[dashed, color=red] (0.1,-2.2) -- (0.1,2);
\draw[very thick, color=red] (-0.2,-2.22) -- (0.1,-2.22);
\node[below, color=red] at (-0.35,-2.18) {\scriptsize$\I_{\is_j, \theta_j} = [m_{\is_j}(\theta_j),M_{\is_j}(\theta_j)]$};

\draw[dashed, color=red] (-2, 0.8) -- (2.05, 0.8); \draw[dashed, color=red] (-2, 0.7) -- (2.05, 0.7);
\draw[very thick, color=red] (2.05,0.8) -- (2.05, 0.7);
\node[right, color=red] at (2.05,0.7) {\scriptsize$\I_{\is_j+1,R_\omega(\theta_j)}$};

\draw[very thick, color=red] (-2,2) -- (-0.2, 0.8) -- (0.1, 0.7) -- (2,-2);
\node[color=red, left] at (-0.8,1.2) {$f_{0,\theta_j}$};
\draw[very thick, color=blue] (-2,2) -- (0, 0.95) -- (2,-2);
\node[color=blue, right] at (0.3,0.5) {$f_{0,\theta}$};
\end{tikzpicture}
\caption{A symbolic representation of the maps
$f_{0,\theta}$ and $f_{0,\theta_j}$ in Case~2 of
the proof of Proof of Proposition~\ref{T0mapProperties} for $m=0$.
The map \textcolor{blue}{$f_{0,\theta}$} and the points
\textcolor{blue}{$0$} and \textcolor{blue}{$\gamma(R_\omega(\theta))$}
are drawn in \textcolor{blue}{blue}.
The map \textcolor{red}{$f_{0,\theta_j}$} and the corresponding intervals
\textcolor{red}{$\I_{\is_j, \theta_j} $} and \textcolor{red}{$\I_{\is_j+1,R_\omega(\theta_j)}$}
are drawn in \textcolor{red}{red}.}\label{fig:mapsF0}
\end{figure}
Therefore, the claim holds provided that
\[
 \lim_{j\to\infty} \diam\left(f_{0,\theta_j}\left(\I_{\is_j,\theta_j}\right)\right) = 0
\]
(see again Figure~\ref{fig:mapsF0}).

When $\theta_j \in \basint{\is_j} \setminus \OBG{\is_j}{\is_j+1}$ and $\is_j \ge 0,$
by Definitions~\ref{T0mapDefi} and \ref{defi-gi-positiva},
\[
   \diam\left(f_{0,\theta_j}\left(\I_{\is_j,\theta_j}\right)\right) =
   \diam\left(g_{_{\is_j,\theta_j}}\left(\I_{\is_j,\theta_j}\right)\right) =
   \diam\left(\{\gams{\is_j+1}(R_\omega(\theta_j)\}\right) = 0.
\]
Otherwise, by Definition~\ref{T0mapDefi}, and
Lemmas~\ref{gpositiva}(b) and \ref{gnegativa}(b),
\begin{align*}
\{R_\omega(\theta_j)\} \times f_{0,\theta_j}\left(\I_{\is_j,\theta_j}\right)
  & = \{R_\omega(\theta_j)\} \times g_{_{\is_j,\theta_j}}\left(\I_{\is_j,\theta_j}\right)
    = G_{_{\is_j}}\bigl(\setfibpt{\basicbox{\is_j}}{\theta_j}\bigr) \\
  & \subset \setfibpt{\basicbox{\is_j + 1}}{R_\omega(\theta_j)}.
\end{align*}
So, by Remark~\ref{propiedadesR}(2),
\[
   \diam\left(f_{0,\theta_j}\left(\I_{\is_j,\theta_j}\right)\right) \le
   \diam\left(\basicbox{\is_j + 1}\right) \le
   2 \cdot 2^{-n_{\abs{\is_j + 1}}}.
\]
In any case,
\[
   0 \le \diam\left(f_{0,\theta_j}\left(\I_{\is_j,\theta_j}\right)\right) \le
   2 \cdot 2^{-n_{\abs{\is_j + 1}}}
   \andq[for every] j \in \N
\]
and, by \eqref{limit-ns-zero},
$
\lim_{j\to\infty} \diam\left(f_{0,\theta_j}\left(\I_{\is_j,\theta_j}\right)\right) = 0.
$
This ends the proof of the claim.

By the last claim, to end the proof of the proposition
in the case $m = 0$ it is enough to show that \mbox{(\ref{ElqueCalDemostrar1}--\ref{ElqueCalDemostrar2})} hold.
We start by proving \eqref{ElqueCalDemostrar1}.
By Lemma~\ref{f0alesvores},
\[
  m_{\is_j}(\Bd(\wbasint{\is_j})) = M_{\is_j}(\Bd(\wbasint{\is_j})) = 0,
\]
and from the definition of the maps $m_{\is_j}$ and $M_{\is_j},$
Definition~\ref{GenericBoxes} (or Lemma~\ref{voresdelescaixesalesvores})
and Remark~\ref{propiedadesR}(2),
for every $s \in \wobasint{\is_j}$ we get
\begin{equation}\label{claimfibres}
\begin{split}
& -1 \le m_{\is_j}(s) < 0 < M_{\is_j}(s) \le 1,\text{ and}\\
& M_{\is_j}(s) - m_{\is_j}(s)
  = \diam\bigl(\I_{{\is_j},s}\bigr)
  \le 2\cdot2^{-n_{\abs{\is_j}}}.
\end{split}
\end{equation}
So, \eqref{ElqueCalDemostrar1} holds by \eqref{limit-ns-zero}.
Now we prove \eqref{ElqueCalDemostrar2}.

By \eqref{gammathetaproperties},
\eqref{alphadeltaintervals} and
\eqref{alphadeltaintervalsineg},
it follows that
\begin{align*}
 & m_{\is_j}(\theta_j)
   < \gams{\abs{\is_j}}(\theta_j) < M_{\is_j}(\theta_j)
      && \text{if $\theta_j \ne \sstar{\is_j},$ and}\\
 & m_{\is_j}(\theta_j)
   < \gams{\abs{\is_j}-1}(\theta_j) = 0 < M_{\is_j}(\theta_j)
      && \text{if $\theta_j = \sstar{\is_j}.$}
\end{align*}
Also, from Definition~\ref{T0mapDefi},
the definitions of $G_i$ and $g_{_{i,\theta}}$
(Definitions~\ref{defi-gi-positiva} and \ref{defi-gi-negativa}),
and Lemmas~\ref{gpositiva}(c) and \ref{gnegativa}(c)
we get
\begin{align*}
 & f_{0,\theta_j}(\gams{\abs{\is_j}}(\theta_j))
   = g_{_{\is_j,\theta_j}}(\gams{\abs{\is_j}}(\theta_j))
   = \gams{\abs{\is_j+1}}(R_\omega(\theta_j))
   && \text{if $\theta_j \ne \sstar{\is_j},$}\\
 & f_{0,\theta_j}(\gams{\is_j-1}(\theta_j))
   = g_{_{\is_j,\theta_j}}(\gams{\is_j-1}(\theta_j))
   = \gams{\is_j}(R_\omega(\theta_j))
   && \text{if $\theta_j = \sstar{\is_j}$ and $\is_j \ge 0,$ and}\\
 & f_{0,\theta_j}(\gams{\abs{\is_j}-1}(\theta_j))
   = g_{_{\is_j,\theta_j}}(\gams{\abs{\is_j+1}}(\theta_j))
   = \gams{\abs{\is_j+2}}(R_\omega(\theta_j))
   && \text{if $\theta_j = \sstar{\is_j}$ and $\is_j < 0.$}
\end{align*}
Thus, to prove \eqref{ElqueCalDemostrar2}, we have to show that
\begin{equation}\label{thelimits}
\begin{cases}
 \lim_{j\to\infty} \gams{\abs{\is_j+1}}(R_\omega(\theta_j))
   = \gamma(R_\omega(\theta))
   &\text{if $\theta_j \ne \sstar{\is_j},$}\\
 \lim_{j\to\infty} \gams{\is_j}(R_\omega(\theta_j))
   = \gamma(R_\omega(\theta))
   &\text{if $\theta_j = \sstar{\is_j}$ and $\is_j \ge 0,$ and}\\
 \lim_{j\to\infty} \gams{\abs{\is_j+2}}(R_\omega(\theta_j))
   = \gamma(R_\omega(\theta))
   &\text{if $\theta_j = \sstar{\is_j}$ and $\is_j < 0$}
\end{cases}
\end{equation}
(that is, we take
$x_j := \gams{\abs{\is_j}}(\theta_j)$ if $\theta_j \ne \sstar{\is_j},$
$x_j := \gams{\is_j-1}(\theta_j)$ if $\theta_j = \sstar{\is_j}$ and $\is_j \ge 0,$ and
$x_j := \gams{\abs{\is_j}-1}(\theta_j)$ if $\theta_j = \sstar{\is_j}$ and $\is_j < 0$).

Let $\varepsilon > 0.$
By Lemma~\ref{denso}(c) and Definition~\ref{PCgenerators}(\tsfR.1)
we have that $\theta \notin\Orbom$ and, hence, $R_\omega(\theta) \notin \Orbom.$
By the continuity of $\gamma$ on $\SI\setminus\Orbom$
and the fact that $\lim_{i\to\infty}\gams{i} = \gamma,$
there exist $\delta > 0$ and $L \in \N$
such that
\begin{align*}
 & \abs{\gamma(R_\omega(\theta)) - \gamma(\widehat{\theta})} < \varepsilon/2
   \quad\text{for every
      $\widehat{\theta} \in \ball{R_\omega(\theta)}{\delta} \setminus \Orbom$,
   and}\\
 & d_{\infty}\left(\gamma,\gams{i}\right) < \varepsilon/2
   \quad\text{for every $i \ge L.$}
\end{align*}
Then, since
$\lim_{j\to\infty} \theta_j = \theta$ and
$\lim_{j\to\infty} \abs{\is_j} = \infty,$
there exists $N \in \N$ such that
$\abs{\theta - \theta_j} < \delta/2,$
and $\abs{\is_j} \ge L+2$
for every $j \ge N.$

First we will show that
\[
 \abs{\gamma(R_\omega(\theta)) - \gams{\abs{\is_j+1}}(R_\omega(\theta_j))}
    \le \varepsilon
\]
for every $j \ge N$ such that  $\theta_j \ne \sstar{\is_j}.$
To see it observe that,
by Definition~\ref{PCgenerators}(\tsfR.2) and
Remark~\ref{PCgeneratorsExplicitConsequences}(\tsfR.2),
$\theta_j, R_\omega(\theta_j) \notin \Zstar_{\abs{\is_j+1}}$
whenever $\theta_j \ne \sstar{\is_j}.$
Thus,
$\gams{\abs{\is_j+1}}$ is continuous at $R_\omega(\theta_j)$
by Lemma~\ref{Propertiesvarphi}(a).

Also, there exists a sequence
$
\{\widehat{\theta}_{j_{\ell}}\}_{\ell=1}^\infty \subset
\bigl(\ball{\theta_j}{\delta/2} \cap \wobasint{\is_j}\bigr) \setminus \Orbom
$
converging to $\theta_j,$
because $\SI \setminus \Orbom$ is dense in $\SI.$
Clearly, for every $j \ge N,$ we have
$
  \{R_\omega(\widehat{\theta}_{j_{\ell}})\}_{\ell=1}^\infty \subset
  \ball{R_\omega(\theta)}{\delta} \setminus \Orbom
$
and
$
\lim_{\ell\to\infty} R_\omega(\widehat{\theta}_{j_{\ell}})
  = R_\omega(\theta_j).
$
Moreover, since
$
\{R_\omega(\widehat{\theta}_{j_{\ell}})\}_{\ell=1}^\infty
  \subset \SI \setminus \Orbom \subset \SI \setminus \Zstar_{\abs{\is_j+1}},\
\gams{\abs{\is_j+1}}
$
is defined for every $R_\omega(\widehat{\theta}_{j_{\ell}}).$
Then, for every $j \ge N$ and $\ell \in \N$, we have
\begin{align*}
\abs{\gamma(R_\omega(\theta)) -
     \gams{\abs{\is_j+1}}(R_\omega(\widehat{\theta}_{j_{\ell}}))
}\ \le\
    & \abs{\gamma(R_\omega(\theta)) -
           \gamma(R_\omega(\widehat{\theta}_{j_{\ell}}))
      } + \\
    & \hspace*{5em}
      \abs{\gamma(R_\omega(\widehat{\theta}_{j_{\ell}})) -
           \gams{\abs{\is_j+1}}(R_\omega(\widehat{\theta}_{j_{\ell}}))
      }\\
<\  & \tfrac{\varepsilon}{2} + d_{\infty}\Bigl(\gamma,\gams{\abs{\is_j+1}}\Bigr)
      < \varepsilon.
\end{align*}
Consequently,
\[
 \abs{\gamma(R_\omega(\theta)) - \gams{\abs{\is_j+1}}(R_\omega(\theta_j))}
   = \lim_{\ell\to\infty}
           \abs{\gamma(R_\omega(\theta)) -
                \gams{\abs{\is_j+1}}(R_\omega(\widehat{\theta}_{j_{\ell}}))
           } \le \varepsilon
\]
This ends the proof of the first equality of \eqref{thelimits}.
The second and third equalities of \eqref{thelimits} follow as above
by replacing $\gams{\abs{\is_j+1}}$ by
$\gams{\is_j}$ (respectively $\gams{\abs{\is_j+2}}$),
and noting that
\[
R_\omega(\theta_j) = R_\omega(\sstar{\is_j}) =
 \begin{cases}
    \sstar{(\is_j+1)} \notin \Zstar_{\is_j}
       & \text{ if $\is_j \ge 0,$ and}\\
    \sstar{(-(\abs{\is_j}-1))} \notin \Zstar_{\abs{\is_j}-2}
       & \text{ if $\is_j < 0.$}
 \end{cases}
\]
This ends the proof of the continuity of $T_0,$
and the proposition for the case $m=0.$
\end{proof}

\section{Proof of Proposition~\ref{seqTmProperties} for $m > 0$}\label{proofofseqTmProperties}

This section is the second technical counterpart of Section~\ref{skew-product}
and is devoted to prove Proposition~\ref{seqTmProperties} for every map $T_m$ with $m > 0.$
To do this we will need some more technical results.
Also we will use the notion of fibre map function introduced
in the previous section.

The next two lemmas establish some basic properties
of the maps $T_m\evalat{\wIVD}$ and clarify some aspects
of Definition~\ref{seqTmDefi}.

\begin{lemma}\label{seqTmDefirevisited}
For every $m \in \N$ and for every $\theta \in \wIBD,$
\[
  f_{m,\theta}\evalat{\I_{\is,\theta}} =
  g_{_{\is,\theta}}\evalat{\I_{\is,\theta}},
\]
where $\is = \bt{\theta}.$
Moreover, assume that
$
\theta \in \WB \setminus \WIB.
$
Then,
\[
f_{m,\theta}(x) = \begin{cases}
   g_{_{\is,\theta}}(x)
        & \text{if $x \in \I_{\is,\theta}$}, \\[0.75ex]
   \frac{2 - g_{_{\is,\theta}}\left(m_{\is}(\theta)\right)\hfill}{
         2 - f_{m-1,\theta}\left(m_{\is}(\theta)\right)
        } (f_{m-1,\theta}(x) - 2) + 2
        & \text{if $x \in [-2,m_{\is}(\theta)]$},\\[1ex]
   \frac{2 + g_{_{\is,\theta}}\left(M_{\is}(\theta)\right)\hfill}{
         2 + f_{m-1,\theta}\left(M_{\is}(\theta)\right)
        } (f_{m-1,\theta}(x) + 2) - 2
        & \text{if $x \in [M_{\is}(\theta),2]$}.
\end{cases}
\]
\end{lemma}

\begin{proof}
We start by proving the first statement.
When $\theta\in\IBD$ there is nothing to prove.
So, assume that $\theta \in \wIBD \setminus \IBD.$
By Definition~\ref{curvesinthewings}, $\theta \in \WB,$
$\is < 0$ and $\theta \in \wbasint{\is} \setminus \obasintabs{\is}.$
By Lemma~\ref{VerticalIntervalsIntheWings}(b),
\[
  \I_{\is, \theta} = \{\gams{\abs{\is}}(\theta)\} \subset \IW{\theta}.
\]
Consequently, by Definition~\ref{seqTmDefi} and the definition of
the maps $g_{_{\is,\theta}}$ for $\is < 0$
(Definition~\ref{defi-gi-negativa} ---
notice that $\I_{\is, \theta} \subset \wbasicbox{\is}$ by definition),
\[
 f_{m,\theta}\left(\gams{\abs{\is}}(\theta)\right) =
     \gams{\abs{\is+1}}\left(R_\omega(\theta)\right) =
     g_{_{\is,\theta}}\left(\gams{\abs{\is}}(\theta)\right).
\]
So, the first statement holds.
Now we prove the second one.
By Lemma~\ref{VerticalIntervalsIntheWings}(b),
\[
    \I_{\is, \theta} =
    \{m_{\is}(\theta)\} = \{M_{\is}(\theta)\} =
    \{\gams{\abs{\is}}(\theta)\} =
    \{\lambda_m(\theta)\} = \{\tau_m(\theta)\} =
    \IW{\theta}.
\]
Thus, by the part already proven, the formulas
\begin{align*}
&\begin{cases}
   g_{_{\is,\theta}}(x)
        & \text{if $x \in \I_{\is,\theta}$}, \\[0.75ex]
   \frac{2 - g_{_{\is,\theta}}\left(m_{\is}(\theta)\right)\hfill}{
         2 - f_{m-1,\theta}\left(m_{\is}(\theta)\right)
        } (f_{m-1,\theta}(x) - 2) + 2
        & \text{if $x \in [-2,m_{\is}(\theta)]$},\\[1ex]
   \frac{2 + g_{_{\is,\theta}}\left(M_{\is}(\theta)\right)\hfill}{
         2 + f_{m-1,\theta}\left(M_{\is}(\theta)\right)
        } (f_{m-1,\theta}(x) + 2) - 2
        & \text{if $x \in [M_{\is}(\theta),2]$},
\end{cases}\\
\intertext{and}\\
&\begin{cases}
   \gams{\abs{\is+1}}\left(R_\omega(\theta)\right)
        & \text{if $x \in \IW{\theta}$}, \\[0.75ex]
   \frac{2 - \gams{\abs{\is+1}}\left(R_\omega(\theta)\right)\hfill}{
         2 - f_{m-1,\theta}\left(\lambda_m(\theta)\right)
        } (f_{m-1,\theta}(x) - 2) + 2
        & \text{if $x \in [-2,\lambda_m(\theta)]$},\\[1ex]
   \frac{2 + \gams{\abs{\is+1}}\left(R_\omega(\theta)\right)\hfill}{
         2 + f_{m-1,\theta}\left(\tau_m(\theta)\right)
        } (f_{m-1,\theta}(x) + 2) - 2
        & \text{if $x \in [\tau_m(\theta),2]$},
\end{cases}
\end{align*}
coincide.
\end{proof}

\begin{lemma}\label{seqTmPropsInBasicintervals}
The following statements hold for every $m \in \N$ and $i \in \DS:$
\begin{enumerate}[(a)]
\item The map $T_m\evalat{\wbasband{i}}$ is well defined and continuous.
\item For every $\theta\in\wbasint{i},$
\begin{enumerate}[(b.i)]
      \item $f_{m,\theta}(2) = -2$ and $f_{m,\theta}(-2) = 2,$
      \item $f_{m,\theta}$ is piecewise affine and non-increasing, and
      \item $-1 \le f_{m,\theta}\left(M_{i}(\theta)\right) \le
              f_{m,\theta}\left(m_{i}(\theta)\right) \le 1.$
\end{enumerate}
\item $T_m\evalat{\wbasicbox{i}} = G_i$ and
      $T_m\left(\setfibpt{\A_{\ai}}{\istar}\right) = \setfibpt{\A_{\aii}}{\sstar{i+1}}.$
\end{enumerate}
\end{lemma}

\begin{proof}
Clearly, $T_m\evalat{\wbasband{i}}$ is well defined and continuous if and only if
so is $f_m\evalat{\wbasband{i}}.$

We will prove by induction on $m \in \Z^+$ that,
(a), (b) and
\begin{enumerate}[(b.i)]\setcounter{enumi}{3}
\item $f_{m, \theta}\evalat{[-2, -1]}$ and $f_{m, \theta}\evalat{[1, 2]}$
      are affine, $f_{m, \theta}(-1) < 2$ and $f_{m, \theta}(1) > -2$
\end{enumerate}
hold for every $\theta\in\wbasint{i}.$

First we will show that (a), (b) and (b.iv) hold for $m=0$ and $i \in \DS[0]$
(we are including the map $f_0$ studied earlier to correctly
start the induction process).
By Proposition~\ref{T0mapProperties}(a,b) for $m=0$ we have that
$T_{0}\evalat{\wbasband{i}}$ is well defined and continuous and
(b) holds. By Definition~\ref{T0mapDefi},
we also know that
$f_{m, \theta}\evalat{[-2, m_{i}(\theta)]}$
and
$f_{m, \theta}\evalat{[M_{i}(\theta), 2]}$
are affine. Then, (b.iv) follows from
$
-1 \le m_{i}(\theta) \le M_{i}(\theta) \le 1
$
(see Lemma~\ref{voresdelescaixesalesvores}(a)) and (b.iii).

Assume now that (a), (b) and (b.iv) hold for some
$m-1 \in \Z +$ and prove it for $m$ and $i \in \DS.$
By Lemma~\ref{Dsets}(a),
$\theta \in \wbasint{i} \varsubsetneq \wbasint{k}$
for some $k \in \DS[m-1].$
Consequently, $\wbasband{i} \subset \wbasband{k}$
and $f_{m-1}\evalat{\wbasband{i}}$ is well defined and continuous.

By Lemma~\ref{voresdelescaixesalesvores}(a) and Definition~\ref{curvesinthewings},
\begin{equation}\label{colocacioms}
\begin{split}
  -1 \le m_{i}(\theta) \le M_{i}(\theta) \le 1&
     \qquad\text{for $\theta\in\wbasint{i}$, and}\\
  -1 \le \lambda_m(\theta) \le \tau_m(\theta) \le 1&
      \qquad\text{for $\theta\in\wbasint{i}\setminus\obasintabs{i} \subset \WB$ ($i < 0$).}
\end{split}
\end{equation}
Consequently, by (b.ii) and (b.iv) for $m-1,$
\[
 -2 < f_{m-1,\theta}(1) \le f_{m-1,\theta}\left(M_{i}(\theta)\right)
    \le f_{m-1,\theta}\left(m_{i}(\theta)\right) \le f_{m,\theta}(-1) < 2
\]
for every $\theta\in\wbasint{i}$, and
\[
-2 < f_{m-1,\theta}(1) \le f_{m-1,\theta}\left(\tau_m(\theta)\right)
    \le f_{m-1,\theta}\left(\lambda_m(\theta)\right) \le f_{m,\theta}(-1) < 2
\]
for $\theta\in\wbasint{i}\setminus\obasintabs{i} \subset \WB$ when $i < 0$.

On the other hand, as it was observed in Definition~\ref{seqTmDefi},
$f_{m,\theta}$ is defined in two different ways when
$\theta \in \WB \cap \IBD.$
In such a case,
by Lemmas~\ref{VerticalIntervalsIntheWings}(a,b)
and \ref{seqTmDefirevisited}, $\theta \notin \WIB$ and
both definitions for $f_{m,\theta}$ coincide.
Hence, $f_{m}\evalat{\wbasband{i}}$ is well defined.

Now we prove that $f_{m}\evalat{\wbasband{i}}$ is continuous
by using the continuity of $f_{m-1}\evalat{\wbasband{i}}.$
Since $\basintabs{i} \subset \IBD,$
by Definition~\ref{seqTmDefi},
the continuity of the maps $m_i$ and $M_i$ (see Lemma~\ref{voresdelescaixesalesvores}(b)),
and the continuity of the maps $g_i$ (Lemmas~\ref{gpositiva}(a) and \ref{gnegativa} (a)),
$f_m\evalat{\setsilift{\basintabs{i}}}$ is continuous.
Now we assume that $i < 0$ and we study the continuity of
$f_m\evalat{\setsilift{U}}$ on a connected component $U$ of
$\wbasint{i}\setminus\obasintabs{i}.$
Observe that, by Definition~\ref{curvesinthewings} and Lemma~\ref{Dsets}(b),
$U$ is a connected component of $\WB.$
Then, again by Definition~\ref{seqTmDefi},
the continuity of the maps $\lambda_m\evalat{U}$ and $\tau_m\evalat{U}$
(Lemma~\ref{continuouscurvesinthewings}),
and the continuity of the map $\gams{\ai}\evalat{U}$
(Lemma~\ref{Propertiesvarphi}(a) and
Definition~\ref{PCgenerators}(\tsfR.2) and
Remark~\ref{PCgeneratorsExplicitConsequences}(\tsfR.2)),
$f_m\evalat{\setsilift{U}}$ is continuous.
Therefore, $f_{m}\evalat{\wbasband{i}}$ is continuous because it is well defined
on $\setsilift{\left(\left(\wbasint{i}\setminus\obasintabs{i}\right) \cap \basintabs{i}\right)}.$

Let $\theta \in \basintabs{i} \subset \IBD$.
By Definition~\ref{seqTmDefi},
and the definition of the maps $g_{i,\theta}$
(Definitions~\ref{defi-gi-positiva} and \ref{defi-gi-negativa}),
$f_{m, \theta}\evalat{\I_{i,\theta}}$
is piecewise affine and non-increasing.
So, by Lemma~\ref{seqTmDefirevisited} for $m-1$ and Definition~\ref{seqTmDefi},
$f_{m,\theta}(2) = -2,$ $f_{m,\theta}(-2) = 2,$ and
$f_{m, \theta}\evalat{[-2, m_i(\theta)]}$ and
$f_{m, \theta}\evalat{[M_i(\theta), 2]}$
are affine transformations of the map $f_{m-1, \theta}$
with positive slope.
Hence, (b.i,ii) hold for $f_{m, \theta}$ in this case.
Moreover, (b.iv) is verified
by \eqref{colocacioms} and (b.iv) for $m-1.$

Consider $\theta \in \wbasint{i}\setminus\obasintabs{i} \subset \WB$.
Again by Definition~\ref{seqTmDefi},
$f_{m, \theta}\evalat{\IW{\theta}}$ is constant.
Then, (b.i,ii) and (b.iv) hold for $f_{m, \theta}$ as above
by replacing $m_i(\theta)$ and $M_i(\theta)$ by
$\lambda_m(\theta)$ and $\tau_m(\theta),$ respectively.

By (b.ii) and \eqref{colocacioms} we have
$
  f_{m,\theta}\left(M_{i}(\theta)\right) \le
  f_{m,\theta}\left(m_{i}(\theta)\right).
$
Hence, (b.iii) follows from Lemma~\ref{seqTmDefirevisited},
Definition~\ref{seqTmDefi},
Lemmas~\ref{gpositiva}(b) and \ref{propiedadesA}(c),
Definition~\ref{PCgenerators}(\tsfR.2) and
Remark~\ref{PCgeneratorsExplicitConsequences}(\tsfR.2),
Lemma~\ref{gnegativa}(b) and Lemma~\ref{Propertiesvarphi}(b).

\inidemopart{c}
In a similar way to the proof of Proposition~\ref{T0mapProperties} for the case $m = 0$,
\[
 \wbasicbox{i} =
    \LSleftlimits{\bigcup}{\theta \in \wbasint{i}} \{\theta\} \times \I_{i,\theta}
    \subset \wbasband{i} \subset \wIVD
\]
and, by Definition~\ref{seqTmDefi}, Lemma~\ref{seqTmDefirevisited}
and the definition of $G_i$
(Definitions~\ref{defi-gi-positiva} and \ref{defi-gi-negativa})
it follows that
\[
T_m(\theta,x) = \bigl(R_\omega(\theta), f_m(\theta, x)\bigr)
              = \bigl(R_\omega(\theta),  g_{_{i,\theta}}(x)\bigr)
              = G_i(\theta, x),
\]
for every $(\theta,x) \in \wbasicbox{i}.$
Thus,
$
T_m\left(\setfibpt{\A_{\ai}}{\istar}\right)
  = \setfibpt{\A_{\aii}}{\sstar{i+1}}
$
from Lemmas~\ref{propiedadesA}(b), \ref{gpositiva}(c) and \ref{gnegativa}(c).
\end{proof}

The next technical lemma compares the images of
$f_{m,\theta}$ and $f_{m-1,\theta}$ on a point.
It is an extension of Lemma~\ref{QuePassaALesAles}.

\begin{lemma}\label{fmfm-1alesales}
Assume that $\wbasint{i} \subset \wbasint{k}$
for some $i \in \DS,\ k \in \DS[m-1]$ and $m \in \N.$
Then, for every $\theta \in \wbasint{i} \setminus \obasintabs{i},$
$m_i(\theta) = M_i(\theta) = \gams{i}(\theta)$ and
\begin{align*}
f_{m,\theta}\left(m_i(\theta)\right) & = g_{_{i, \theta}}\left(m_i(\theta)\right) =
   \gams{\aii}\left(R_\omega(\theta)\right),\text{ and}\\
f_{m-1,\theta}\left(m_i(\theta)\right) & = g_{_{k, \theta}}\left(m_i(\theta)\right) =
   \gams{\akk}\left(R_\omega(\theta)\right).
\end{align*}
\end{lemma}

\begin{proof}
The fact that $m_i(\theta) = M_i(\theta) = \gams{i}(\theta)$ follows directly from the definitions.
The first equation follows from Lemma~\ref{seqTmDefirevisited},
and the definition of the map $g_{_{i, \theta}}$
(Definitions~\ref{defi-gi-positiva} and \ref{defi-gi-negativa}).

By Lemma~\ref{QuePassaALesAles},
$
\I_{i, \theta} = \{m_i(\theta)\} = \{\gams{\ak}(\theta)\} \subset \I_{k, \theta}.
$
Moreover, as in the proof of Lemma~\ref{QuePassaALesAles},
$\theta \ne \kstar.$
Consequently, by Definition~\ref{T0mapDefi}, Lemma~\ref{seqTmDefirevisited},
Lemmas~\ref{gpositiva}(c) and \ref{gnegativa}(c) and \eqref{gammathetaproperties}
(alternatively, for the last equality check directly the proofs
of the Lemmas~\ref{gpositiva}(c) and \ref{gnegativa}(c)),
\[
f_{m-1,\theta}\left(m_i(\theta)\right) =
   g_{_{k, \theta}}\left(m_i(\theta)\right) =
   g_{_{k, \theta}}\left(\gams{\ak}(\theta)\right) =
   \gams{\akk}\left(R_\omega(\theta)\right).
\]
\end{proof}

The following lemma is the analogue of Lemma~\ref{f0alesvores} for $m \ge 1.$
To state it we will use the set
\[
 \setsilift{\wEIBD} = \wEIBD \times \I \subset \wIVD.
\]

\begin{lemma}\label{voresiguals}
$
  T_m\evalat{\setsilift{\wEIBD}} =
  T_{m-1}\evalat{\setsilift{\wEIBD}}
$
for every $m \in \N.$
Equivalently, $f_{m,\theta} = f_{m-1,\theta}$
for every $m \in \N$ and $\theta \in \wEIBD.$
\end{lemma}

\begin{proof}
Fix $m \in \N$ and $\theta \in \wEIBD \subset \wIBD.$
By Lemma~\ref{Dsets}(a,b), there exist $i \in \DS$ and $k \in \DS[m-1]$ such that
$\theta \in \Bd\left(\wbasint{i}\right) \subset \wbasint{i} \varsubsetneq \wbasint{k}.$
So, we are in the assumptions of
Lemmas~\ref{QuePassaALesAles} and \ref{fmfm-1alesales}
and, hence,
\begin{align*}
\I_{i, \theta} & =\{m_i(\theta)\} = \{\gams{\ai}(\theta)\} = \{\gams{\ak}(\theta)\} \subset \I_{k, \theta},\\
f_{m,\theta}\left(m_i(\theta)\right) & = g_{_{i, \theta}}\left(m_i(\theta)\right) =
   \gams{\aii}\left(R_\omega(\theta)\right),\text{ and}\\
f_{m-1,\theta}\left(m_i(\theta)\right) & = g_{_{k, \theta}}\left(m_i(\theta)\right) =
   \gams{\akk}\left(R_\omega(\theta)\right).
\end{align*}

Thus, if $i \ge 0,$ $\theta \in \IBD$ and,
by Definition~\ref{seqTmDefi} and
Lemma~\ref{seqTmPropsInBasicintervals}(a),
to prove that $f_{m,\theta} = f_{m-1,\theta}$
we only have to show that
\[
g_{_{i, \theta}}\left(m_i(\theta)\right) =
   \gams{\aii}\left(R_\omega(\theta)\right) =
   \gams{\akk}\left(R_\omega(\theta)\right) =
   f_{m-1,\theta}\left(m_i(\theta)\right).
\]

When $i < 0,$ $\theta \in \WB \cap \wEIBD$
and, by Lemma~\ref{VerticalIntervalsIntheWings}(a),
$\theta \notin \WIB.$
Then, by Lemma~\ref{seqTmDefirevisited},
we get again that
\[
g_{_{i, \theta}}\left(m_i(\theta)\right) =
   \gams{\aii}\left(R_\omega(\theta)\right) =
   \gams{\akk}\left(R_\omega(\theta)\right) =
   f_{m-1,\theta}\left(m_i(\theta)\right).
\]
implies $f_{m,\theta} = f_{m-1,\theta}.$

If $\akk = \aii$ there is nothing to prove.
So, by Lemma~\ref{QuePassaALesAles}, we can assume that $\akk < \aii$
and we have
\[
 \gams{\akk}\left(R_\omega(\theta)\right) =
 \gams{\akk+1}\left(R_\omega(\theta)\right) = \dots =
 \gams{\aii-1}\left(R_\omega(\theta)\right).
\]
Hence, we have to show that
$
 \gams{\aii}\left(R_\omega(\theta)\right) = \gams{\aii-1}\left(R_\omega(\theta)\right).
$
If $i \ge 0$ we get
\[
 \gams{\aii}\left(R_\omega(\theta)\right) =
 \gams{i+1}\left(R_\omega(\theta)\right) =
 \gams{i}\left(R_\omega(\theta)\right) =
 \gams{\aii-1}\left(R_\omega(\theta)\right)
\]
by Lemma~\ref{Propertiesvarphi}(e).
Otherwise we have $i < 0,$
$\theta \in \Bd\left(\wbasint{i}\right) = \Bd\left(\BSG{i}{\aii}\right)$
and, consequently, $R_\omega(\theta) \in \Bd\left(\BSG{i+1}{\aii}\right).$
Again by Lemma~\ref{Propertiesvarphi}(e) for $j = \aii$,
\[
 \gams{\aii}\left(R_\omega(\theta)\right) =
 \gams{\aii-1}\left(R_\omega(\theta)\right).
\]
This ends the proof of the lemma.
\end{proof}

Now we aim at computing two different kind of
upper bounds for
$\norm{f_{m,\theta} - f_{m-1,\theta}}$
(Lemma~\ref{distTmTm-1smallboxes} and Proposition~\ref{distTmTm-1}).
This will be a key tool in the proof of
Propositions~\ref{seqTmProperties} for $m > 0$ and \ref{distTmTm-1}.
The next two lemmas and remark will be useful to automate
and simplify the proofs of these two results.

\begin{lemma}\label{normainterna}
\[
 \norm{f_{m,\theta} - f_{m-1,\theta}} =
 \begin{cases}
  \norm{f_{m,\theta}\evalat{\I_{\bt{\theta}, \theta}} - f_{m-1,\theta}\evalat{\I_{\bt{\theta}, \theta}}}
        & \text{when $\theta \in \wIBD\setminus\WIB$, and}\\[10pt]
  \norm{f_{m,\theta}\evalat{\IW{\theta}} - f_{m-1,\theta}\evalat{ \IW{\theta}}}
        & \text{when $\theta \in \WIB$,}
\end{cases}
\]
for every $m \ge 2$ and $\theta \in \wIBD.$
\end{lemma}


\begin{proof}
Set $i = \bt{\theta} \in \DS$, so that $\theta \in \wbasint{i}$.

When $\theta \in \wIBD\setminus\WIB = \IBD \cup \WB \setminus \WIB,$
by Definition~\ref{seqTmDefi} and Lemma~\ref{seqTmDefirevisited},
it is enough to show that
\[
 \abs{f_{m,\theta}(x) - f_{m-1,\theta}(x)} \le
 \abs{f_{m,\theta}(m_i(\theta)) - f_{m-1,\theta}(m_i(\theta))}
\]
for every $x \in [-2, m_i(\theta)]$, and
\[
 \abs{f_{m,\theta}(x) - f_{m-1,\theta}(x)} \le
 \abs{f_{m,\theta}(M_i(\theta)) - f_{m-1,\theta}(M_i(\theta))}
\]
for every $x \in [M_i(\theta),2]$.
We will prove the first statement. The second one follows similarly.

Definition~\ref{seqTmDefi} and Lemma~\ref{seqTmDefirevisited} give
\begin{align*}
f_{m,\theta}(x) - f_{m-1,\theta}(x) &=
     \frac{2 - g_{_{i,\theta}}\left(m_{i}(\theta)\right)\hfill}{2 - f_{m-1,\theta}\left(m_{i}(\theta)\right)}
                       (f_{m-1,\theta}(x) - 2) + 2 - f_{m-1,\theta}(x)\\
  &= \frac{2 - f_{m,\theta}\left(m_{i}(\theta)\right)\hfill}{2 - f_{m-1,\theta}\left(m_{i}(\theta)\right)}
                       (f_{m-1,\theta}(x) - 2) - (f_{m-1,\theta}(x) - 2)\\
  &=(f_{m-1,\theta}(x) - 2)\left(
            \frac{2 - f_{m,\theta}\left(m_{i}(\theta)\right)\hfill}{2 - f_{m-1,\theta}\left(m_{i}(\theta)\right)}
     - 1\right)\\
  &=(2 - f_{m-1,\theta}(x)) \frac{f_{m,\theta}\left(m_{i}(\theta)\right) - f_{m-1,\theta}\left(m_{i}(\theta)\right)}{2 - f_{m-1,\theta}\left(m_{i}(\theta)\right)}.
\end{align*}

By Lemma~\ref{seqTmPropsInBasicintervals}(b),
$2 \ge f_{m-1,\theta}(x) \ge f_{m-1,\theta}\left(m_{i}(\theta)\right)$
and $1 \ge f_{m-1,\theta}\left(m_{i}(\theta)\right).$ Hence,
\begin{align*}
\abs{f_{m,\theta}(x) - f_{m-1,\theta}(x)} &=
      (2 - f_{m-1,\theta}(x)) \frac{\abs{f_{m,\theta}\left(m_{i}(\theta)\right) - f_{m-1,\theta}\left(m_{i}(\theta)\right)}}{2 - f_{m-1,\theta}\left(m_{i}(\theta)\right)}\\
  &\le \abs{f_{m,\theta}\left(m_{i}(\theta)\right) - f_{m-1,\theta}\left(m_{i}(\theta)\right)}.
\end{align*}

Now assume that $\theta \in \WIB \subset \WB.$
By Definition~\ref{seqTmDefi} it is enough to show that
\[
 \abs{f_{m,\theta}(x) - f_{m-1,\theta}(x)} \le
 \abs{f_{m,\theta}(\lambda_m(\theta)) - f_{m-1,\theta}(\lambda_m(\theta))}
\]
for every $x \in [-2, \lambda_m(\theta)]$, and
\[
 \abs{f_{m,\theta}(x) - f_{m-1,\theta}(x)} \le
 \abs{f_{m,\theta}(\tau_m(\theta)) - f_{m-1,\theta}(\tau_m(\theta))}
\]
for every $x \in [\tau_m(\theta),2]$.
As before, we will prove the first statement. The second one follows similarly.
We have
\[
f_{m,\theta}(x) - f_{m-1,\theta}(x) =
    (2 - f_{m-1,\theta}(x)) \frac{f_{m,\theta}\left(\lambda_m(\theta)\right) - f_{m-1,\theta}\left(\lambda_m(\theta)\right)}{2 - f_{m-1,\theta}\left(\lambda_m(\theta)\right)}.
\]

By Lemma~\ref{seqTmPropsInBasicintervals}(b),
$2 \ge f_{m-1,\theta}(x) \ge f_{m-1,\theta}\left(\lambda_m(\theta)\right)$
and hence,
\[
\abs{f_{m,\theta}(x) - f_{m-1,\theta}(x)}
      \le \abs{f_{m,\theta}\left(m_{i}(\theta)\right) - f_{m-1,\theta}\left(m_{i}(\theta)\right)}
\]
provided that $2 - f_{m-1,\theta}\left(\lambda_m(\theta)\right) \ne 0.$
Assume by way of contradiction that we have
$f_{m-1,\theta}\left(\lambda_m(\theta)\right) = 2.$
Then, by Definition~\ref{curvesinthewings} and
Lemma~\ref{seqTmPropsInBasicintervals}(b),
$-1 \le \lambda_m(\theta)$ and
\[
 2 \ge f_{m-1,\theta}(-1) \ge f_{m-1,\theta}\left(\lambda_m(\theta)\right) = 2;
\]
which contradicts statement (b.iv) from the proof of
Lemma~\ref{seqTmPropsInBasicintervals}.
\end{proof}

Next we compute an upper bound for
$\norm{f_{m,\theta} - f_{m-1,\theta}}$
for every $\theta \in \wbasint{i}$ and
$i \in \DS$ such that $\diam(\wbasint{i})$ is small enough.

\begin{lemma}\label{distTmTm-1smallboxes}
Assume that $T_{m-1}$ is continuous for some $m \ge 2$
and let $\varepsilon$ be positive. Then, there exist
$\varrho_m(\varepsilon) \in \N$ such that
\[
 \norm{f_{m,\theta} - f_{m-1,\theta}} \le \varepsilon
\]
for every $\theta \in \wbasint{i}$ and
$i \in \DS$ (that is, $\wbasint{i} \subset \wIBD$)
such that $\ai \ge \varrho_m(\varepsilon).$
\end{lemma}

\begin{proof}
Since $T_{m-1}$ is uniformly continuous, there exists
$\delta_{m-1} = \delta_{m-1}(\varepsilon) > 0$ such that
$\dom(T_{m-1}(\theta,x), T_{m-1}(\nu, y)) < \varepsilon$
provided that $\dom((\theta,x), (\nu, y)) < \delta_{m-1}.$
We choose $\varrho_m = \varrho_m(\varepsilon) \in \N$ such that
\[
 3 \cdot 2^{-\varrho_m} < \min\{ \delta_{m-1}(\varepsilon/2), \varepsilon/2 \}.
\]

Assume that $i \in \DS$ verifies $\ai \ge \varrho_m(\varepsilon)$
and let $(\theta,x) \in \wbasband{i} = \wbasint{i} \times \I.$
When $\theta \in \wbasint{i} \setminus \WIB$
we can use Lemma~\ref{normainterna} with $\I_{i, \theta}$
to compute $\norm{f_{m,\theta} - f_{m-1,\theta}}.$
We have to show that
$\abs{f_{m,\theta}(x) - f_{m-1,\theta}(x)} < \varepsilon$
for every $x \in  \I_{i, \theta}.$

Let $\nu \in \Bd\left(\wbasint{i}\right) \subset \wEIBD$.
We have $(\theta,x),(\nu, m_i(\nu)) \in \wbasicbox{i}$ and,
by Lemmas~\ref{seqTmPropsInBasicintervals}(c) and \ref{Propertiesvarphi}(f),
\begin{align*}
\dom(T_m(\theta, x),T_m(\nu, m_i(\nu)))
   & = \dom(G_i(\theta, x),G_i(\nu, m_i(\nu)))\\
   & \le \diam\left(G_i\left(\wbasicbox{i}\right)\right),\text{ and}\\
\dom((\theta, x),(\nu, m_i(\nu))
   & \le \diam\left(\wbasicbox{i}\right)
     \le 2 \cdot 2^{-\ai} < 3 \cdot 2^{-\varrho_m} < \delta_{m-1}(\varepsilon/2).
\end{align*}
Thus,
\[ \dom(T_{m-1}(\theta, x),T_{m-1}(\nu, m_i(\nu)) < \varepsilon/2. \]
Consequently, by Lemma~\ref{voresiguals},
\begin{align*}
\abs{f_{m,\theta}(x) - f_{m-1,\theta}(x)}
   & = \dom(T_m(\theta, x),T_{m-1}(\theta, x))\\
   & \le \dom(T_m(\theta, x),T_{m-1}(\nu, m_i(\nu)))\ +\\
   &\hspace*{4em}\dom(T_{m-1}(\nu, m_i(\nu)), T_{m-1}(\theta, x))\\
   & <   \dom(T_m(\theta, x),T_{m}(\nu, m_i(\nu))) +  \varepsilon/2\\
   & <   \diam\left(G_i\left(\wbasicbox{i}\right)\right) +  \varepsilon/2.
\end{align*}

Now we look at the size of $G_i\left(\wbasicbox{i}\right).$
When $i < 0$, from Lemmas~\ref{gnegativa}(b) and \ref{Propertiesvarphi}(f),
we obtain
\begin{equation}\label{diamineg}
 \diam\left(G_i\left(\wbasicbox{i}\right)\right)
   \le \diam\left(\basicbox{i+1}\right) \le 2^{-(\ai-1)} < 2 \cdot 2^{-\ai}.
\end{equation}

When $i \ge 0,$ from Lemma~\ref{gpositiva}(b) we get
\[
 G_i\left(\wbasicbox{i}\right) = G_i\left(\basicbox{i}\right) \subset
   \basicbox{i+1} \cup \setfib{\A_{i+1}}{\left(\BSG{i+1}{i} \setminus \obasint{i+1}\right)}.
\]
Moreover, as in the proof of Lemma~\ref{Propertiesvarphi}(f) for $\ell < 0,$
the set
\[
\basicbox{i+1} \cup \setfib{\A_{i+1}}{\left(\BSG{i+1}{i} \setminus \obasint{i+1}\right)}
\]
is connected. So, by Lemma~\ref{Propertiesvarphi}(f),
\begin{align*}
\diam\left(G_i\left(\wbasicbox{i}\right)\right)
    &\le \diam\left(\basicbox{i+1} \cup \setfib{\A_{i+1}}{\left(\BSG{i+1}{i} \setminus \obasint{i+1}\right)}\right)\\
    &\le \diam\left(\basicbox{i+1}\right) + \diam\left(\setfib{\A_{i+1}}{\left(\BSG{i+1}{i} \setminus \obasint{i+1}\right)}\right)\\
    &\le 2^{-(i+1)} + \diam\left(\setfib{\A_{i+1}}{\left(\BSG{i+1}{i} \setminus \obasint{i+1}\right)}\right).
\end{align*}
As noticed earlier,
$\BSG{i+1}{i} \setminus \obasint{i+1}$ is disjoint from
\[ \obasint{i+1} \cup \wbasint{-(i+1)} \cup \Zstar_{i+1} \]
by Definition~\ref{PCgenerators}(\tsfR.2) and
Remark~\ref{PCgeneratorsExplicitConsequences}(\tsfR.2).
So, by Lemma~\ref{propiedadesA}(c),
Definition~\ref{PCgenerators} and Lemma~\ref{Propertiesvarphi}(a),
\begin{align*}
\setfibpt{\A_{i+1}}{\nu}
    &= \{(\nu, \gams{i+1}(\nu)\} = \{(\nu, \gams{i}(\nu)\}\\
    &\in \{\nu\} \times \left[\gams{i}(\sstar{i+1})-2^{-n_{i}}, \gams{i}(\sstar{i+1}) + 2^{-n_{i}}\right].
\end{align*}
for every $\nu \in \BSG{i+1}{i} \setminus \obasint{i+1}.$
On the other hand,
$
\gams{i}(\sstar{i+1}) \in \I_{i+1,\sstar{i+1}}
$
by Lemma~\ref{Propertiesvarphi}(c).
Hence, by Remark~\ref{propiedadesR}(2),
Definition~\ref{PCgenerators}(\tsfR.1) and
Remark~\ref{PCgeneratorsExplicitConsequences}(\tsfR.1),
\begin{align*}
\diam&\left(\setfib{\A_{i+1}}{\left(\BSG{i+1}{i} \setminus \obasint{i+1}\right)}\right)\\
    &\le \max\left\{
               \diam\left(\BSG{i+1}{i} \setminus \obasint{i+1}\right),
               2\cdot (2^{-n_{i}}+2^{-n_{i+1}})
         \right\}\\
    &\le 2\cdot \max\left\{ \alpha_i, 2^{-n_{i}}+2^{-n_{i+1}}
         \right\} = 2\cdot (2^{-n_{i}}+2^{-n_{i+1}})\\
    &< 4\cdot 2^{-n_{i}} \le 2 \cdot 2^{-i}.
\end{align*}
Summarizing, when $i \ge 0$,
\[
 \diam\left(G_i\left(\wbasicbox{i}\right)\right)
    \le 2^{-(i+1)} + 2 \cdot 2^{-i} < 3 \cdot 2^{-i}
\]
and, from \eqref{diamineg},
\[
   \diam\left(G_i\left(\wbasicbox{i}\right)\right)
   <   3 \cdot 2^{-\ai} \le  3 \cdot 2^{-\varrho_m} < \varepsilon/2
\]
for every $i \in \Z^+$.
Thus, for every $x \in  \I_{i, \theta}$,
\[
\abs{f_{m,\theta}(x) - f_{m-1,\theta}(x)}
   < \diam\left(G_i\left(\wbasicbox{i}\right)\right) +  \varepsilon/2 < \varepsilon.
\]

Now assume that $\theta \in \wbasint{i} \cap \WIB.$
We can use Lemma~\ref{normainterna} with $\IW{\theta}$
to compute $\norm{f_{m,\theta} - f_{m-1,\theta}}.$
We have to show that
$\abs{f_{m,\theta}(x) - f_{m-1,\theta}(x)} < \varepsilon$
for every $x \in  \IW{\theta}.$
Since $\theta \in \WIB,$
by Definition~\ref{curvesinthewings} and Lemma~\ref{VerticalIntervalsIntheWings}(b),
$i < 0$, $\theta \in \WB$ and
\[
 \I_{i, \theta} = \left\{\gams{\ai}(\theta)\right\}
 \subset \IW{\theta} = \I_{\ell, \theta} \ni x
\]
with $\ell = \bt[\led{\theta}]{\theta} \in \WDS.$
In this case we will consider the points
$(\theta,x) \in \basicbox{\ell}$ and
$(\nu, m_i(\nu)),(\theta,\gams{\ai}(\theta))  \in \wbasicbox{i}$ with
$\nu \in \Bd\left(\wbasint{i}\right) \subset \wEIBD$.
By Lemma~\ref{DepthintheWings}(b), Remark~\ref{propiedadesR}(2)
and Lemma~\ref{Propertiesvarphi}(f), $\ai < \all$ and
\begin{align*}
\dom((\theta, x),(\nu, m_i(\nu))
   & \le \dom((\theta, x),(\theta,\gams{\ai}(\theta)) +
         \dom((\theta,\gams{\ai}(\theta)),(\nu, m_i(\nu))\\
   & \le \abs{x - \gams{\ai}(\theta)} +
         \diam\left(\wbasicbox{i}\right)\\
   & \le \diam\left(\basicbox{\ell}\right) +
         \diam\left(\wbasicbox{i}\right)\\
   & \le 2^{-\all} + 2 \cdot 2^{-\ai} < 3 \cdot 2^{-\ai}
     \le 3 \cdot 2^{-\varrho_m} < \delta_{m-1}(\varepsilon/2).
\end{align*}
Thus,
\[ \dom(T_{m-1}(\theta, x),T_{m-1}(\nu, m_i(\nu)) < \varepsilon/2. \]
On the other hand,
by Lemma~\ref{seqTmPropsInBasicintervals}(c), Definition~\ref{seqTmDefi}
and \eqref{diamineg},
\begin{align*}
\dom(T_m(\theta&, x),T_m(\nu, m_i(\nu)))\\
   & \le \dom(T_m(\theta, x),T_m(\theta,\gams{\ai}(\theta))) +
         \dom(T_m(\theta,\gams{\ai}(\theta)),T_m(\nu, m_i(\nu)))\\
   & \le \abs{f_{m, \theta}(x) - f_{m, \theta}(\gams{\ai}(\theta))} +
         \dom(G_i(\theta,\gams{\ai}(\theta)),G_i(\nu, m_i(\nu)))\\
   & = \dom(G_i(\theta, x),G_i(\nu, m_i(\nu)))
     \le \diam\left(G_i\left(\wbasicbox{i}\right)\right) < 2 \cdot 2^{-\ai}\\
   & \le  3 \cdot 2^{-\varrho_m} < \varepsilon/2.
\end{align*}
So, in a similar way as before, Lemma~\ref{voresiguals} gives
\begin{align*}
\abs{f_{m,\theta}(x) - f_{m-1,\theta}(x)}
   & = \dom(T_m(\theta, x),T_{m-1}(\theta, x))\\
   & \le \dom(T_m(\theta, x),T_{m-1}(\nu, m_i(\nu)))\ +\\
   &\hspace*{4em}\dom(T_{m-1}(\nu, m_i(\nu)), T_{m-1}(\theta, x))\\
   & < \varepsilon.
\end{align*}
\end{proof}

\begin{proof}[Proof of Proposition~\ref{seqTmProperties} for $m > 0$]
\inidemopart{a}
We start by proving by induction on $m$ that $T_m$ is continuous for every $m \in \Z^+.$

By Proposition~\ref{T0mapProperties}(a) for $m=0$, $T_0$ is continuous.
So, we may assume that $T_{m-1}$ is continuous for some $m\in \N$
and prove that $T_m$ is continuous.

Let $\varepsilon > 0$ be fixed but arbitrary, and
let $(\theta,x), (\nu, y) \in \Omega.$
We have to show that there exists $\delta(\varepsilon) > 0$
such that
\[ \dom(T_m(\theta,x), T_m(\nu, y)) < \varepsilon
\andq[when]
\dom((\theta,x), (\nu, y)) < \delta.
\]

We start by defining $\delta(\varepsilon)$. To this end we need to
introduce some more notation and establish some facts
about the maps $T_m$ and $T_{m-1}.$

Since $T_{m-1}$ is uniformly continuous, we know that
\begin{equation}\label{Tmm1UC}
\text{\parbox{0.9\textwidth}{
there exists $\delta_{m-1} = \delta_{m-1}(\varepsilon) > 0$ such that
$\dom(T_{m-1}(\theta,x), T_{m-1}(\nu, y)) < \varepsilon$
provided that $\dom((\theta,x), (\nu, y)) < \delta_{m-1}.$
}}\end{equation}

On the other hand, Lemma~\ref{seqTmPropsInBasicintervals}(a)
tells us that
$T_m\evalat{\wbasband{i}}$ is uniformly continuous
for every $i \in \DS.$
So, for every $i \in \DS,$
\begin{equation}\label{TmUCalesbandes}
\text{\parbox{0.9\textwidth}{
there exists $\delta_{m,i} = \delta_{m,i}(\varepsilon) > 0$
such that
$\dom(T_m(\theta,x), T_m(\nu, y)) < \varepsilon$
for every $(\theta,x), (\nu, y) \in \wbasband{i} \subset \wIVD$
verifying $\dom((\theta,x), (\nu, y)) < \delta_{m,i}(\varepsilon).$
}}\end{equation}

Then, by using the numbers
$\delta_{m-1}(\varepsilon/7)$ given by \eqref{Tmm1UC},
$\delta_{m,i}(\varepsilon/7)$ given by \eqref{TmUCalesbandes} and
$\varrho_m(\varepsilon/7)$ given by Lemma~\ref{distTmTm-1smallboxes},
we set
\[
\delta = \delta(\varepsilon) := \min\left\{
    \delta_{m-1}(\varepsilon/7),
    \min \set{\delta_{m,i}(\varepsilon/7)}{i \in \DS \cap Z_{\varrho_m(\varepsilon/7)}}
\right\}.
\]
Clearly, $\delta > 0$ because the set $\DS \cap Z_{\varrho_m(\varepsilon/7)}$ is finite.

Now we will show that if $\dom((\theta,x), (\nu, y)) < \delta,$
then $\dom(T_m(\theta,x), T_m(\nu, y)) < \varepsilon.$

Assume first that
$(\theta,x),(\nu,y) \in \wbasband{\ell}$
for some $\ell \in \DS \cap Z_{\varrho_m(\varepsilon/7)}.$
We have
\[
 \dom((\theta,x), (\nu, y)) < \delta
    \le \min \set{\delta_{m,i}(\varepsilon/7)}{i \in \DS \cap Z_{\varrho_m(\varepsilon/7)}}
    \le \delta_{m,\ell}(\varepsilon/7).
\]
Hence, by \eqref{TmUCalesbandes},
\[
 \dom(T_m(\theta,x), T_m(\nu, y)) < \varepsilon/7 < \varepsilon.
\]

Next we assume that
$(\theta,x),(\nu,y) \in \wbasband{\ell}$
for some $\ell \in \DS$ such that $\all > \varrho_m(\varepsilon/7)$
(in particular, $\theta,\nu \in \wbasint{\ell}$).
In this situation we have
\[
 \dom((\theta,x), (\nu, y)) < \delta \le \delta_{m-1}(\varepsilon/7)
\]
and, by \eqref{Tmm1UC} and Lemma~\ref{distTmTm-1smallboxes},
\begin{align*}
\dom(T_m(\theta,x), T_m(\nu, y))
   &\le \dom(T_m(\theta,x), T_{m-1}(\theta,x)) +
         \dom(T_{m-1}(\theta,x), T_{m-1}(\nu, y))\ +\\
   &    \hspace*{4em}\dom(T_{m-1}(\nu, y), T_m(\nu, y)) \\
   &=   \abs{f_{m,\theta}(x) - f_{m-1,\theta}(x)} +
        \dom(T_{m-1}(\theta,x), T_{m-1}(\nu, y))\ +\\
   &    \hspace*{4em}\abs{f_{m,\nu}(y) - f_{m-1,\nu}(y)}\\
   &\le \norm{f_{m,\theta} - f_{m-1,\theta}} +
        \dom(T_{m-1}(\theta,x), T_{m-1}(\nu, y))\ +\\
   &    \hspace*{4em}\norm{f_{m,\nu} - f_{m-1,\nu}}\\
   &<   \tfrac{3}{7} \varepsilon < \varepsilon .
\end{align*}

In summary, we have proved that
\[ \dom(T_m(\theta,x), T_m(\nu, y)) < \tfrac{3}{7} \varepsilon \]
when $\dom((\theta,x), (\nu, y)) < \delta$
and $(\theta,x),(\nu,y) \in \wbasband{\ell}$ for some $\ell \in \DS.$

Next we assume that $(\theta,x),(\nu,y) \in \wIVD$ but
$(\theta,x),(\nu,y) \notin \wbasband{\ell}$ for every $\ell \in \DS.$
By Lemma~\ref{Dsets}(a,b), there exist
$i = \bt{\theta},k=\bt{\nu} \in \DS,$ $i \ne k,$ such that
$\theta \in \wbasint{i},$
$(\theta,x) \in \wbasband{i},$
$\nu \in \wbasint{k}$ and
$(\nu,y) \in \wbasband{k}.$
Then, there exist
\[
  \widetilde{\theta} \in A \cap \Bd\left(\wbasint{i}\right) \subset \wEIBD
  \andq
  \widetilde{\nu} \in A \cap \Bd\left(\wbasint{k}\right) \subset \wEIBD,
\]
where $A$ denotes the closed arc of $\SI$ such that
\[ \diam(A) = \dSI(\theta,\nu) \andq \Bd(A) = \{\theta, \nu\}. \]
Clearly we have,
$(\theta,x), \bigl(\widetilde{\theta}, x\bigr) \in \wbasband{i},$
$(\nu, y), \bigl(\widetilde{\nu}, y\bigr) \in \wbasband{k}$ and,
by the previous case,
\begin{align*}
\dom\left((\theta,x), \bigl(\widetilde{\theta}, x\bigr)\right)
   &= \dSI\bigl(\theta,\widetilde{\theta}\bigr)
      \le \dSI(\theta,\nu) \le \dom((\theta,x), (\nu, y)) < \delta,\\
   &  \dom\left(T_m(\theta,x), T_m\bigl(\widetilde{\theta}, x\bigr)\right) < \tfrac{3}{7} \varepsilon\\
\dom\left((\nu,y), \bigl(\widetilde{\nu}, y\bigr)\right)
   &= \dSI\bigl(\nu,\widetilde{\nu}\bigr)
      \le \dSI(\theta,\nu) \le \dom((\theta,x), (\nu, y)) < \delta,\text{ and}\\
   &  \dom\left(T_m(\nu, y), T_m\bigl(\widetilde{\nu}, y\bigr)\right) < \tfrac{3}{7} \varepsilon.
\end{align*}
On the other hand,
$
 \bigl(\widetilde{\theta}, x\bigr),
 \bigl(\widetilde{\nu}, y\bigr) \in \setsilift{\wEIBD}
 \subset \wIVD  \subset \wIVD[m-1]
$
and, by Lemma~\ref{voresiguals} and \eqref{Tmm1UC},
\begin{align*}
\dom\left(\bigl(\widetilde{\theta},x\bigr), \bigl(\widetilde{\nu}, y\bigr)\right)
   &=   \max\left\{\dSI\bigl(\widetilde{\theta},\widetilde{\nu}\bigr), \abs{x-y}\right\}
    \le \max\left\{\dSI(\theta,\nu), \abs{x-y}\right\}\\
   &=   \dom((\theta,x), (\nu, y)) < \delta \le \delta_{m,i}(\varepsilon/7),\text{ and}\\
\dom(T_m(\theta,x), T_m(\nu, y))
   &\le \dom\left(T_m(\theta,x), T_m\bigl(\widetilde{\theta}, x\bigr)\right) +
        \dom\left(T_m\bigl(\widetilde{\theta}, x\bigr), T_m\bigl(\widetilde{\nu}, y\bigr)\right)\ +\\
   &    \hspace*{4em}\dom\left(T_m\bigl(\widetilde{\nu}, y\bigr), T_m(\nu, y)\right)\\
   &< \tfrac{3}{7} \varepsilon +
      \dom\left(T_{m-1}\bigl(\widetilde{\theta}, x\bigr), T_{m-1}\bigl(\widetilde{\nu}, y\bigr)\right) +
      \tfrac{3}{7} \varepsilon = \varepsilon.
\end{align*}

If $(\theta,x),(\nu,y) \notin \wIVD$ then,
by Definition~\ref{seqTmDefi} and \eqref{Tmm1UC} ,
\[
 \dom(T_m(\theta,x), T_m(\nu, y)) = \dom(T_{m-1}(\theta,x), T_{m-1}(\nu, y))
     < \varepsilon/7 < \varepsilon
\]
because
$
 \dom((\theta,x), (\nu, y)) < \delta \le \delta_{m-1}(\varepsilon/7).
$

Lastly, assume that $(\nu,y) \notin \wIVD$ but
$(\theta,x) \in \wbasband{i} \subset \wIVD,$
for some $i\in\DS$ (that is, $\theta \in \wbasint{i}$).
In this situation, as before, there exists
$
  \widetilde{\theta} \in \Bd\left(\wbasint{i}\right) \subset \wEIBD
$
such that, by Lemma~\ref{voresiguals} and Definition~\ref{seqTmDefi}
($\bigl(\widetilde{\theta}, x\bigr) \in  \setsilift{\wEIBD} \subset \wIVD  \subset \wIVD[m-1]$), and \eqref{Tmm1UC},
\begin{align*}
\dom\left((\theta,x), \bigl(\widetilde{\theta}, x\bigr)\right) & < \delta,\\
\dom\left(\bigl(\widetilde{\theta},x\bigr), (\nu, y\bigr)\right) & < \delta  \le \delta_{m-1}(\varepsilon/7),\\
\dom\left(T_m(\theta,x), T_m\bigl(\widetilde{\theta}, x\bigr)\right) & < \tfrac{3}{7} \varepsilon,\text{ and}\\
\dom(T_m(\theta,x), T_m(\nu, y))
   &\le \dom\left(T_m(\theta,x), T_m\bigl(\widetilde{\theta}, x\bigr)\right) +
        \dom\left(T_m\bigl(\widetilde{\theta}, x\bigr), T_m(\nu, y)\right)\\
   &< \tfrac{3}{7} \varepsilon +
      \dom\left(T_{m-1}\bigl(\widetilde{\theta}, x\bigr), T_{m-1}(\nu, y)\right) < \varepsilon.
\end{align*}

This ends the proof of the continuity of $T_m$
and, hence, of (a).

\inidemopart{b}
When $\theta \in \wIBD$ the statement follows from Lemma~\ref{seqTmPropsInBasicintervals}(b).
When $\theta \in \SI\setminus\wIBD,$
it follows from the part already proven and
the continuity of $T_m.$

\inidemopart{c}
The first two statements follow from
Lemma~\ref{seqTmPropsInBasicintervals}(c) and
statement (a). On the other hand,
as in the proof of Proposition~\ref{T0mapProperties}(c) for $m=0$,
Lemma~\ref{denso}(b) implies that
$\istar \in \wIBD$ but $\istar \notin \wIBD[k]$
for every $k > m.$
Then, we get $f_{k,\istar} = f_{m,\istar}$
from Definition~\ref{seqTmDefi}.
\end{proof}

\section{Proof of Proposition~\ref{distTmTm-1}}\label{proofofdistTmTm-1}

This section is devoted to prove
Proposition~\ref{distTmTm-1}.
It is the third technical counterpart of Section~\ref{skew-product}.
In contrast to Lemma~\ref{distTmTm-1smallboxes}
the bound given by Proposition~\ref{distTmTm-1}.
is valid for every $\theta \in \wIBD$.

Before starting the proof of this proposition we will state and
prove a number of very simple lemmas that will help in
automating the proof of Proposition~\ref{distTmTm-1}.

\begin{lemma}\label{pointnormbound}
Assume that $\wbasint{i} \subset \wbasint{k}$
for some $i \in \DS,\ k \in \DS[m-1]$ and $m \ge 2,$
and assume that either
\[
   i < 0 \text{ and } \theta \in \wbasint{i} \setminus \istarset
      \text{ or }
   i \ge 0 \text{ and } \theta \in \basint{i} \setminus \OBG{i}{i+1}.
\]
Then,
\[
\abs{\gams{\aii}\left(R_\omega(\theta)\right) - \gams{\akk}\left(R_\omega(\theta)\right)}
    \le 2^{-\ak}.
\]
\end{lemma}

\begin{proof}
The lemma holds trivially when $\akk = \aii.$
Thus, we may assume that $\akk \ne \aii.$
Then by Lemma~\ref{QuePassaALesAles},
$\ak < \ai,$ $\akk < \aii$ and
\[
 \gams{\akk}\left(R_\omega(\theta)\right) =
 \gams{\akk+1}\left(R_\omega(\theta)\right) = \dots =
 \gams{\aii-1}\left(R_\omega(\theta)\right).
\]
By assumption we have
\[
 \theta \in \begin{cases}
    \basint{i} \setminus \OBG{i}{i+1} &\text{when $i \ge 0$, and}\\
    \wbasint{i} \setminus \istarset =
         \BSG{i}{\aii}\setminus \istarset
                 &\text{when $i < 0,$}
\end{cases}
\]
and, hence,
\[
R_\omega(\theta) \in \begin{cases}
    \BSG{i+1}{i} \setminus \obasint{i+1}     &\text{when $i \ge 0$, and}\\
    \basintabs{i+1} \setminus \sstarset{i+1} &\text{when $i < 0$}.
\end{cases}
\]
Thus, in the case $i \ge 0$ we have
\[ R_\omega(\theta) \notin \obasint{i+1} \cup \wbasint{-(i+1)} \cup \Zstar_{i+1} \]
by Definition~\ref{PCgenerators}(\tsfR.2) and
Remark~\ref{PCgeneratorsExplicitConsequences}(\tsfR.2).
So, by Definition~\ref{PCgenerators},
\[
\gams{i+1}\left(R_\omega(\theta)\right)
   = \gams{i}\left(R_\omega(\theta)\right)
   = \gams{\akk}\left(R_\omega(\theta)\right).
\]
This ends the proof of the lemma in this case.

Assume now that $i < 0.$
By Lemma~\ref{Propertiesvarphi}(c,d,f) and
Definition~\ref{PCgenerators}(\tsfR.2) and
Remark~\ref{PCgeneratorsExplicitConsequences}(\tsfR.2),
\begin{align*}
\abs{\gams{\aii}\left(R_\omega(\theta)\right) - \gams{\akk}\left(R_\omega(\theta)\right)}
  &= \abs{\gams{\aii}\left(R_\omega(\theta)\right) - \gams{\aii-1}\left(R_\omega(\theta)\right)}\\
  &\le \diam\left(\basicbox{i+1}\right) \le 2^{-\aii} \le 2^{-\ak}
\end{align*}
(observe that $\aii > \akk \ge \ak - 1$).
\end{proof}

\begin{lemma}\label{intervalsnormbound}
Let $s, t \in \Z,$ $s \ne t$ be such that
$\theta \in \wobasint{s} \setminus \obasintabs{s},$
and either $t < 0$ and $\theta \in \obasintabs{t}$
or $t \ge 0$ and $\theta \in \OBG{t}{t+1}.$
Then, the following statements hold:
\begin{enumerate}[(a)]
\item $R_\omega(\theta) \in \obasintabs{s+1} \cap \obasintabs{t+1}.$
\item Let $u,v \in \Z$ be such that $\{u,v\} = \{ s, t\}$
      and $\abs{u + 1} \le \abs{v + 1}$. \\
      Then, $\I_{v+1, R_\omega(\theta)} \subset \I_{u+1, R_\omega(\theta)}.$
\item \[
       \abs{x - y} \le 2\cdot 2^{-\abs{u}}
      \]
      for every $x \in \I_{t+1, R_\omega(\theta)}$ and $y \in \I_{s+1, R_\omega(\theta)}.$
\end{enumerate}
\end{lemma}

\begin{proof}
By assumption we have
\[
 \theta \in \begin{cases}
    \OBG{t}{t+1} &\text{when $t \ge 0$, and}\\
    \obasintabs{t} \subset \wobasint{t} = \OBG{t}{\abs{t+1}}
                 &\text{when $t < 0.$}
\end{cases}
\]
Hence,
$
R_\omega(\theta) \in \obasintabs{t+1}.
$
Moreover, as in the proof of Lemma~\ref{pointnormbound},
$s < 0$ and
$
R_\omega(\theta) \in \obasintabs{s+1}.
$
This proves (a).

Now we prove (b). From (a) we have
\begin{align*}
R_\omega(\theta) &\in \obasintabs{u+1} \cap \obasintabs{v+1}\\
                 &\subset \obasintabs{u+1} \cap \wbasint{v+1}.
\end{align*}
Moreover,
$s \ne t$ implies $u+1 \ne v+1$ and we have $\abs{u+1} \le \abs{v+1}$
by assumption. Consequently, by Lemma~\ref{Propertiesvarphi}(g,d) and
Definition~\ref{PCgenerators}(\tsfR.2) and
Remark~\ref{PCgeneratorsExplicitConsequences}(\tsfR.2),
$\abs{u+1} < \abs{v+1}$ and
\[
  \basicbox{v+1} \subset \Int\left(\basicbox{u+1} \setminus \setsilift{\sstarset{u+1}}\right)
\]
which implies (b).

Thus, $x,y \in \I_{u+1, R_\omega(\theta)}$ and, by Lemma~\ref{Propertiesvarphi}(f),
\[
  \abs{x - y} \le \diam\left(\basicbox{u+1}\right)
              \le 2^{-\abs{u+1}} \le 2^{-(\abs{u}-1)} = 2\cdot 2^{-\abs{u}}.
\]
\end{proof}

Now we are ready to start the proof of Proposition~\ref{distTmTm-1}.

\begin{proof}[Proof of Proposition~\ref{distTmTm-1}]
We start by showing that
$\{T_m\}_{k = 0}^{\infty}$ is a Cauchy sequence,
assuming that the bound \eqref{fitanorma}
holds for every $m \ge 2$ and $\theta \in \SI.$

We start by estimating $\dinf(T_m, T_{m+1})$ for every $m \in \N.$
From \eqref{fitanorma} and the definition of $\mu_m$
\[
 \dinf(T_m, T_{m+1})
   = \sup_{\theta \in \SI}  \norm{f_{m,\theta} - f_{m+1,\theta}}
   \le 2\cdot \sup_{\theta \in \SI} 2^{-\abs{\bt{\theta}}}
   \le 2\cdot 2^{-\mu_m}.
\]

By Lemma~\ref{denso}(a) $\{\mu_m\}_{m=0}^\infty$
is strictly increasing (and $\lim_{m\to\infty} \mu_m = \infty$).
Therefore, for every $\varepsilon > 0,$
there exists $N \ge 2,$ such that
$4\cdot 2^{-\mu_m} < \varepsilon$ for every $m \ge N.$
Hence,
\begin{align*}
\dinf(T_m, T_{m+i})
   &\le \sum_{\ell=m}^{m+i-1} \dinf(T_{\ell}, T_{\ell + 1})
    \le  2\cdot \sum_{\ell=m}^{m+i-1} 2^{-\mu_{\ell}}\\
   &\le  2\cdot 2^{-\mu_m} \sum_{\ell=0}^{\infty} 2^{-\ell}
    = 4\cdot 2^{-\mu_m} \le 4\cdot 2^{-\mu_N} < \varepsilon
\end{align*}
for every $m \ge N$ and $i \in \N.$
So, $\{T_m\}_{k = 0}^{\infty}$ is a Cauchy sequence.

Now we prove \eqref{fitanorma}. That is,
\[
 \norm{f_{m,\theta} - f_{m-1,\theta}} \le 2 \cdot 2^{-\abs{\bt[m-1]{\theta}}}
\]
for every $m \ge 2$ and $\theta \in \SI.$

From Definition~\ref{seqTmDefi} and Lemma~\ref{voresiguals}
we know that $f_{m,\theta} = f_{m-1,\theta}$ for every
$\theta \in \left(\SI \setminus \wIBD\right) \cup \wEIBD.$
Then, \eqref{fitanorma} holds in this case.

In the rest of the involved proof we assume that
$\theta \in \wIBD \setminus \wEIBD.$
Thus, by Lemmas~\ref{Dsets}(a,b), \ref{Propertiesvarphi}(g) and \ref{QuePassaALesAles},
\begin{quotation}\itshape
$
  \theta \in \wobasint{i} \subset
  \wobasint{k} \setminus \left( \Bd\left(\basintabs{k}\right) \cup \kstarset \right)
$
where\newline
$i = \bt{\theta} \in \DS,$ $k = \bt[m-1]{\theta} \in \DS[m-1],$\newline
$\ak < \ai,$ and $\akk \le \aii.$
\end{quotation}
Moreover, $\wbasband{i} \subset \wbasband{k} \subset \wIVD[m-1].$
Consequently, by Lemma~\ref{seqTmPropsInBasicintervals}(a,b),
the maps $f_{m, \theta}$ and $f_{m-1, \theta}$ are well defined, continuous,
piecewise affine and non-increasing, and
$f_{m,\theta}(2) = f_{m-1,\theta}(2) = -2$ and
$f_{m,\theta}(-2) = f_{m-1,\theta}(-2) = 2$
(see Figures~\ref{fig:Cas1.3}, \ref{fig:Cas2} and \ref{fig:Cas3.1}
 for some examples in generic cases).

We split the proof into three cases according to whether $\theta$
belongs to
\[
  \wobasint{i} \setminus \obasintabs{i},\
  \obasintabs{i} \subset \wobasint{k} \setminus \basintabs{k}
  \text{ or }
  \obasintabs{i} \subset \obasintabs{k}.
\]

\begin{autocase}{1}
$\theta \in \wobasint{i} \setminus \obasintabs{i}.$
\end{autocase}
We have $i < 0$ because $\wobasint{i} = \obasint{i}$ for $i \ge 0.$
Moreover, by Definition~\ref{curvesinthewings}, $\theta \in \WB.$

To deal with this case we consider three subcases.

\begin{autocase}[Subcase]{1.1}
$\theta \in \left(\wobasint{i} \setminus \obasintabs{i}\right) \setminus \WIB.$
\end{autocase}
By Lemmas~\ref{QuePassaALesAles}, \ref{fmfm-1alesales},
\ref{normainterna} and \ref{pointnormbound},
\begin{align*}
\I_{i, \theta} & = \{m_i(\theta)\} = \{\gams{\ai}(\theta)\} = \{\gams{\ak}(\theta)\} \subset \I_{k, \theta},\\
f_{m,\theta}\left(m_i(\theta)\right) & = \gams{\aii}\left(R_\omega(\theta)\right),\\
f_{m-1,\theta}\left(m_i(\theta)\right) & = \gams{\akk}\left(R_\omega(\theta)\right), \text{ and}\\
\norm{f_{m,\theta} - f_{m-1,\theta}} &=
 \norm{f_{m,\theta}\evalat{\I_{i, \theta}} - f_{m-1,\theta}\evalat{\I_{i, \theta}}} =
 \abs{f_{m,\theta}\left(m_{i}(\theta)\right) - f_{m-1,\theta}\left(m_{i}(\theta)\right)}\\
&= \abs{\gams{\aii}\left(R_\omega(\theta)\right) - \gams{\akk}\left(R_\omega(\theta)\right)}
   \le 2^{-\abs{\bt[m-1]{\theta}}}.
\end{align*}

\begin{autocase}[Subcase]{1.2}
$\theta \in \left(\wobasint{i} \setminus \obasintabs{i}\right) \cap \WIB$
and $\wobasint{i} \subset \wobasint{k} \setminus \basintabs{k}.$
\end{autocase}
In this subcase, by Definition~\ref{curvesinthewings} we have
\[
\theta \in \wobasint{k} \setminus \basintabs{k} \subset \WB[m-1]
\]
(recall that $i < 0$).
Then, by Lemmas~\ref{QuePassaALesAles} and \ref{VerticalIntervalsIntheWings}(b,c),
Definition~\ref{seqTmDefi} and Lemmas~\ref{normainterna} and \ref{pointnormbound},
\begin{align*}
\I_{i, \theta} & = \{\gams{\ai}(\theta)\} = \{\gams{\ak}(\theta)\} \subset \IW{\theta} = \IW[m-1]{\theta},\\
f_{m,\theta}(x) & = \gams{\aii}\left(R_\omega(\theta)\right) \text{ for every $x \in \IW{\theta}$,}\\
f_{m-1,\theta}(x) & = \gams{\akk}\left(R_\omega(\theta)\right) \text{ for every $x \in \IW[m-1]{\theta}$, and}\\
\norm{f_{m,\theta} - f_{m-1,\theta}}
   &= \norm{f_{m,\theta}\evalat{\IW{\theta}} - f_{m-1,\theta}\evalat{\IW{\theta}}} \\
   &= \abs{\gams{\aii}\left(R_\omega(\theta)\right) - \gams{\akk}\left(R_\omega(\theta)\right)}
    \le 2^{-\abs{\bt[m-1]{\theta}}}.
\end{align*}

Observe that since $\wobasint{i}$ is connected and
\[
  \wobasint{i} \subset
  \wobasint{k} \setminus \left( \Bd\left(\basintabs{k}\right) \cup \kstarset \right),
\]
$\wobasint{i} \not\subset \wobasint{k} \setminus \basintabs{k}$
implies $\wobasint{i} \subset \obasintabs{k} \setminus \kstarset.$

\begin{autocase}[Subcase]{1.3}
$\theta \in \left(\wobasint{i} \setminus \obasintabs{i}\right) \cap \WIB$
and $\wobasint{i} \subset \obasintabs{k} \setminus \kstarset$
\upshape  (see Figure~\ref{fig:Cas1.3}  for a symbolic representation of this case).
\end{autocase}
\begin{figure}
\begin{tikzpicture}[scale=2]
\draw (-2,-2) rectangle (2,2);
\foreach \c in {-2, 2} { \node[below] at (\c,-2) {$\c$}; \node[left] at (-2,\c) {$\c$}; }

\draw[dashed, color=blue] (-1,-2.3) -- (-1,2); \draw[dashed, color=blue] (1,-2.3) -- (1,2);
\draw[decorate, very thick, decoration={brace,amplitude=5pt, mirror, raise=2pt}, color=blue] (-1,-2.3) -- (1,-2.3);
\node[below, color=blue] at (0,-2.4) {\scriptsize$\I_{k,\theta}$};
\node[below, color=blue] at (-1,-2.4) {\scriptsize$m_k(\theta)$};
\node[below, color=blue] at (1,-2.4) {\scriptsize$M_k(\theta)$};

\draw[dashed, color=blue] (-2, 0.3) -- (2, 0.3); \draw[dashed, color=blue] (-2, 1) -- (2, 1);
\draw[decorate, very thick, decoration={brace,amplitude=5pt,raise=2pt}, color=blue]  (-2, 0.3) -- (-2, 1);
\node[left, color=blue] at (-2.1,0.65) {\scriptsize$\I_{k+1,R_\omega(\theta)}$};
\node[left, color=blue] at (-2.05, 0.3) {\scriptsize$m_{k+1}\bigl(R_\omega(\theta)\bigr)$};
\node[left, color=blue] at (-2.05, 1) {\scriptsize$M_{k+1}\bigl(R_\omega(\theta)\bigr)$};

\draw[very thick, color=blue] (-2,2) -- (-1.6, 1.85) -- (-1.1, 1.4) --(-1, 1)
                 -- (1, 0.3) -- (1.2, -1) -- (1.4, -1.1) -- (1.7, -1.8) -- (2,-2);
\node[color=blue, above right] at (1.15,-1) {$f_{m-1,\theta}$};
\node[color=blue] at (-0.5,0.95) {$g_{_{k, \theta}}$};

\draw[dashed, color=red] (-0.8,-2) -- (-0.8,2); \draw[dashed, color=red] (-0.1,-2) -- (-0.1,2);
\draw[decorate, very thick, decoration={brace,amplitude=5pt, mirror, raise=2pt}, color=red] (-0.8,-2) -- (-0.1,-2);
\node[below, color=red] at (-0.45,-2.1) {\scriptsize $\IW{\theta}$};
\node[below left, color=red] at (-0.75,-2.01) {\scriptsize$\lambda_m(\theta)$};
\node[below right, color=red] at (-0.15,-2.01) {\scriptsize$\tau_m(\theta)$};

\draw[dashed, color=red] (-2, 0.4) -- (2, 0.4);
\node[right, color=red] at (2, 0.4) {\scriptsize$\gams{\aii}\left(R_\omega(\theta)\right)$};

\draw[very thick, color=red] (-2,2) -- (-1.6, 1.8) -- (-1.1, 1.2266) --(-1, 0.711)
                 -- (-0.8, 0.4) -- (-0.1, 0.4) -- (1, -0.04848)
                 -- (1.2, -1.1515) -- (1.4, -1.2363) -- (1.7, -1.8303) -- (2,-2);
\node[color=red, left] at (-1,0.711) {$f_{m,\theta}$};
\end{tikzpicture}
\caption{A symbolic representation of the maps
$f_{m,\theta}$ and $f_{m-1,\theta}$ in
Subcase~1.3 of Proposition~\ref{distTmTm-1}
($\theta \in \left(\protect\wobasint{i} \setminus \protect\obasintabs{i} \right) \cap \WIB$
and $\protect\wobasint{i} \subset \protect\obasintabs{k} \setminus \protect\kstarset$).
The map \textcolor{blue}{$f_{m-1,\theta}$} and the corresponding intervals
\textcolor{blue}{$\I_{k,\theta}$} and \textcolor{blue}{$\I_{k+1,R_\omega(\theta)}$}
are drawn in \textcolor{blue}{blue}.
The map \textcolor{red}{$f_{m,\theta}$}, the interval
\textcolor{red}{$\IW{\theta}$} and the point \textcolor{red}{$\gams{\aii}\left(R_\omega(\theta)\right)$}
are drawn in \textcolor{red}{red}.}\label{fig:Cas1.3}
\end{figure}
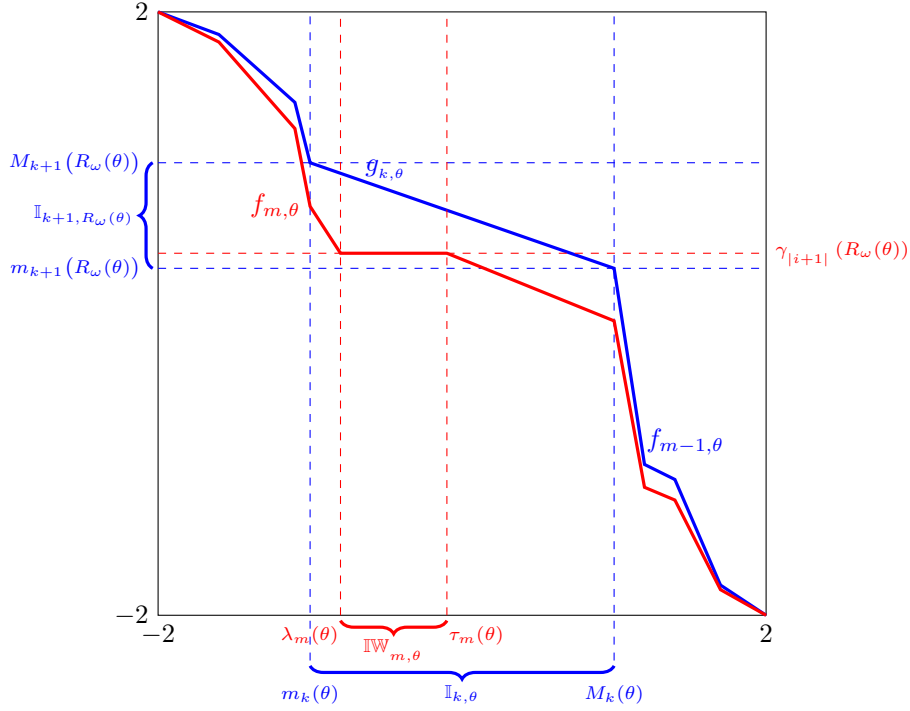
By Lemmas~\ref{QuePassaALesAles} and \ref{VerticalIntervalsIntheWings}(b)
and Definition~\ref{seqTmDefi},
\begin{align*}
\I_{i, \theta}  & = \{\gams{\ai}(\theta)\} = \{\gams{\ak}(\theta)\} \subset \IW{\theta},\text{ and}\\
f_{m,\theta}(x) & = \gams{\aii}\left(R_\omega(\theta)\right)\text{ for every $x \in \IW{\theta}$.}
\end{align*}

On the other hand,
by Definition~\ref{curvesinthewings} and Lemma~\ref{DepthintheWings}(a,b),
$\theta \in \WIB \subset \WDB,$ and
\[
\theta \in \basintabs{\ell} \subset
    \wobasint{i} \setminus \basintabs{i} \subset
    \obasintabs{k} \setminus \kstarset
\]
with $\ell = \bt[\led{\theta}]{\theta} \in \WDS$ and $\all > \ai > \ak.$
Then, by Lemma~\ref{Propertiesvarphi}(g) and Definition~\ref{curvesinthewings},
$\basicbox{\ell} \subset \Int\left(\basicbox{k} \setminus \setsilift{\kstar}\right)$
and
\[
  \IW{\theta} = \I_{\ell, \theta} \subset \I_{k, \theta}.
\]
Moreover, since $\theta \in \obasintabs{k} \subset \IBD,$
Definition~\ref{seqTmDefi},
Lemmas~\ref{gpositiva}(b) and \ref{gnegativa}(b),
and the definition of the maps $g_{_{i, \theta}}$ for $i \ge 0$
(Definition~\ref{defi-gi-positiva}) give
\begin{align*}
f_{m-1,\theta}\left(\IW{\theta}\right)
   &\subset f_{m-1,\theta}\left(\I_{k, \theta}\right)\\
   &\subset \begin{cases}
      \I_{k+1, R_\omega(\theta)} & \text{if $k < 0$ or $k \ge 0$ and $\theta \in \OBG{k}{k+1},$}\\
      \{\gams{k+1}\left(R_\omega(\theta)\right)\} & \text{if $k \ge 0$ and $\theta \in \basint{k} \setminus \OBG{k}{k+1}$.}
\end{cases}
\end{align*}

Now, as before, we will use Lemma~\ref{normainterna} to bound
$\norm{f_{m,\theta} - f_{m-1,\theta}}.$
We start with the simplest case:
$k \ge 0$ and $\theta \in \basint{k} \setminus \OBG{k}{k+1}.$
By Lemma~\ref{pointnormbound},
\begin{align*}
\norm{f_{m,\theta} - f_{m-1,\theta}}
   &= \norm{f_{m,\theta}\evalat{\IW{\theta}} - f_{m-1,\theta}\evalat{\IW{\theta}}} \\
   &= \abs{\gams{\aii}\left(R_\omega(\theta)\right) - \gams{\akk}\left(R_\omega(\theta)\right)}
    \le 2^{-\abs{\bt[m-1]{\theta}}}.
\end{align*}

Now we assume that $k < 0$ or $k \ge 0$ and $\theta \in \OBG{k}{k+1}.$
In this case Lemma~\ref{intervalsnormbound} applies.
By Lemmas~\ref{intervalsnormbound}, \ref{Propertiesvarphi}(d) and
Definition~\ref{PCgenerators}(\tsfR.2) and
Remark~\ref{PCgeneratorsExplicitConsequences}(\tsfR.2),
and Lemma~\ref{normainterna} we have
\begin{align*}
\gams{\aii}\left(R_\omega(\theta)\right)
   & \in \I_{i+1, R_\omega(\theta)} \subset \I_{k+1, R_\omega(\theta)},\\
f_{m-1,\theta}(x) &\in \I_{k+1, R_\omega(\theta)}\quad\text{for every $x \in \IW{\theta}$.}
\end{align*}
and
\begin{align*}
\norm{f_{m,\theta} - f_{m-1,\theta}}
   &= \LSleftlimits{\sup}{x \in \IW{\theta}} \abs{f_{m,\theta}(x) - f_{m-1,\theta}(x)} \\
   &= \LSleftlimits{\sup}{x \in \IW{\theta}} \abs{\gams{\aii}\left(R_\omega(\theta)\right) - f_{m-1,\theta}(x)}\\
   &\le 2\cdot 2^{-\ak} = 2\cdot 2^{-\abs{\bt[m-1]{\theta}}}.
\end{align*}

This ends the proof of the proposition in this case.

\begin{autocase}{2}
$\theta \in \obasintabs{i} \subset \wobasint{k} \setminus \basintabs{k}$
\upshape (see Figure~\ref{fig:Cas2} for a symbolic representation of this case).
\end{autocase}
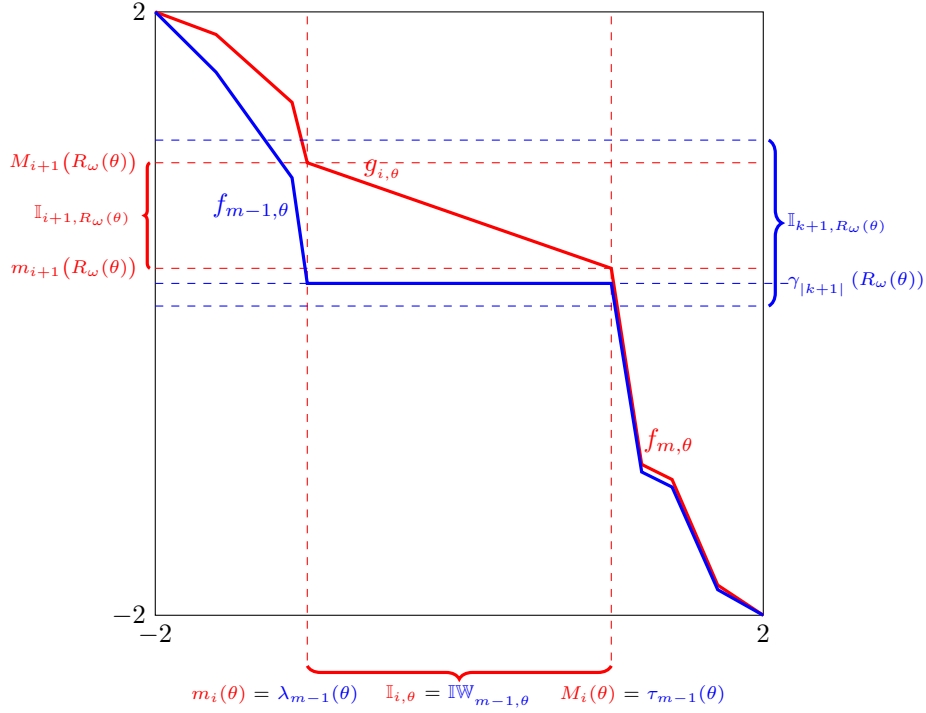
\begin{figure}
\begin{tikzpicture}[scale=2]
\draw (-2,-2) rectangle (2,2);
\foreach \c in {-2, 2} { \node[below] at (\c,-2) {$\c$}; \node[left] at (-2,\c) {$\c$}; }

\draw[dashed, color=red] (-1,-2.3) -- (-1,2); \draw[dashed, color=red] (1,-2.3) -- (1,2);
\draw[decorate, very thick, decoration={brace,amplitude=5pt, mirror, raise=2pt}, color=red] (-1,-2.3) -- (1,-2.3);
\node[below] at (0,-2.4) {\scriptsize\textcolor{red}{$\I_{i,\theta}$} = \textcolor{blue}{$\IW[m-1]{\theta}$}};
\node[below left, color=red] at (-0.6,-2.4) {\scriptsize\textcolor{red}{$m_i(\theta)$} = \textcolor{blue}{$\lambda_{m-1}(\theta)$}};
\node[below right, color=red] at (0.6,-2.4) {\scriptsize\textcolor{red}{$M_i(\theta)$} = \textcolor{blue}{$\tau_{m-1}(\theta)$}};

\draw[dashed, color=red] (-2, 1) -- (2, 1); \draw[dashed, color=red] (-2, 0.3) -- (2, 0.3);
\draw[decorate, very thick, decoration={brace,amplitude=2pt,raise=2pt}, color=red]  (-2, 0.3) -- (-2, 1);
\node[left, color=red] at (-2.1,0.65) {\scriptsize$\I_{i+1,R_\omega(\theta)}$};
\node[left, color=red] at (-2.05, 0.3) {\scriptsize$m_{i+1}\bigl(R_\omega(\theta)\bigr)$};
\node[left, color=red] at (-2.05, 1) {\scriptsize$M_{i+1}\bigl(R_\omega(\theta)\bigr)$};

\draw[very thick, color=red] (-2,2) -- (-1.6, 1.85) -- (-1.1, 1.4) --(-1, 1)
                 -- (1, 0.3) -- (1.2, -1) -- (1.4, -1.1) -- (1.7, -1.8) -- (2,-2);
\node[color=red, above right] at (1.15,-1) {$f_{m,\theta}$};
\node[color=red] at (-0.5,0.95) {$g_{_{i, \theta}}$};

\draw[dashed, color=blue] (-2, 0.2) -- (2.2, 0.2);
\node[right, color=blue] at (2.1, 0.2) {\scriptsize$\gams{\akk}\left(R_\omega(\theta)\right)$};
\draw[dashed, color=blue] (-2, 0.05) -- (2, 0.05); \draw[dashed, color=blue] (-2, 1.15) -- (2, 1.15);
\draw[decorate, very thick, decoration={brace, mirror, amplitude=5pt,raise=2pt}, color=blue]  (2, 0.05) -- (2, 1.15);
\node[right, color=blue] at (2.1,0.6) {\scriptsize$\I_{k+1,R_\omega(\theta)}$};

\draw[very thick, color=blue] (-2,2) -- (-1.6, 1.6) -- (-1.1, 0.9) --(-1, 0.2)
                 -- (1, 0.2) -- (1.2, -1.05) -- (1.4, -1.15) -- (1.7, -1.8303) -- (2,-2);
\node[color=blue, left] at (-1.05,0.711) {$f_{m-1,\theta}$};
\end{tikzpicture}
\caption{A symbolic representation of the maps
$f_{m,\theta}$ and $f_{m-1,\theta}$ in Case~2
($\theta \in \protect\obasintabs{i} \subset \protect\wobasint{k} \setminus \protect\basintabs{k}$)
of Proposition~\ref{distTmTm-1}.
The map \textcolor{blue}{$f_{m-1,\theta}$} and the corresponding intervals
\textcolor{blue}{$\IW[m-1]{\theta}$} and \textcolor{blue}{$\I_{k+1,R_\omega(\theta)}$}
are drawn in \textcolor{blue}{blue}.
The map \textcolor{red}{$f_{m,\theta}$} and the corresponding intervals
\textcolor{red}{$\I_{i,\theta}$} = \textcolor{blue}{$\IW[m-1]{\theta}$} and \textcolor{red}{$\I_{i+1,R_\omega(\theta)}$}
are drawn in \textcolor{red}{red}.}\label{fig:Cas2}
\end{figure}
In this case we will use Lemma~\ref{normainterna} with $\I_{i, \theta}.$
Thus, we need to compare the maps $f_{m,\theta}\evalat{\I_{i, \theta}}$
and $f_{m-1,\theta}\evalat{\I_{i, \theta}}.$

Directly from the definitions we get $k < 0,$
$\basintabs{i} \subset \IBD$ and
$\basintabs{k} \subset \IBD[m-1].$
Consequently, by Lemma~\ref{Dsets}(b) and Definition~\ref{curvesinthewings},
\[
  \theta \in \IBD
  \andq
  \theta \in \wIBD[m-1] \setminus \IBD[m-1] \subset \WDB[m-1] \subset \WB[m-1].
\]
Moreover, $\led[m-1]{\theta} = m,$
$i = \bt{\theta} = \bt[\led[m-1]{\theta}]{\theta} \in \WDS[m-1]$
and, by Definition~\ref{curvesinthewings},
$\theta \in \WIB[m-1],$ and
\[
  \IW[m-1]{\theta} = \I_{i, \theta}.
\]
Furthermore, since $k < 0,$ as in the proof of Lemma~\ref{pointnormbound},
$
R_\omega(\theta) \in \obasintabs{k+1}.
$
Thus, Definition~\ref{seqTmDefi},
Lemma~\ref{Propertiesvarphi}(d) and
Definition~\ref{PCgenerators}(\tsfR.2) and
Remark~\ref{PCgeneratorsExplicitConsequences}(\tsfR.2),
give
\[
 f_{m-1,\theta}(x) = \gams{\akk}\left(R_\omega(\theta)\right) \in \I_{k+1, R_\omega(\theta)}
\]
for every $x \in \I_{i, \theta} = \IW[m-1]{\theta}.$

Now we will use Lemma~\ref{normainterna}
to bound the norm $\norm{f_{m,\theta} - f_{m-1,\theta}}.$
By Definition~\ref{curvesinthewings} and Lemma~\ref{normainterna},
$\theta \in \IBD \subset \wIBD \setminus \WIB,$ and
\begin{align*}
\norm{f_{m,\theta} - f_{m-1,\theta}}
   &= \LSleftlimits{\sup}{x \in \I_{i, \theta}} \abs{f_{m,\theta}(x) - f_{m-1,\theta}(x)} \\
   &= \LSleftlimits{\sup}{x \in \I_{i, \theta}} \abs{f_{m,\theta}(x) - \gams{\akk}\left(R_\omega(\theta)\right)}.
\end{align*}

Next we will compute $f_{m,\theta}(\I_{i, \theta}).$
We start with the simplest case:
$i \ge 0$ and $\theta \in \obasint{i} \setminus \OBG{i}{i+1}.$
By Definition~\ref{seqTmDefi},
the definition of the maps $g_{_{i, \theta}}$ for $i \ge 0$
(Definition~\ref{defi-gi-positiva}) and Lemma~\ref{pointnormbound},
\begin{align*}
\norm{f_{m,\theta} - f_{m-1,\theta}}
   &= \LSleftlimits{\sup}{x \in \I_{i, \theta}} \abs{f_{m,\theta}(x) - \gams{\akk}\left(R_\omega(\theta)\right)}\\
   &= \abs{\gams{i+1}\left(R_\omega(\theta)\right) - \gams{\akk}\left(R_\omega(\theta)\right)}
    \le 2^{-\abs{\bt[m-1]{\theta}}}.
\end{align*}

Assume that $i < 0$ or $i \ge 0$ and $\theta \in \OBG{i}{i+1}.$
Then, again by Definition~\ref{seqTmDefi} and
Lemmas~\ref{gpositiva}(b), \ref{gnegativa}(b) and \ref{intervalsnormbound},
\[
f_{m,\theta}(x) \in \I_{i+1, R_\omega(\theta)} \subset \I_{k+1, R_\omega(\theta)}
\andq[for every]
x \in \I_{i, \theta},
\]
and
\begin{align*}
\norm{f_{m,\theta} - f_{m-1,\theta}}
   &= \LSleftlimits{\sup}{x \in \I_{i, \theta}} \abs{f_{m,\theta}(x) - \gams{\akk}\left(R_\omega(\theta)\right)}\\
   &\le 2\cdot 2^{-\ak} = 2\cdot 2^{-\abs{\bt[m-1]{\theta}}}.
\end{align*}

This ends the proof of the proposition in Case~2.

\begin{autocase}{3}
$\theta \in \obasintabs{i} \subset \obasintabs{k}.$
\end{autocase}
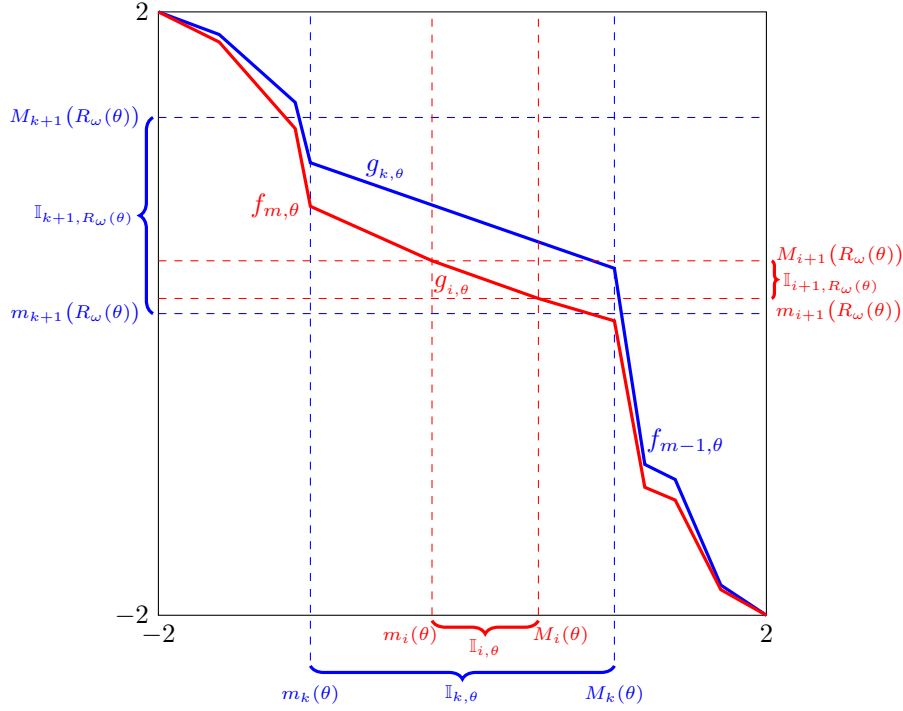
\begin{figure}
\begin{tikzpicture}[scale=2]
\draw (-2,-2) rectangle (2,2);
\foreach \c in {-2, 2} { \node[below] at (\c,-2) {$\c$}; \node[left] at (-2,\c) {$\c$}; }

\draw[dashed, color=blue] (-1,-2.3) -- (-1,2); \draw[dashed, color=blue] (1,-2.3) -- (1,2);
\draw[decorate, very thick, decoration={brace,amplitude=5pt, mirror, raise=2pt}, color=blue] (-1,-2.3) -- (1,-2.3);
\node[below, color=blue] at (0,-2.4) {\scriptsize$\I_{k,\theta}$};
\node[below, color=blue] at (-1,-2.4) {\scriptsize$m_k(\theta)$};
\node[below, color=blue] at (1,-2.4) {\scriptsize$M_k(\theta)$};

\draw[dashed, color=blue] (-2, 0) -- (2, 0); \draw[dashed, color=blue] (-2, 1.3) -- (2, 1.3);
\draw[decorate, very thick, decoration={brace,amplitude=5pt,raise=2pt}, color=blue]  (-2, 0) -- (-2, 1.3);
\node[left, color=blue] at (-2.1,0.65) {\scriptsize$\I_{k+1,R_\omega(\theta)}$};
\node[left, color=blue] at (-2.05, 0) {\scriptsize$m_{k+1}\bigl(R_\omega(\theta)\bigr)$};
\node[left, color=blue] at (-2.05, 1.3) {\scriptsize$M_{k+1}\bigl(R_\omega(\theta)\bigr)$};

\draw[very thick, color=blue] (-2,2) -- (-1.6, 1.85) -- (-1.1, 1.4) --(-1, 1)
                 -- (1, 0.3) -- (1.2, -1) -- (1.4, -1.1) -- (1.7, -1.8) -- (2,-2);
\node[color=blue, above right] at (1.15,-1) {$f_{m-1,\theta}$};
\node[color=blue] at (-0.5,0.95) {$g_{_{k, \theta}}$};

\draw[dashed, color=red] (-0.2,-2) -- (-0.2,2); \draw[dashed, color=red] (0.5,-2) -- (0.5,2);
\draw[decorate, very thick, decoration={brace,amplitude=5pt, mirror, raise=2pt}, color=red] (-0.2,-2) -- (0.5,-2);
\node[below, color=red] at (0.15,-2.1) {\scriptsize$\I_{i,\theta}$};
\node[below left, color=red] at (-0.1,-2.01) {\scriptsize$m_i(\theta)$};
\node[below right, color=red] at (0.4,-2.01) {\scriptsize$M_i(\theta)$};

\draw[dashed, color=red] (-2, 0.1) -- (2, 0.1); \draw[dashed, color=red] (-2, 0.35) -- (2, 0.35);
\draw[decorate, very thick, decoration={brace, mirror, amplitude=2pt,raise=2pt}, color=red]  (2, 0.1) -- (2, 0.35);
\node[right, color=red] at (2.05,0.2) {\scriptsize$\I_{i+1,R_\omega(\theta)}$};
\node[below right, color=red] at (2, 0.15) {\scriptsize$m_{i+1}\bigl(R_\omega(\theta)\bigr)$};
\node[above right, color=red] at (2, 0.25) {\scriptsize$M_{i+1}\bigl(R_\omega(\theta)\bigr)$};

\draw[very thick, color=red] (-2,2) -- (-1.6, 1.8) -- (-1.1, 1.2266) --(-1, 0.711)
                 -- (-0.2, 0.35) -- (0.5, 0.1) -- (1, -0.04848)
                 -- (1.2, -1.1515) -- (1.4, -1.2363) -- (1.7, -1.8303) -- (2,-2);
\node[color=red, left] at (-1,0.711) {$f_{m,\theta}$};
\node[color=red, above right] at (-0.25,0.05) {$g_{_{i, \theta}}$};
\end{tikzpicture}
\caption{A symbolic representation of the maps
$f_{m,\theta}$ and $f_{m-1,\theta}$ in
Subcase~3.1 from the proof of Proposition~\ref{distTmTm-1}
($\theta \in \protect\obasintabs{i}$ and $\I_{i, \theta} \subset \I_{k, \theta}$
 and either $k < 0$ or $k\ge 0$ and $\istar \in \BSG{k}{k+1}$).
The map \textcolor{blue}{$f_{m-1,\theta}$} and the corresponding intervals
\textcolor{blue}{$\I_{k,\theta}$} and \textcolor{blue}{$\I_{k+1,R_\omega(\theta)}$}
are drawn in \textcolor{blue}{blue}.
The map \textcolor{red}{$f_{m,\theta}$} and the corresponding intervals
\textcolor{red}{$\I_{i,\theta}$} and \textcolor{red}{$\I_{i+1,R_\omega(\theta)}$}
are drawn in \textcolor{red}{red}.}\label{fig:Cas3.1}
\end{figure}
In this case we have $\obasintabs{i} \subset \IBD$ and
$\obasintabs{k} \subset \IBD[m-1]$
so that, $\theta \in \IBD \cap \IBD[m-1].$
Moreover, by Lemma~\ref{Propertiesvarphi}(g),
$\basicbox{i} \subset \Int\left(\basicbox{k} \setminus \setsilift{\kstar}\right)$
and, hence,
\[
  \I_{i, \theta} \subset \I_{k, \theta}.
\]

Since $\theta \in \IBD,$ by Definition~\ref{curvesinthewings} and
Lemma~\ref{normainterna},
$\theta \in \wIBD \setminus \WIB,$ and
\[
\norm{f_{m,\theta} - f_{m-1,\theta}}
   = \norm{f_{m,\theta}\evalat{\I_{i, \theta}} - f_{m-1,\theta}\evalat{\I_{i, \theta}}}
   = \LSleftlimits{\sup}{x \in \I_{i, \theta}} \abs{f_{m,\theta}(x) - f_{m-1,\theta}(x)}.
\]
Thus, we need to compare the maps $f_{m,\theta}\evalat{\I_{i, \theta}}$
and $f_{m-1,\theta}\evalat{\I_{i, \theta}}.$
To do this we consider two subcases.

\begin{autocase}[Subcase]{3.1} Either $k < 0$ or $k\ge 0$ and $\theta \in \OBG{k}{k+1}$\\
\upshape (see Figure~\ref{fig:Cas3.1} for a symbolic representation of this case).
\end{autocase}
In this situation we aim at proving that
\[
 f_{m-1,\theta}\left(\I_{i, \theta}\right),
 f_{m,\theta}\left(\I_{i, \theta}\right) \subset
     \I_{k+1, R_\omega(\theta)}.
\]

We start with $f_{m-1,\theta}\left(\I_{i, \theta}\right).$
By Definition~\ref{seqTmDefi} and
Lemmas~\ref{gpositiva}(b) and \ref{gnegativa}(b)
we obtain
\[
f_{m-1,\theta}\left(\I_{i, \theta}\right)
 \subset f_{m-1,\theta}\left(\I_{k, \theta}\right) =
 g_{_{k, \theta}}\left(\I_{k, \theta}\right) \subset \I_{k+1, R_\omega(\theta)}.
\]

Next we show that
$f_{m,\theta}\left(\I_{i, \theta}\right) \subset \I_{k+1, R_\omega(\theta)}.$

Since $k < 0$ or $k\ge 0$ and $\theta \in \OBG{k}{k+1},$
by Definition~\ref{PCgenerators}(\tsfR.1) we obtain
\begin{equation}\label{eq:RotTheta}
R_\omega(\theta) \in
 \begin{cases}
      R_\omega\left(\obasintabs{k}\right) = \OBG{k+1}{\ak} \subset \obasintabs{k+1}\ \text{if $k < 0,$}\\
      R_\omega\left(\OBG{k}{k+1}\right) = \obasint{k+1}\ \text{if $k\ge 0$ and $\theta \in \OBG{k}{k+1}$.}
 \end{cases}
\end{equation}

Assume that $i < 0$ or $i \ge 0$ and $\theta \in \OBG{i}{i+1}.$
By \eqref{eq:RotTheta} with $k$ replaced by $i$,
\[
 R_\omega(\theta) \in \obasintabs{i+1} \cap \obasint{k+1}
                 \subset \wbasint{i+1} \cap \wbasint{k+1}.
\]
Therefore, since $\akk \le \aii$ and $k+1 \ne i+1,$
from Lemma~\ref{Propertiesvarphi}(g) we obtain $\akk < \aii,$
\begin{align*}
\basintabs{i+1} &\subset \obasintabs{k+1} \setminus \sstarset{k+1},\text{ and}\\
\basicbox{i+1}  &\subset \Int\left(\basicbox{k+1} \setminus \setsilift{\sstar{k+1}}\right).
\end{align*}
Thus, by Definition~\ref{seqTmDefi} and
Lemmas~\ref{gpositiva}(b) and \ref{gnegativa}(b),
\[
f_{m,\theta}\left(\I_{i, \theta}\right)
   = g_{_{i, \theta}}\left(\I_{i, \theta}\right)
   \subset \I_{i+1, R_\omega(\theta)} \subset \I_{k+1, R_\omega(\theta)}.
\]

Now we will consider the case
$i \ge 0$ and $\theta \in \obasint{i} \setminus \OBG{i}{i+1}.$
The fact that $\ak < \ai=i$ implies $\akk \le \ak+1 \le i.$
We claim that
\[
  \OBG{i+1}{i} \subset \obasintabs{k+1} \setminus \sstarset{k+1}.
\]
To prove the claim note that, by \eqref{eq:RotTheta},
\[
 R_\omega(\theta) \in R_\omega\left(\obasint{i}\right) \cap \obasintabs{k+1}
                  \subset  \OBG{i+1}{i} \cap \wbasint{k+1}.
\]
Moreover, the interval $\OBG{i+1}{i}$
is disjoint from $\wbasint{i}$ and $\wbasint{-i}$
by Definition~\ref{PCgenerators}(\tsfR.2).
Thus, $i \ne k+1, -(k+1)$ and, hence, $\akk < i$
(that is, $k+1 \in Z_{i-1}$).
So, there exists $q\in Z_{i-1}$ such that
$\BSG{i+1}{i} \cap \wbasint{q} \ne \emptyset$ and
$ \aq \ge \akk$ is maximal verifying these conditions.
By Definition~\ref{PCgenerators}(\tsfR.4),
\[
  \OBG{i+1}{i} \subset \wobasint{q} \setminus
     \left(\Bd\left(\basintabs{q}\right) \cup \qstarset\right).
\]
So, the claim holds when $q = k+1.$
Assume that $q \ne k+1.$ Then,
\[
R_\omega(\theta) \in \OBG{i+1}{i} \cap \obasintabs{k+1} \subset \wobasint{q} \cap \obasintabs{k+1}.
\]
Hence, by Lemma~\ref{Propertiesvarphi}(g), $ \aq > \akk$ and
\[
  \OBG{i+1}{i} \subset \wbasint{q} \subset \obasintabs{k+1} \setminus \sstarset{k+1}.
\]
This ends the proof of the claim.

On the other hand,
by Definition~\ref{PCgenerators}(\tsfR.2) and
Remark~\ref{PCgeneratorsExplicitConsequences}(\tsfR.2),
\[
  \left(\BSG{i+1}{i} \setminus \obasint{i+1}\right) \cap Z_{i+1} = \emptyset.
\]
Thus, by the claim,
\begin{align*}
R_\omega(\theta)
    &\in  R_\omega\left(\obasint{i} \setminus \OBG{i}{i+1}\right)
      = \OBG{i+1}{i} \setminus \obasint{i+1}\\
    &\subset \obasintabs{k+1} \setminus Z_{i+1}.
\end{align*}
By Definition~\ref{seqTmDefi},
the definition of the maps $g_{_{i, \theta}}$ for $i \ge 0$
(Definition~\ref{defi-gi-positiva})
and Lemma~\ref{Propertiesvarphi}(d) (with $\ell = k+1$ and $n = i+1$),
\[
   f_{m,\theta}\left(\I_{i, \theta}\right) =
   g_{_{i, \theta}}\left(\I_{i, \theta}\right) =
   \left\{\gams{i+1}\left(R_\omega(\theta)\right)\right\}
   \subset \I_{k+1, R_\omega(\theta)}.
\]

Summarizing, we have proved that
\[
 f_{m-1,\theta}\left(\I_{i, \theta}\right),
 f_{m,\theta}\left(\I_{i, \theta}\right) \subset
     \I_{k+1, R_\omega(\theta)}.
\]
So, by Lemma~\ref{Propertiesvarphi}(f) (and the fact that $\akk \ge \ak -1$),
\begin{align*}
\norm{f_{m,\theta} - f_{m-1,\theta}}
    &= \LSleftlimits{\sup}{x \in \I_{i, \theta}} \abs{f_{m,\theta}(x) - f_{m-1,\theta}(x)}
       \le \diam\left(\I_{k+1, R_\omega(\theta)}\right)\\
    &\le \diam\left(\basicbox{k+1}\right) \le 2^{-\akk}
     \le 2 \cdot 2^{-\ak} =  2 \cdot 2^{-\abs{\bt[m-1]{\theta}}}.
\end{align*}
This ends the proof of the proposition in this subcase.

\begin{autocase}[Subcase]{3.2} $k \ge 0$ and $\theta \in \obasint{k} \setminus \OBG{k}{k+1}.$ \end{autocase}
We start by computing $f_{m-1,\theta}\left(\I_{i, \theta}\right).$
By Definition~\ref{seqTmDefi} and
the definition of the maps $g_{_{k, \theta}}$ for $k \ge 0$
(Definition~\ref{defi-gi-positiva}),
\[
   f_{m-1,\theta}\left(\I_{i, \theta}\right)
       \subset f_{m-1,\theta}\left(\I_{k, \theta}\right)
       =  g_{_{k, \theta}}\left(\I_{k, \theta}\right)
       = \{\gams{k+1}\left(R_\omega(\theta)\right)\}.
\]

Analogously,
if $i \ge 0$ and $\theta \in \obasint{i} \setminus \OBG{i}{i+1},$
\[
   f_{m,\theta}\left(\I_{i, \theta}\right)
       =  g_{_{i, \theta}}\left(\I_{i, \theta}\right)
       = \{\gams{i+1}\left(R_\omega(\theta)\right)\}.
\]
Then, by Lemma~\ref{pointnormbound},
\begin{align*}
\norm{f_{m,\theta} - f_{m-1,\theta}}
&= \norm{f_{m,\theta}\evalat{\I_{i, \theta}} - f_{m-1,\theta}\evalat{\I_{i, \theta}}}\\
&= \abs{\gams{i+1}\left(R_\omega(\theta)\right) - \gams{k+1}\left(R_\omega(\theta)\right)}
   \le 2^{-\abs{\bt[m-1]{\theta}}}.
\end{align*}

Assume now that $i < 0$ or $i \ge 0$ and $\theta \in \OBG{i}{i+1}.$
By \eqref{eq:RotTheta},
Definition~\ref{seqTmDefi} and Lemmas~\ref{gpositiva}(b) and
\ref{gnegativa}(b)
\begin{align*}
R_\omega(\theta) & \in \obasintabs{i+1},\text{ and}\\
f_{m,\theta}\left(\I_{i, \theta}\right)
   & = g_{_{i, \theta}}\left(\I_{i, \theta}\right)
     \subset \I_{i+1, R_\omega(\theta)}.
\end{align*}
Moreover, if $k+1 < \aii,$
by Lemmas~\ref{QuePassaALesAles}(a) and \ref{Propertiesvarphi}(c),
we have
\[
f_{m-1,\theta}\left(\I_{i, \theta}\right)
    = \left\{\gams{k+1}\left(R_\omega(\theta)\right)\right\}
    = \left\{\gams{\aii-1}\left(R_\omega(\theta)\right)\right\}
    \subset \I_{i+1, R_\omega(\theta)}.
\]
Therefore, by Lemma~\ref{Propertiesvarphi}(f),
\begin{align*}
\norm{f_{m,\theta} - f_{m-1,\theta}}
   &= \LSleftlimits{\sup}{x \in \I_{i, \theta}} \abs{f_{m,\theta}(x) - f_{m-1,\theta}(x)}\\
   &= \LSleftlimits{\sup}{x \in \I_{i, \theta}} \abs{f_{m,\theta}(x) - \gams{\aii-1}\left(R_\omega(\theta)\right)}\\
   &\le \diam\left(\I_{i+1, R_\omega(\theta)}\right)
    \le \diam\left(\basicbox{i+1}\right) \le 2^{-\aii}\\
   &< 2^{-(k+1)} < 2^{-\abs{\bt[m-1]{\theta}}}.
\end{align*}

So, to end the proof of the proposition we have to show that,
in this subcase, $k+1 < \aii.$
To prove this, notice that when $i \ge 0,$
$k + 1 = \ak +1 < \ai + 1 = \aii.$
So, assume by way of contradiction that
$i < 0$ and $k+1 = \aii$ (recall that $k+1 \le \aii$).
Then, $k+1 = - (i+1)$ and, hence,
\begin{align*}
R_\omega(\theta)
    &\in  R_\omega\left(\obasint{k}\right) = \OBG{k+1}{k}, \text{ and}\\
R_\omega(\theta)
    &\in \obasintabs{i+1} = \OBG{-(k+1)}{k+1} \subset \wobasint{-(k+1)},
\end{align*}
which is a contradiction by  Definition~\ref{PCgenerators}(\tsfR.2).
\end{proof}

\def\cprime{$'$}

\end{document}